\documentclass[]{interact}

\usepackage{mathrsfs}
\usepackage{makecell,multirow,diagbox}
\usepackage{amsmath}
\usepackage{amssymb}
\usepackage{amsthm}
\usepackage{color}
\usepackage{xcolor}
\usepackage{graphicx}
\usepackage[caption=false]{subfig}
\usepackage[inline]{enumitem}
\usepackage{array}

\usepackage{makecell}
\usepackage{booktabs}
\usepackage{arydshln}
\usepackage{bm}
\usepackage{upgreek}
\usepackage{comment}

\usepackage[ruled]{algorithm2e}

\usepackage[colorlinks = true,linkcolor = blue,urlcolor  = blue,citecolor = blue, anchorcolor = blue]{hyperref}
\usepackage{lineno}

\usepackage{graphicx}
\newcommand{\tabstyle}{%
  \scriptsize
  \setlength{\tabcolsep}{3pt}
  \renewcommand{\arraystretch}{1.05}
}

\usepackage[numbers,sort&compress]{natbib}
\bibpunct[, ]{[}{]}{,}{n}{,}{,}

\theoremstyle{plain}

\numberwithin{equation}{section}

\newtheorem{theorem}{Theorem}[section]
\newtheorem{proposition}{Proposition}[section]
\newtheorem{lemma}{Lemma}[section]

\newtheorem{remark}{Remark}[section]

\newtheorem{definition}{Definition}[section]

\usepackage{lineno}

\begin{document}

\articletype{ORIGINAL ARTICLE}

\title{Reinforcement Learning for Dividend Optimization in Partially Observed Regime-Switching Diffusion Model}

\author{
\name{Zhongqin Gao\textsuperscript{a}, Yan Lv\textsuperscript{a} and Jingmin He\textsuperscript{b}
\thanks{CONTACT Email: Zhongqin Gao, zhongqingaox@126.com; Yan Lv, lvyan1998@aliyun.com; Jingmin He, nkjmhe\_2002@163.com.}}
\affil{\textsuperscript{a}School of Mathematics and Statistics, Nanjing University of Science and Technology, Nanjing 210094, China; \textsuperscript{b}School of Science, Tianjin University of Technology, Tianjin 300384, China}
}

\maketitle

\begin{abstract}
This paper studies the optimal dividend problem with a bounded payout rate in a partially observed regime-switching diffusion model, where, in practice, the market regime is unobserved and key model parameters are unknown. 
To address this partial-information setting, we propose a continuous-time reinforcement learning (RL) approach within an exploratory (entropy-regularized) stochastic control framework for discounted dividends under regime switching.
The associated exploratory Hamilton-Jacobi-Bellman (HJB) system admits semi-analytical characterizations of the value function and the optimal exploratory dividend policy, determined by two unknown functions solving two ordinary differential equations (ODEs) together with positive real roots of the induced quadratic equations.
Exploiting this structure, we introduce parametric families for both the value function and the policy, using low-degree polynomial approximations to the ODE solutions.
We then develop an actor-critic RL algorithm to learn the optimal exploratory policy through interactions with the market environment: it performs belief-state filtering from observed data and iterates policy evaluation and policy improvement online to refine the policy.
Numerical experiments demonstrate strong out-of-sample performance of the learned dividend policies.
\end{abstract}

\begin{keywords}
Optimal dividend problem, Reinforcement Learning, Regime-Switching, Actor-critic, Partial information, Exploratory
\end{keywords}

\begin{amscode}
93E20, 91B28, 93E11
\end{amscode}

\normalsize
\section{Introduction}

Dividend optimization is a critical focus in actuarial finance, aimed at maximizing the expected discounted present value of dividends paid to shareholders up to ruin. Initiated by \cite{De1957} in a binomial framework, the problem has since been studied extensively.
In the drifted Brownian model, \cite{Gerber2004} established the optimality of a barrier strategy. Many diffusion-based extensions (see, among others, \cite{Asmussen2000}, \cite{Paulsen2003}, \cite{Decamps2007} and \cite{gao2024}) typically lead to barrier- or band-type optimal policies and their variants. For a comprehensive survey, we refer to \cite{Avanzi2009}.
Motivated by regulatory capital requirements, liquidity constraints or internal payout rules, many studies impose an upper bound on the dividend rate, under which the optimal policy typically takes a threshold form. For instance, \cite{asmussen1997} established threshold optimality for a drifted Brownian model with bounded dividend rates.

Most of the dividend-optimization literature assumes that the model (environment) parameters are known. In this setting, researchers typically adopted the classical control approach: they modeled the surplus as a controlled Markov process and, by the dynamic programming principle, characterized the value function as the solution to a Hamilton-Jacobi-Bellman (HJB) equation (or variational inequality) with appropriate boundary conditions, which in turn yields a time-homogeneous optimal feedback control (often of barrier, band or threshold type). However, in practice the model parameters are not directly observable. Implementing such a policy requires estimating them from historical data and substituting the estimates into the control law. This plug-in procedure can be fragile: parameters may be difficult to estimate accurately, and optimal controls can be extremely sensitive to estimation error, so that even small errors may lead to materially different or potentially poor decisions, as emphasized in the mean-variance (MV) problem by \cite{luenberger2014}.
Thus, while the classical approach delivers a theoretically optimal policy, its plug-in implementation with estimated parameters is generally no longer optimal and may even be significantly inferior.

In recent years, reinforcement learning (RL) has become increasingly popular in financial applications, see \cite{jaimungal2022} and \cite{hambly2023} for surveys.
Rooted in machine learning, RL focuses on learning an optimal policy through trial-and-error interaction with an unknown environment. In particular, a stochastic policy is used to generate different trial actions, the learner implements a sampled action in the environment and receives feedback in the form of states and rewards to evaluate the current policy. This policy is then updated based on this feedback to improve a prescribed objective. Repeating this interaction-update loop progressively refines the policy and ultimately drives convergence toward optimality.
In this sense, RL decision-making is governed by the exploration-exploitation trade-off: exploitation follows a policy that is optimal relative to the learner's current value estimates based on the information collected so far;
whereas exploration deliberately samples alternative actions through trial and error to improve the policy toward optimality. As new information arrives, RL updates and adapts the policy accordingly, thereby improving overall performance.
Compared with the classical approach, RL offers two potential advantages: first, stochastic exploration may overcome the traps of local optima during learning, as argued in \cite{zhou2021}; second, continual learning from observed market data can mitigate sensitivity to parameter-estimation error, as evidenced by \cite{wang20202} and \cite{wang2023}.

To extend RL to continuous time and continuous state spaces under general diffusion dynamics, \cite{wang2020} proposed an entropy-regularized exploratory stochastic control framework, in which the control is allowed to be stochastic and the exploration weight is tuned by an entropy coefficient (often termed the temperature parameter). Building on this perspective, \cite{jia20222, jia2022, jia20223} developed a systematic actor-critic theory for policy evaluation (PE), policy improvement (PI) and policy-gradient (PG) computation.
More broadly, continuous-time RL algorithms have been developed for a range of financial applications, including MV portfolio selection (see e.g., \cite{wang20202}, \cite{jia2022}, \cite{huang2022}, \cite{wu2024}), mean-field games (see e.g., \cite{guo2022} and \cite{firoozia2022}), and optimal execution (see \cite{wang2023}).
In very recent work, \cite{bai2023} presented the first examination of the optimal dividend problem in a continuous-time, entropy-regularized RL framework, for a drifted Brownian model with a bounded payout rate, assuming that the payout cap exceeds the temperature parameter (thereby ensuring an increasing value function), and constructing an approximating sequence for the optimal policy via PI and PE.
Subsequently, \cite{hu2024} removed this ordering restriction and instead classified the value function into three cases according to the relative magnitude of these two quantities.
Across these studies, the state variables (e.g., the surplus) are typically assumed fully observable, and the resulting RL algorithms are model-free in the sense that value functions and policies are learned directly from trajectory data without explicit parameter estimation.

Moreover, regime-switching models, introduced by \cite{hamilton1989new}, capture market cyclicality (e.g., bull/bear states) by modeling the surplus as a Markov-modulated process with regime-dependent coefficients (e.g., drift/volatility). Incorporating such dynamics can yield more robust control policies.
In dividend optimization, the canonical regime-switching diffusion model has been studied extensively; for example, \cite{sotomayor2011} and \cite{zhu2014} analyzed the bounded-rate case and derived semi-analytical characterizations of the optimal dividend policy via the classical approach, while \cite{jiang2012} investigated the unconstrained-rate counterpart.
A common feature of these works is the assumption that the regime is directly observed.
In practice, however, insurers typically observe only surplus outcomes, and the regime is unobserved and must be inferred from the observed data, leading to a substantially more challenging partial-information control problem.
A natural idea is to apply filtering to augment the state by the regime belief and thereby recast the problem as a complete-information control problem, in the spirit of the separation principle in \cite{xiong2007}.
Concretely, the task decomposes into filtering and optimization. One first constructs the belief (posterior) process for the underlying regime from the observed data, and when model parameters are unknown, this step entails parameter estimation (now including the regime transition intensities), which is an essential difference from model-free RL algorithms. Conditional on the belief state, the optimization step reduces to a complete-information control problem, to which continuous-time RL techniques such as the exploratory control framework of \cite{wang2020} and the actor-critic PI/PE theory of \cite{jia20222, jia2022, jia20223} become directly applicable.

Applying RL to dividend optimization remains an active and evolving area, especially under regime switching and partial information. In this paper, we study an exploratory optimal dividend problem for a regime-switching diffusion model with unobserved market regimes, unknown environment parameters, and a bounded dividend rate, extending \cite{bai2023} and \cite{hu2024} to the regime-switching setting.
To enable RL algorithm design, we adopt the entropy-regularized exploratory control framework of \cite{wang2020}, and derive semi-analytical representations of the value function and the optimal exploratory dividend policy, characterized by two unknown functions solving two ordinary differential equations (ODEs) along with the positive real roots of the induced quadratic equations.
We establish ODE-based structural properties of these functions and approximate them with low-degree polynomials.
The resulting optimal policy is a truncated exponential (Gibbs) distribution, whose parameters depend on the surplus level and the temperature, which motivates restricting the learned policy class to Gibbs distributions.
Building on this structure, we then introduce parametric families for both the value function and the policy and develop an actor-critic RL algorithm.
For each episode, the algorithm performs filtering at every interaction step to construct the regime belief, which serves as the observable (filtered) state. Conditional on this belief, it alternates PE and PI to learn the optimal policy, with PE implemented via two martingale-based updates: Martingale-Loss (ML) minimization and episodic online CTD.
Additionally, ruin is approximated on a discrete time grid by a first-hitting rule over a sufficiently long finite horizon: a path is absorbed at zero upon the first nonpositive surplus and is otherwise terminated at the horizon. Empirically, our algorithm is largely insensitive to this discretized ruin-time approximation.

The rest of the paper is organized as follows. Section \ref{sec2} revisits the complete-information regime-switching optimal dividend problem with a bounded dividend rate, and constructs the candidate value function and optimal policy.
Section \ref{sec3} introduces belief-state filtering to recast the partial-information problem as a complete-information control problem, and formulates it within an entropy-regularized exploratory control framework for continuous-time RL.
Section \ref{sec4} solves the associated exploratory HJB system, characterizes the value function and optimal policy and establishes their key analytical properties.
Section \ref{sec5} develops an actor-critic RL algorithm that exploits the semi-analytical structure from Section~\ref{sec4} to solve the partial-information dividend problem in Section~\ref{sec3}.
Numerical analyses are provided in Section \ref{sec6}, showing that the learned RL-based dividend policies perform well out of sample, the mean value estimate closely matches the finite-difference (FD) benchmark while the cross-path variance is substantially reduced.

\section{Regime-switching optimal dividend problem under completely observable}\label{sec2}
In this section, under complete information, we revisit the optimal dividend problem for a regime-switching diffusion model with a bounded dividend rate.
Let $(\Omega, \mathcal{F}, (\mathcal{F}_t^{X})_{t\geq 0}, \mathbb{P})$ be a filtered probability space upon which a stochastic process $X\equiv (X_t)_{t\geq 0}$ is defined and adapted to the filtration $\mathbb{F}^X:=(\mathcal{F}_t^{X})_{t\geq 0}$, where the filtration is augmented in the usual way so as to satisfy the usual conditions (i.e., completeness and right-continuity). For $\mathbf{a}> 0$, denote $([0, \mathbf{a}], \mathscr{B}([0, \mathbf{a}]), \mathscr{P}([0, \mathbf{a}]))$ as the filtered probability space generated by $[0, \mathbf{a}]$, with $\mathscr{B}([0, \mathbf{a}])$ the Borel $\sigma$-algebra and $\mathscr{P}([0, \mathbf{a}])$ the space of all probability measure.
Given a metric space $\mathbb{X}$ with the Borel sets $\mathscr{B}(\mathbb{X})$, we define $\mathbb{L}^0(\mathbb{X})$ as the collection of all $\mathscr{B}(\mathbb{X})$-measurable functions, and for any $p \geq 1$, $\mathbb{L}^p(\mathbb{X})$ denotes the space of $p$-th integrable functions on $\mathbb{X}$, and $\mathbb{L}^{\infty}(\mathbb{X})$ the space of functions that are essentially bounded ($\|\cdot\|_{\infty} := \|\cdot\|_{\mathbb{L}^{\infty}(\mathbb{X})}$). Moreover, for a domain $\mathbb{D} \subset \mathbb{R}$, let $\mathbb{C}^k(\mathbb{D})$ denote the space of $k$-times continuously differentiable functions, where $\mathbb{C}(\mathbb{D}) = \mathbb{C}^0(\mathbb{D})$. The subspace $\mathbb{C}_b^k(\mathbb{D})\subset \mathbb{C}^k(\mathbb{D})$ consists of functions whose values and all derivatives up to order $k$ are bounded. For $\mathbb{D}_1, \mathbb{D}_2 \subset \mathbb{R}$, $\mathbb{C}^{2}(\mathbb{D}_1 \times \mathbb{D}_2)$ is the set of functions on $\mathbb{D}_1 \times \mathbb{D}_2$ that are twice continuously differentiable with respect to (w.r.t.) each variable, and whose mixed partial derivatives up to order two are also continuous. Denote by $\bar{\mathbb{D}}$ the closure of $\mathbb{D}$. In general, $\mathbb{R}_+=[0, \infty)$, and $(x)_+:=\max\{x, 0\}$ for $x\in\mathbb{R}$.
The law and the expectation are denoted as $\mathbb{P}$ and $\mathbb{E}$, respectively.
\footnote{\tiny{Let $f:\mathbb{D}_1 \times \mathbb{D}_2 \to \mathbb{D}_3$ with $f=f(x,y)$. For clarity, the subscripts $f_x$ and $f_{xx}$ denote the first- and second-order partial derivatives of $f$ w.r.t. the first argument, respectively; while $f_y$ and $f_{yy}$ denote those w.r.t. the second argument; second-order mixed partial derivatives are written as $f_{xy}$, with the dependence on $(x,y)$ omitted unless otherwise specified.}}

\subsection{Problem formulation}
Let $W=(W_t)_{t\geq 0}$ be a one-dimensional standard Brownian motion, and $I= (I_t)_{t \geq 0}$ a continuous-time Markov chain with finite state space $E$ and generator $Q=(q_{ij})_{i,j\in{E}}$, independent of $W$, where transition rates $q_{ij}> 0$ for $i\neq j$ and $\sum_{j\in E}q_{ij}=0$ for each $i\in E$.
Let $I$ be the economic regime process, and $X=(X_t)_{t\geq 0}$ the insurer's surplus process, modeled as a regime-switching drifted Brownian motion satisfying
\begin{equation}\label{sde1}
    \mathrm{d}X_t := \mu_{I_t} \mathrm{d}t + \sigma \mathrm{d}W_t,
\end{equation}
with initial capital $X_0=x\in\mathbb{R}_+$ and initial regime $I_0=i\in E$, where the regime-dependent drift $\mu_{I_t}>0$ represents the premium rate, and $\sigma > 0$ is a constant volatility. We assume complete information: both $X$ and $I$ are observable, and $\mu_{I_t}$ $(I_t\in E)$ and $\sigma$ are known constants, with $\mu_i\neq \mu_j$ whenever $i,j\in E$ and $i\neq j$.
We restrict attention to the two-state case $E=\{1, 2\}$, corresponding to the standard bull-bear interpretation.
Throughout, the environment parameters are specified as $\mu_i$, $\sigma$ and $q_{ij}$ for $i,j\in E$.

Under a threshold dividend strategy, \eqref{sde1} becomes the controlled surplus process $X^{\alpha}= (X_t^{\alpha})_{t\geq 0}$,
\begin{align}\label{sde2}
	\mathrm{d} X_t^{\alpha}
:=\big(\mu_{I_t} - \alpha_t\big)\mathrm{d} t + \sigma \mathrm{d} W_t,
\end{align}
with $X_0^{\alpha}=x\in\mathbb{R}_+$, where $\alpha_t$ is the dividend rate at time $t\geq 0$. A dividend control $\alpha= (\alpha_t)_{t\geq 0}$ is admissible if it is $\mathbb{F}$-adapted and satisfies $\alpha_t\in [0, \mathbf{a}]$. We denote by $\mathscr{A}_{x,i}[0, \mathbf{a}]$ the set of all admissible dividend policies given the initial state $(x, i)$.
Define the ruin time by $T^{\alpha}_{x,i}:=\inf\{t\geq 0: X_t^{\alpha}\leq0\}$ with the usual convention $\inf \emptyset = \infty$. Clearly, $X_{T^{\alpha}_{x,i}}^{\alpha}=0$, the company is ruined at time $T_{x,i}^{\alpha}$, no dividends are paid thereafter.

The expected total discounted dividends paid until ruin under $\alpha$ is given by
\begin{eqnarray*}
	J(x,i; \alpha):=\mathbb{E}_{x,i}\Big[\int_0^{T^{\alpha}_{x,i}} e^{-\int_0^{t} \delta_{I_s} \mathrm{d}s} \alpha_t \mathrm{d}t\Big],
\end{eqnarray*}
where $\delta: E\rightarrow(0,\infty)$ is the regime-dependent discount rate (with $\delta_i$ constant for each $i\in E$), and $\mathbb{E}_{x,i}[\cdot]:=\mathbb{E}[\cdot | X_0^{\alpha}=x, I_0 = i]$. The regime-switching optimal dividend problem is to maximize $J$ over all admissible policies $\alpha\in \mathscr{A}_{x,i}[0, \mathbf{a}]$, leading to
\begin{align}\label{def V}
	&V(x,i):=\sup_{\alpha\in \mathscr{A}_{x,i}[0, \mathbf{a}]} J(x, i; \alpha)
=\sup_{\alpha\in \mathscr{A}_{x,i}[0, \mathbf{a}]} \mathbb{E}_{x,i}\Big[\int_0^{T^{\alpha}_{x,i}} e^{-\int_0^{t} \delta_{I_s}\mathrm{d}s} \alpha_t \mathrm{d}t \Big], \ \ i\in E,
\end{align}
termed the value function, with $V(0,i)=0$. A dividend policy $\alpha^*\in \mathscr{A}_{x,i}[0, \mathbf{a}]$ is optimal if $V(x,i)=J(x,i; \alpha^*)$. It also satisfies $0\leq V(x,i) \leq \tfrac{\mathbf{a}}{\delta_1 \wedge \delta_2}$.

\subsection{HJB equation and verification theorem}
Following Theorem~3.1 of \cite{sotomayor2011}, we state the HJB equation and verification theorem for the problem \eqref{def V} as follows.
\begin{theorem}[Verification theorem]\label{verthem00}
Suppose that $\nu: \mathbb{R}_+ \times E \rightarrow \mathbb{R}_+$ satisfies $\nu(\cdot, i)\in \mathbb{C}^{2}((0, \infty)\backslash \mathcal{P}_i) \cap \mathbb{C}^{1}((0, \infty))$ for a finite subset $\mathcal{P}_i \in \mathbb{R}_+$ and $i\in E$.
If $\nu(x, i)$ satisfies the following HJB equation
\begin{align}\label{hjbeq1122}
&\delta_i \nu(x, i)  = \sup_{\alpha \in [0, \mathbf{a}]} \Big(\alpha \big(1-\nu_x(x, i)\big)\Big) + \tfrac{\sigma^2}{2} \nu_{xx}(x, i)+ \mu_i \nu_x(x, i) + \sum_{j\neq i} q_{ij}\big(\nu(x, j) - \nu(x,i)\big),
\end{align}
with boundary condition $\nu(0,i)=0$, and $|\nu|$, $|\nu_x|$ are bounded,
then $\nu(x,i)$ is the value function of the regime-switching optimal dividend problem \eqref{def V}, the associated optimal dividend policy $\alpha^* =(\alpha_t^*)_{t\geq0}$ is induced by a deterministic feedback control $\alpha^*: \mathbb{R}_+ \times E \to [0, \mathbf{a}]$ (we use $\alpha^*$ for both the map and the induced control process), formally,  $\alpha_t^* = \alpha^*(X_t^{\alpha^*},I_t)$ with
\begin{eqnarray}\label{umax11}
		&\alpha^*(X_t^{\alpha^*},I_t)=\arg \sup_{\alpha_t\in[0, \mathbf{a}]}\Big(\alpha_t \big(1-\nu_x(X_t^{\alpha^*},I_t)\big)\Big)=\left \{
		\begin{array}{ll}
			0 & \mathrm{for}\ \nu_x(X_t^{\alpha^*}, I_t)>1 , \\
\mathbf{a} & \mathrm{for}\ \nu_x(X_t^{\alpha^*}, I_t)< 1,
		\end{array} \right.
\end{eqnarray}
and when $\nu_x(X_t^{\alpha^*},I_t)=1$, $\alpha_t^*$-insensitive in the sense that any $\alpha_t^*\in[0, \mathbf{a}]$ is optimal.
\end{theorem}

Let $\mathcal{B}_1, \mathcal{B}_2 \subset \mathbb{R}_+$ be such that $\nu_x(x,i) \geq 1$ $(\alpha^*=0)$ for $x \in \mathcal{B}_1$, and $\nu_x(x,i) \leq 1$ $(\alpha^* = \mathbf{a})$ for $x \in \mathcal{B}_2$, with $b_i: = \mathcal{B}_1 \cap \mathcal{B}_2$ denoting the threshold (or switching) point where the behavior of $\nu_x$ changes. By \eqref{hjbeq1122} and \eqref{umax11}, we have, for $i=1,2$ and $i\neq j$,
\begin{align}
& \tfrac{\sigma^2}{2} \nu_{xx}(x, i)+ \mu_i \nu_x(x, i) - (\delta_i -q_{ii}) \nu(x, i) + q_{ij} \nu(x, j) = 0, \ \ \ \ \ \ \ \ \ \ \ \ \ \mathrm{for} \ x\in \mathcal{B}_1, \label{valuehg1}\\
&\tfrac{\sigma^2}{2} \nu_{xx}(x, i)+ (\mu_i - \mathbf{a}) \nu_x(x, i) - (\delta_i - q_{ii}) \nu(x, i) + q_{ij} \nu(x, j) + \mathbf{a} = 0, \mathrm{for} \ x\in \mathcal{B}_2, \label{valuehg2}
\end{align}
with $\lim_{x\rightarrow \infty} \nu(x, i)<\infty$ and the boundary conditions
\begin{align}\label{barr45}
&\nu(0,i)=0, \ \nu(b_i -, i) = \nu(b_i +, i), \ \nu_x(b_i -, i) = \nu_x(b_i +, i)=1,  \nonumber\\
&\nu(b_i -, j) = \nu(b_i +, j), \ \nu_x(b_i -, j) = \nu_x(b_i +, j).
\end{align}
The ODE$_i^{(0)}$ in \eqref{valuehg1} characterizes the value dynamics in the absence of dividend payments, while the ODE$_i^{(\mathbf{a})}$ in \eqref{valuehg2} describes an inhomogeneous linear equation for $\nu(\cdot,i)$ in the dividend-paying region.

\subsection{Value function and optimal dividend policy}\label{sec23}
The value function and the associated optimal dividend policy for the problem \eqref{def V} are analyzed in Section 3 of \cite{sotomayor2011}. The analysis is based on three possible configurations of the boundary slope $\nu_x(0, i)$: $\nu_x(0,i)>1$ for both $i=1,2$; $\nu_x(0, i) >1$ and $\nu_x(0, j)\leq1$ for $i\in\{1, 2\}$ and $i\neq j$; and $\nu_x(0,i)\leq 1$ for both $i=1,2$. In each case, the candidate value function is constructed via \eqref{valuehg1}-\eqref{valuehg2} to satisfy the coupled HJB system, yielding a closed-form expression up to unknown coefficients and regime-dependent optimal dividend thresholds.
These unknown quantities are then determined by the smooth-fit conditions in \eqref{barr45}; while their explicit forms may be unwieldy and existence/uniqueness may not be established in general, they can be computed numerically in concrete instances.

Consequently, the optimal control is a regime-switching threshold policy $\alpha^{\mathbf{b}^*}:=(\alpha_t^{b_{I_t}^*})_{t\geq 0, I_t\in E}$ with $\mathbf{b}^*=(b_1^*, b_2^*)$ and $b_i^*\geq 0$, where the dividend rate $\alpha_t^{b_{I_t}^*}$ is characterized by the regime-dependent threshold $b_{I_t}^*$, formally, $\alpha_t^{b_{I_t}^*} = \mathbf{a} \cdot \mathbf{1}_{\{X_t^{\alpha^{\mathbf{b}^*}} \geq b_{I_t}^*, t< T_{x,i}^{\alpha^{\mathbf{b}^*}}\}}$.
It is further shown that $\nu_x(\cdot, i)$ is continuous, positive and bounded on $\mathbb{R}_+$, with $\nu_x(\infty,i)=0$ for $i= 1,2$, and that $\nu_x(\cdot, i)$ is non-increasing in $x$ so that $\nu(x, i)$ is concave. In particular, the value function $\nu(x, i)$ satisfies the boundary conditions
\begin{align}\label{defc3c4}
&\lim_{x\to \infty}\nu(x, 1) = \tfrac{(\delta_2 +q_{12} + q_{21}) \mathbf{a}}{(\delta_1 + q_{12}) (\delta_2 + q_{21}) - q_{12} q_{21}}, \ \ \
\lim_{x\to \infty}\nu(x, 2) = \tfrac{(\delta_1 +q_{12} + q_{21}) \mathbf{a}}{(\delta_1 + q_{12}) (\delta_2 + q_{21}) - q_{12} q_{21}}.
\end{align}

\section{Regime-switching optimal dividend problem under partially observable}\label{sec3}
In practice, neither the market regime nor the driving Brownian motion is observable, so regime-dependent primitives (e.g., the premium rate) are also not directly observed. The insurer observes only the surplus process $X$, and the available information is represented by the observation filtration $\mathbb{F}^{X}:=(\mathcal{F}_t^{X})_{t\geq 0}$. The dividend policy $\alpha$ is $\mathbb{F}^{X}$-adapted, whereas the regime process $I$ and the Brownian motion $W$ are not. This yields a partial-information dividend control problem.
To the best of our knowledge, this problem has not been addressed in the existing literature. The main challenge arises from the resulting variational formulation in terms of partial differential equations (PDEs), for which explicit solutions are generally unattainable.
In this section, following the separation principle of \cite{xiong2007}, we first separate the problem into filtering and optimization, which reduces to a complete-information control problem conditional on the filtered belief state, and then formulate it within an exploratory stochastic control framework for continuous-time RL.

\subsection{Filtering and innovation}
To estimate the underlying regime from the observed data, define the belief state (filtered probability) that the market is in regime 1 at time $t$, conditional on $\mathcal{F}_t^{X}$, by $p_t := \mathbb{P}(I_t = 1 | \mathcal{F}_t^X)$.
The filtered drift and discount rate are defined as, respectively,
\begin{align}
    & \widehat{\mu}_t := \mathbb{E}[\mu_{I_t} | \mathcal{F}_t^X] = (\mu_1 - \mu_2)p_t + \mu_2, \label{hatmu} \\
    & \widehat{\delta}_t := \mathbb{E}[\delta_{I_t} | \mathcal{F}_t^X] = (\delta_1 - \delta_2)p_t + \delta_2. \label{hatmu11}
\end{align}
Consider the innovation process $\widehat{W}\equiv (\widehat{W}_t)_{t\geq 0}$ given by
\begin{align}\label{hatwt}
    & \mathrm{d} \widehat{W}_t := \tfrac{\mathrm{d}X_t - \widehat{\mu}_t \mathrm{d} t}{\sigma}.
\end{align}
As established in \cite{kallianpur2013}, Theorem 8.1.2, $\widehat{W}\equiv (\widehat{W}_t)_{t\geq 0}$ is a one-dimensional standard Brownian motion w.r.t. $\mathbb{F}^X$.
By \eqref{hatmu} and \eqref{hatwt}, the surplus dynamics \eqref{sde2} can be reformulated so as to be adapted to $\mathbb{F}^X$ as follows
\begin{eqnarray}\label{sde3}
	\mathrm{d} X_t^{\alpha}:=(\widehat{\mu}_t - \alpha_t)\mathrm{d} t + \sigma \mathrm{d} \widehat{W}_t.
\end{eqnarray}
Moreover, the belief state process $(p_t)_{t\geq 0}$ is governed by the Wonham filter (see \cite{wonham1964}),
\begin{equation}\label{wondpt}
    \mathrm d p_t = \big(q_{21}(1 - p_t) - q_{12}p_t\big) \, \mathrm d t + \tfrac{\mu_1 - \mu_2}{\sigma}p_t(1 - p_t) \, \mathrm d\widehat{W}_t,
\end{equation}
which is $\mathbb{F}^X$-adapted and evolves as a diffusion on $[0, 1]$.
In \cite{dai2010} it is shown that both 0 and 1 are entrance boundaries, which means that $(p_t)_{t\geq 0}$ starting at these points enters the interval $(0,1)$, and once inside, it will not reach 0 or 1 within a finite time. Our analysis assumes $p_0 \in (0,1)$, thus confining the state space of $(p_t)_{t\geq 0}$ to $(0,1)$.
The drift term in \eqref{wondpt} reflects the prior transition dynamics and is equivalently $(q_{12}+q_{21})\big(\tfrac{q_{21}}{q_{12}+q_{21}}-p_t\big)$, which implies $p_t$ mean-reverts to $q_{21}/(q_{12}+q_{21})$, the stationary probability of being in regime 1, at rate $q_{12}+q_{21}$.
The diffusion term in \eqref{wondpt} captures the flow of new information available to the insurer and thus governs the learning rate, which increases with the model-implied signal-to-noise ratio $(\mu_1 - \mu_2)/\sigma$, attains its maximum at $p_t=1/2$ (maximum uncertainty) and vanishes as $p_t\to 0$ or $p_t\to 1$.

The state process $(X_t^{\alpha}, p_t)$ defined by \eqref{sde3}-\eqref{wondpt} is a two-dimensional controlled Markov diffusion on $\mathbb{R}_+\times (0,1)$. The partial-observation dividend problem is converted into an equivalent complete-information problem with $(X_t^{\alpha}, p_t)$ as the observed state.

\subsection{Exploratory problem formulation}\label{subsec32}
We now adopt the exploratory, entropy-regularized stochastic control framework for continuous-time RL of \cite{wang2020} and \cite{bai2023}. A key modeling feature is to randomize the control action $\alpha_t$ as a probability distribution over the control space $[0, \mathbf{a}]$, i.e. $\alpha_t\in \mathscr{P}([0, \mathbf{a}])$. This approach aligns with the theory of relaxed control, where the classical control approach is recovered as a special case corresponding to a Dirac (point-mass) measure.
Specifically, relaxed control actions are randomized to represent exploration (learning) via a stochastic policy, leading to a distributional (measure-valued) control process whose feedback form is characterized by a density function
$$\pi: (X_t^{\pi}, p_t, \cdot)\in \mathbb{R}_+ \times (0,1)\times [0,\mathbf{a}] \mapsto \pi(\cdot| X_t^{\pi}, p_t)\in \mathbb{L}^1_+([0, \mathbf{a}]),$$
with $X_t^{\pi}$ and $p_t$ denoting the observed surplus level and the filtered belief at time $t$, respectively.
In what follows, a distributional control (exploratory dividend policy) is identified with its density process $\pi= \big(\pi(\cdot| X_t^{\pi}, p_t)\big)_{t\geq 0}$.
The overall level of exploration induced by $\pi$ is captured by its differential entropy
\begin{align*}
&H(\pi(\cdot | X_t^{\pi}, p_t)):=-\int_0^{\mathbf{a}} \pi(u| X_t^{\pi}, p_t)\ln\pi(u| X_t^{\pi}, p_t)\mathrm{d}u,
\end{align*}
which satisfies $H(\pi(\cdot| X_t^{\pi}, p_t))\leq \ln \mathbf{a}$ for any control $\pi$ supported on $[0, \mathbf{a}]$, as implied by Jensen's inequality.

The filtered state model \eqref{sde3} is modified as the following state dynamics
\begin{eqnarray}\label{sde4}
	\mathrm{d} X_t^{\pi}
:=\big(\widehat{\mu}_{t} - \alpha^{\pi}(X_t^{\pi}, p_t)\big)\mathrm{d} t + \sigma \mathrm{d} \widehat{W}_t,
\end{eqnarray}
with $X_0^{\pi}=x$ and $p_0=p$, where $\alpha^{\pi}(X_t^{\pi}, p_t)$ is randomly sampled from $\pi(\cdot | X_t^{\pi}, p_t)$. Notably, $\alpha^{\pi}(X_t^{\pi}, p_t)$ is a function of the evolving state $(X_t^{\pi}, p_t)$ rather than of time $t$ itself.
A standard method for sampling from this distribution is inverse transform sampling: one first sampling $Z_t$ from the uniform distribution on $[0,1]$ (generate $Z_t\sim \mathrm{Unif}[0,1]$), and then inverting $Z_t$ by $\alpha^{\pi}(X_t^{\pi}, p_t) = (F^{\pi})^{-1}(Z_t)$, where $F^{\pi}$ is the cumulative distribution function of $\pi(\cdot| X_t^{\pi}, p_t)$. This construction ensures that the randomness caused by $\pi$ at time $t$ is entirely driven by the exogenous uniform variable $Z_t$.
Define the ruin time for the state process \eqref{sde4} as
$$T^{\pi}_{x,p}:=\inf\{t\geq 0: X_t^{\pi}\leq0\}.$$
We assume that $\pi$ is independent of the innovation process $\widehat{W}$, the random source in the filtered state model \eqref{sde3}.
To expand the original filtered probability space to support the new random source introduced by $\pi$, we adopt the mathematical framework of \cite{jia2022}, in which the filtration is enlarged to incorporate information from action randomization, as follows:
for any $t\in[0, T]$ with given $T\geq 0$, let
$\overline{\mathcal{F}}_t:= \mathcal{F}_t^X \vee \mathcal{Z}_t$, where $\mathcal{Z}_t=\sigma(Z_s:s\in[0,t])$ with $(Z_s)_{s\in[0,t]}$ denoting a process of mutually independent random variables uniformly distributed over $[0, 1]$ (see \cite{sun2006}). We define $\overline{\mathbb{P}}$ as an extension of $\mathbb{P}$ from $\mathcal{F}_T^X$ to $\overline{\mathcal{F}}_{T}$, and they coincide on $\mathcal{F}_T^X$, i.e. $\overline{\mathbb{P}}|_{\mathcal{F}_T^X}=\mathbb{P}$. Accordingly, the original filtered probability space
$(\Omega, \mathcal{F}, (\mathcal{F}_t^{X})_{t\in[0,T]}, \mathbb{P})$
is expanded to
$(\Omega, \mathcal{F}, (\overline{\mathcal{F}}_t)_{t\in[0,T]}, \overline{\mathbb{P}})$.

The entropy-regularized objective function under the exploratory stochastic control framework is defined as,
\begin{eqnarray}\label{obje12}
	J(x, p; \pi):=\overline{\mathbb{E}}_{x,p}\Big[\int_0^{T^{\pi}_{x,p}} e^{-\int_0^{t} \widehat{\delta}_s \mathrm{d}s} \mathcal{H}_{\lambda}^{\pi}(X_t^{\pi}, p_t)\mathrm{d}t\Big],
\end{eqnarray}
where $\overline{\mathbb{E}}$ denotes expectation under $\overline{\mathbb{P}}$, which accounts for both the environmental noise modeled by $\widehat{W}$ and the randomness induced by $\pi$, and
\begin{align*}
&\mathcal{H}_{\lambda}^{\pi}(x, p):=\alpha^{\pi}(x, p) +\lambda H(\pi(\cdot| x, p))
=\alpha^{\pi}(x, p) - \lambda\int_0^{\mathbf{a}} \pi(u| x, p) \ln\pi(u| x, p)\mathrm{d}u,
\end{align*}
with $\lambda>0$ being the exploration weight reflecting the trade-off between exploration and exploitation.

The exploratory optimal dividend problem is to maximize $J(x,p;\pi)$ defined in \eqref{obje12}. To this end, following \cite{wang2020}, we derive the exploratory state dynamics, whose trajectories describe the expected evolution of the state process \eqref{sde4} averaged over infinitely many realizations of randomized control actions.
More precisely, consider the discrete version of \eqref{wondpt} and \eqref{sde4}: for small $\Delta t>0$, $\Delta X_t^{\pi}:=X_{t+\Delta t}^{\pi} - X_t^{\pi}$ and $\Delta p_t:= p_{t+\Delta t} - p_t$. By a law-of-large-numbers argument as in \cite{wang2020}, we have
\begin{align*}
&\overline{\mathbb{E}}[\Delta X_t^{\pi}] = \Big(\widehat{\mu}- \int_0^{\mathbf{a}}u \pi(u)\mathrm{d} u\Big) \Delta t, \ \ \ \ \
\overline{\mathbb{E}}[(\Delta X_t^{\pi})^2] = \sigma^2 \Delta t + o(\Delta t),\\
&\overline{\mathrm{Cov}}(\Delta X_t^{\pi}, \Delta p_t) = (\mu_1 - \mu_2) p (1-p) \Delta t +o(\Delta t),
\end{align*}
where $\lim_{\Delta t\rightarrow 0} o(\Delta t)/\Delta t =0$. In the above expectations, the randomness in the actions is averaged out, and the covariance between $\Delta X_t^{\pi}$ and $\Delta p_t$ must also be taken into account due to the regime-switching structure.
As diffusion processes are essentially determined by the first two moments,
the exploratory version of the state dynamics \eqref{sde4} under $\pi$ takes the form of the following stochastic differential equation (SDE)
\begin{align}\label{sde5}
	\mathrm{d} \widehat{X}_t^{\pi}
:=& \Big(\widehat{\mu}_{t}
- \int_0^{\mathbf{a}} u \pi(u| \widehat{X}_t^{\pi}, p_t)\mathrm{d}u\Big)\mathrm{d} t + \sigma \mathrm{d} \widehat{W}_t,
\end{align}
with $\widehat{X}_0^{\pi}=x$ and $p_0=p$. We also denote
\begin{align}\label{defbpixp}
&b^{\pi}(x, p):=\int_0^{\mathbf{a}} u \pi(u | x, p) \mathrm{d} u.
\end{align}
Define the ruin time for the process \eqref{sde5} as $\widehat{T}^{\pi}_{x,p}:=\inf\{t\geq 0: \widehat{X}_t^{\pi}\leq0\}$. According to \cite{jia2022}, the state process $(X_t^{\pi})_{t\geq 0}$ defined in \eqref{sde4} is identical in law to the exploratory state process $(\widehat{X}_t^{\pi})_{t\geq 0}$ in \eqref{sde5}, formally, $\widehat{X}_t^{\pi} \overset{d}{=} X_t^{\pi}$ and $\mathbb{P}(\widehat{T}_{x,p}^{\pi}<t)=\mathbb{P}(T_{x,p}^{\pi}<t)$.
The entropy-regularized objective function \eqref{obje12} is equivalent to
\begin{eqnarray}\label{obje122}
	J(x, p; \pi):=\mathbb{E}_{x,p}\Big[\int_0^{\widehat{T}_{x,p}^{\pi}} e^{-\int_0^{t} \widehat{\delta}_s\mathrm{d}s} \widehat{\mathcal{H}}_{\lambda}^{\pi}(\widehat{X}_t^{\pi}, p_t)\mathrm{d}t\Big],
\end{eqnarray}
where the expectation is taken w.r.t. the randomness of $\widehat{W}$ only, and
\begin{align}\label{enth1}
&\widehat{\mathcal{H}}_{\lambda}^{\pi}(x, p)
:=\int_0^{\mathbf{a}} u \pi(u| x, p)\mathrm{d}u +\lambda H(\pi(\cdot| x, p))
=\int_0^{\mathbf{a}}\big(u - \lambda \ln\pi(u| x, p)\big) \pi(u| x, p)\mathrm{d}u.
\end{align}

The definition of admissible controls for a fixed $\lambda$ is given below, with the dependence of $\pi$ on $\lambda$ suppressed in the notation for clarity.
\begin{definition}\label{def1}
A distributional control $\pi=\pi(\cdot| \cdot, \cdot)\in \mathbb{L}^0(\mathbb{R}_+ \times (0,1) \times [0,\mathbf{a}])$ is called admissible, if for any $(\widehat{X}_t^{\pi}, p_t)\in \mathbb{R}_+ \times (0,1)$,
\begin{enumerate}[label=(\roman*)]
\item $\pi(\cdot| \widehat{X}_t^{\pi}, p_t)\in \mathbb{L}_+^1([0,\mathbf{a}])$ with $\int_0^{\mathbf{a}} \pi(u| \widehat{X}_t^{\pi}, p_t) \mathrm{d} u =1$, and the process $(\widehat{X}_t^{\pi}, p_t, \cdot)\in \mathbb{R}_+ \times (0,1) \times [0, \mathbf{a}] \mapsto \pi(\cdot | \widehat{X}_t^{\pi}, p_t)$ is progressively measurable;
\item The SDE \eqref{sde5} admits a unique strong solution $\widehat{X}^{\pi}$ for any initial state $(x, p)$;
\item $\mathbb{E}_{x,p}\Big[\int_0^{\widehat{T}_{x,p}^{\pi}} e^{-\int_0^{t} \widehat{\delta}_s \mathrm{d}s} \left| \widehat{\mathcal{H}}_{\lambda}^{\pi}(\widehat{X}_t^{\pi}, p_t) \right| \mathrm{d}t \Big] < \infty$.
\end{enumerate}
Denote by $\Pi_{x,p}[0, \mathbf{a}]$ the set of all admissible controls for given $x$ and $p$.
\end{definition}
The regime-switching exploratory optimal dividend problem amounts to maximizing $J(x,p; \pi)$ in \eqref{obje122},
\begin{align}\label{opdipro}
&V(x,p):=\sup_{\pi\in{\Pi_{x,p}[0, \mathbf{a}]}} J(x, p; \pi)
=\sup_{\pi\in{\Pi_{x,p}[0, \mathbf{a}]}} \mathbb{E}_{x,p}\Big[\int_0^{\widehat{T}_{x,p}^{\pi}} e^{-\int_0^{t} \widehat{\delta}_s \mathrm{d}s} \widehat{\mathcal{H}}_{\lambda}^{\pi}(\widehat{X}_t^{\pi}, p_t)\mathrm{d}t\Big],
\end{align}
the so-called value function, where $\widehat{\mathcal{H}}_{\lambda}^{\pi}(\widehat{X}_t^{\pi}, p_t)$ is given by \eqref{enth1}. A dividend policy $\pi^*\in \Pi_{x,p}[0, \mathbf{a}]$ is optimal if $V(x,p)=J(x, p; \pi^*)$ for all $(x, p)\in \mathbb{R_+}\times (0,1)$.
The following proposition establishes the monotonicity and boundedness of $V(x,p)$ for $\mathbf{a}>1$, its proof is deferred to Appendix \ref{prop21pf}.

\begin{proposition}\label{prop21}
For $\mathbf{a}>1$, the value function $V$ satisfies the following properties:
\begin{enumerate}[label=(\roman*)]
\item $V(x,p)> V(y,p)$, for any $x> y\geq 0$ and $p\in(0, 1)$;
\item $0\leq V(x,p) \leq \tfrac{1}{\delta_1 \wedge \delta_2}(\mathbf{a} + \lambda \ln{\mathbf{a}})$, for any $x\in \mathbb{R}_+$ and $p\in(0, 1)$.
\end{enumerate}
\end{proposition}

\section{Solve to the regime-switching exploratory optimal dividend problem}\label{sec4}
This section focuses on solving the regime-switching exploratory optimal dividend problem \eqref{opdipro}. By applying dynamic programming principle, we derive the associated exploratory HJB equation and obtain semi-analytical representations of the value function and the optimal exploratory dividend policy, characterized by two unknown functions solving two ODEs and two positive real roots of the corresponding quadratic equations.
Let $v: \mathbb{R}_+ \times (0,1) \rightarrow \mathbb{R}$ be a function in $\mathbb{C}^{2}((0, \infty)\times (0,1))$.

\subsection{Exploratory HJB equation and optimal exploratory dividend policy}
The dynamic programming principle is as follows,
\begin{align*}
&V(x, p)=\sup_{\pi\in{\Pi_{x,p}[0, \mathbf{a}]}} \mathbb{E}_{x,p}\Big[\int_0^{\zeta \wedge\widehat{T}^{\pi}_{x,p}} e^{-\Lambda_t} \widehat{\mathcal{H}}_{\lambda}^{\pi}(\widehat{X}_t^{\pi}, p_t)\mathrm{d}t + e^{-\Lambda_{\zeta \wedge\widehat{T}^{\pi}_{x,p}}} V(\widehat{X}_{\zeta\wedge\widehat{T}^{\pi}_{x,p}}^{\pi}, p_{\zeta\wedge\widehat{T}^{\pi}_{x,p}})\Big],
\end{align*}
with $\Lambda_t := \int_0^t \widehat{\delta}_s \mathrm{d}s$ and $\zeta>0$ an $\mathbb{F}$-stopping time, such as the first regime-switching time.
Following this standard arguments, the exploratory HJB equation is
\begin{align}\label{hjbeq}
&\big(\delta_2 + (\delta_1-\delta_2)p \big) v \nonumber \\
& = \sup_{\pi\in{\Pi_{x,p}[0, \mathbf{a}]}} \Big( \int_0^{\mathbf{a}}\big(u (1-v_x) - \lambda \ln \pi(u | x,p)\big) \pi(u | x, p) \mathrm{d} u \Big)
+ \tfrac{\sigma^2}{2}v_{xx} + \big(\mu_2 +(\mu_1 - \mu_2) p \big)v_x\nonumber\\
&\
+\tfrac{1}{2} \tfrac{(\mu_1-\mu_2)^2}{\sigma^2} p^2 (1-p)^2 v_{pp}+\big(q_{21}- (q_{12}+q_{21})p\big)v_p +(\mu_1-\mu_2) p (1-p)v_{xp},
\end{align}
with boundary condition $v(0,p)=0$ for all $p\in(0,1)$.
Combining the method of Lagrange multipliers with the calculus of variations (see \cite{ferguson2004}), subject to the constraint $\int_0^{\mathbf{a}} \pi(u | x, p)\mathrm{d} u =1$,
the optimal distributional control that maximizes the integral on the right-hand side of \eqref{hjbeq} is given by the Gibbs form
\begin{align}\label{densityop}
&\pi^*(u| x, p) = G\big(u, 1-v_x(x,p)\big),
\end{align}
where for $u\in[0, \mathbf{a}]$ and $y\in\mathbb{R}$,
\begin{align}\label{gu1y}
&G(u,1-y):= \tfrac{1-y}{\lambda(e^{\frac{\mathbf{a}}{\lambda}(1-y)} -1)} \cdot e^{\frac{u}{\lambda}(1-y)} \mathbf{1}_{\{y\neq 1\}} + \tfrac{1}{\mathbf{a}} \mathbf{1}_{\{y=1\}}.
\end{align}
To capture the state-dependent behavior, let $\vartheta(x,p):= \big(1- v_x(x,p)\big)/\lambda$, \eqref{densityop} is reformulated as a truncated exponential distribution over $[0, \mathbf{a}]$, $\pi^*(u| x,p) = \tfrac{\vartheta(x,p) e^{\vartheta(x,p) u}}{e^{\vartheta(x,p) \mathbf{a}} -1}$.
This structure shows that the exploratory policy is directly governed by the marginal value $v_x$: when $v_x>1$, $\pi^*(u| x,p)$ is strictly decreasing in $u$ so it assigns higher probability to small dividend rates near zero; when $v_x<1$, $\pi^*(u| x,p)$ is strictly increasing in $u$ so it has a large probability to take a large dividend rate close to $\mathbf{a}$, and at $v_x=1$, it degenerates to the uniform distribution on $[0, \mathbf{a}]$.
\eqref{densityop} may also be viewed as a probabilistic relaxation of the classical bang-bang policy $\alpha^*$ defined by \eqref{umax11}: it assigns higher probabilities to dividend rates near the classical optimal (i.e., 0 or $\mathbf{a}$), with the likelihood to take a certain rate decreasing as it moves away from the classical value.

An equivalent expression of \eqref{densityop} is given by
$\pi^*(u| \widehat{X}_t^{\pi^*}, p_t) = G\big(u, 1-v_x(\widehat{X}_t^{\pi^*}, p_t)\big)$. Substituting it into the exploratory state dynamics \eqref{sde5} yields, for $v_x = v_x(\widehat{X}_t^{\pi^*}, p_t)$,
\begin{align}\label{sde56}
	\mathrm{d} \widehat{X}_t^{\pi^*}
&=\Big(\widehat{\mu}_{t}
- \big( \mathbf{a} - \tfrac{\lambda}{1- v_x}
+ \tfrac{\mathbf{a}}{e^{\frac{\mathbf{a}}{\lambda}(1-v_x)} -1}\big)
\mathbf{1}_{\{v_x\neq 1\}}
-\tfrac{\mathbf{a}}{2} \mathbf{1}_{\{v_x= 1\}} \Big)
\mathrm{d} t + \sigma \mathrm{d} \widehat{W}_t.
\end{align}
SDE \eqref{sde56} features bounded drift and constant volatility, so it admits a unique strong solution, where the boundedness of the drift follows from the fact that $\int_0^{\mathbf{a}} u \pi^*(u| \widehat{X}_t^{\pi^*}, p_t) \mathrm{d} u \in[0, \mathbf{a}]$ and $p_t\in (0,1)$.

\subsection{Verification theorem}
Substituting \eqref{densityop} into exploratory HJB equation \eqref{hjbeq} yields the following second-order nonlinear PDE:
\begin{align}\label{hjb156}
& \tfrac{\sigma^2}{2}v_{xx} + f_{\lambda}(v_x) + \big(\mu_2 +(\mu_1 - \mu_2)p\big) v_x
+\tfrac{1}{2} \tfrac{(\mu_1-\mu_2)^2}{\sigma^2} p^2 (1-p)^2 v_{pp}+\big(q_{21}- (q_{12}+q_{21})p\big)v_p \nonumber\\
& +(\mu_1-\mu_2) p (1-p)v_{xp}
 - \big(\delta_2 + (\delta_1-\delta_2)p \big) v =0,
\end{align}
with $v(0, p)=0$ for all $p\in(0,1)$, where
\begin{align}\label{funcfyp}
& f_{\lambda}(y):=\big( \lambda \ln \tfrac{\lambda (e^{\frac{\mathbf{a}}{\lambda}(1-y)}-1)}{1-y} \big) \mathbf{1}_{\{y\neq 1\}}
+ ( \lambda\ln \mathbf{a} ) \mathbf{1}_{\{y=1\}}.
\end{align}
Particularly, $f_{\lambda}(0) = \lambda \ln\big(\lambda(e^{\frac{\mathbf{a}}{\lambda}} -1)\big)$. By integration by parts, we obtain, for $z=\tfrac{1-y}{\lambda}$,
\begin{align}\label{firdefy}
&f_{\lambda}'(y)
= \tfrac{-\int_0^{\mathbf{a}} r e^{\frac{1-y}{\lambda} r}\mathrm{d} r}{\int_0^{\mathbf{a}} e^{\frac{1-y}{\lambda} r}\mathrm{d} r}
= \tfrac{-(\tfrac{\mathbf{a}}{z} e^{z \mathbf{a}} -\tfrac{e^{z \mathbf{a}} -1 }{z^2})}{\tfrac{e^{z \mathrm{a}} -1}{z}}
=- (\tfrac{\mathbf{a} z -1}{z} +\tfrac{\mathbf{a}}{e^{z \mathbf{a}} -1})
=- (\tfrac{\mathbf{a} (1-y) - \lambda}{1-y} +\tfrac{\mathbf{a}}{e^{\frac{1-y}{\lambda} \mathbf{a}} -1}).
\end{align}
The properties of $f_{\lambda}$ and $f_{\lambda}'$ are summarized below, see Appendix \ref{propfuncf0011} for the proof.

\begin{proposition}\label{propfuncf00}
The function $f_{\lambda}$ defined by \eqref{funcfyp} can be equivalently expressed as, for all $y\in\mathbb{R}$, $f_{\lambda}(y)= \lambda \ln \big(\mathbf{a}+\sum_{n=2}^{\infty} \tfrac{\mathbf{a}^n (1-y)^{n-1}}{n! \lambda^{n-1}}\big)$, and belongs to $\mathbb{C}^{\infty}(\mathbb{R})$.
It satisfies:
\begin{enumerate}[label=(\roman*)]
\item $f_{\lambda}(y)$ is convex and strictly decreasing. Its derivative $f_{\lambda}'(y)$ is continuous, bounded and strictly increasing, with $f_{\lambda}'(y)\in(-\mathbf{a}, 0)$, and $\lim_{y\to 1} f_{\lambda}'(y)= -\mathbf{a}/{2}$;
\item $f_{\lambda}(y)$ is uniformly Lipschitz continuity on $\mathbb{R}$ (hence uniformly continuous), with Lipschitz constant at most $\mathbf{a}>0$, and it is locally bounded on $\mathbb{R}$;
\item $f_{\lambda}(0)>0$ if and only if $\mathbf{a}>\lambda \ln (1+1/\lambda)$;
\item There exists a unique $y_0\in\mathbb{R}$ such that $f_{\lambda}(y_0) =0$, and $y_0\in (1,1+\lambda)$ if $\mathbf{a}>1$;
\item $|f_{\lambda}(y)|< \mathbf{a}|y| + c$ for some constant $c\in\mathbb{R}$, which depends on $\mathbf{a}$ only;
\item $f_\lambda(y)$ converges pointwise to $f_0(y):=\mathbf{a}(1-y)_+$ as $\lambda\downarrow 0$.
\end{enumerate}
\end{proposition}

A verification theorem below identifies the value function and the optimal distributional control for problem \eqref{opdipro}, with the proof deferred to Appendix \ref{verpf}.

\begin{theorem}[Verification theorem]\label{verthem}
Suppose that $v(x, p)\in \mathbb{C}^2((0, \infty)\times(0,1))$. If $v(x, p)$ satisfies the exploratory HJB equation \eqref{hjb156} with boundary condition $v(0, \cdot)=0$, and $|v|$, $|v_x|$, $|v_p|$ are bounded, then $v(x,p)$ is the value function of the regime-switching exploratory optimal dividend problem \eqref{opdipro}, and $\pi^*$ defined by \eqref{densityop} constitutes the optimal distributional control.
\end{theorem}
We analyze the asymptotics of $v(x, p)$ as $x\rightarrow \infty$, see Appendix \ref{vbound22} for the proof.
\begin{proposition}\label{vbound2}
Suppose $v(x,p)$ is the value function of problem \eqref{opdipro}, $\pi^*$ given by \eqref{densityop} is the optimal distributional control.
Letting $x\to \infty$ and $\delta=\delta_1=\delta_2$, we have
\begin{align}\label{asyv00}
\lim_{\substack{x\rightarrow\infty }} v(x,p)
=\tfrac{f_{\lambda}(0)}{\delta}.
\end{align}
\end{proposition}

\begin{remark}\label{remark41}
When the two-regime market collapses to a single regime (so that regime switching is absent), the optimal exploratory dividend policy $\pi^*$ in \eqref{densityop} reduces to that of the drifted Brownian motion studied in \cite{bai2023} and \cite{hu2024}.
\end{remark}

\subsection{Solution to exploratory HJB equation}
An analytical solution to the exploratory HJB equation \eqref{hjb156} is intractable because the coefficient of $v_{pp}$ is a fourth-order polynomial in $p$. We first characterize $v(x, p)$ by several key structural properties:
(i) the boundary condition $v(0, \cdot) = 0$;
(ii) it is bounded and satisfies the asymptotic condition \eqref{asyv00} as $x\rightarrow \infty$ when $\delta_1 =\delta_2$;
(iii) $\lim_{x\rightarrow \infty} v_x(x, p) =0$, reflecting diminishing marginal value in magnitude.
Motivated by (i)-(iii), we adopt a two-exponential saturation ansatz as a flexible approximation
\begin{equation}\label{eq:single_exp_ansatz}
v(x,p) = g_1(p) (1 - e^{-\kappa_1 x}) + g_2(p) (1 - e^{-\kappa_2 x}), \qquad \kappa_1, \kappa_2 > 0,
\end{equation}
which separates the $x$- and $p$-dependence and ensures $\lim_{x\to\infty} v(x,p) = g_1(p) + g_2(p):=g(p)$, where $g_i(p)$, $i=1,2$ are chosen to satisfy the boundary conditions
\begin{align}
&g(0) = g_1(0) + g_2(0)
= \tfrac{(\delta_1 +q_{12} +q_{21}) f_{\lambda}(0)}{(\delta_1+q_{12})(\delta_2 +q_{21})- q_{12}q_{21}}, \label{boucong1}\\
&g(1) = g_1(1) + g_2(1)
= \tfrac{(\delta_2 +q_{12} +q_{21}) f_{\lambda}(0)}{(\delta_1+q_{12})(\delta_2 +q_{21})- q_{12}q_{21}},
\label{boucong2}
\end{align}
consistent with \eqref{asyv00} and as $\lambda\downarrow 0$, with the completely observed limit \eqref{defc3c4}.
Initialize admissible endpoint splits $\varpi_i^{(j)}\in[0,1]$ $(i=1,2,j=0,1)$ with $\varpi_1^{(j)}+\varpi_2^{(j)}=1$, and impose the conditions
\begin{equation}\label{addicondi}
g_i(0)=\varpi_i^{(0)}\,g(0),\qquad
g_i(1)=\varpi_i^{(1)}\,g(1).
\end{equation}

Differentiating \eqref{eq:single_exp_ansatz} and substituting into the PDE \eqref{hjb156} yields
\begin{align*}
& -\tfrac{\sigma^2}{2} \big(\kappa_1^2 g_1(p) e^{-\kappa_1 x} +\kappa_2^2 g_2(p) e^{-\kappa_2 x}\big)
+ f\big(\kappa_1 g_1(p) e^{-\kappa_1 x} +\kappa_2 g_2(p) e^{-\kappa_2 x}\big)\\
&
+ D(p) \, \big(\kappa_1 g_1(p) e^{-\kappa_1 x} +\kappa_2 g_2(p) e^{-\kappa_2 x}\big)
+ \tfrac12 A(p) \big(g_1''(p) (1 - e^{-\kappa_1 x}) +g_2''(p) (1 - e^{-\kappa_2 x})\big) \\
&
+ B(p) \big( g_1'(p) (1 - e^{-\kappa_1 x}) +g_2'(p) (1 - e^{-\kappa_2 x})\big)
 \\
& + (\mu_1 - \mu_2) p (1-p) \big(\kappa_1 g_1'(p) e^{-\kappa_1 x}+\kappa_2 g_2'(p) e^{-\kappa_2 x}\big)\\
&
- C(p) \big(g_1(p) (1 - e^{-\kappa_1 x}) +g_2(p) (1 - e^{-\kappa_2 x})\big) = 0,
\end{align*}
where, for notational brevity,
\begin{align}\label{apbpcpdp}
& A(p) := \tfrac{(\mu_1 - \mu_2)^2}{\sigma^2} \, p^2 (1-p)^2, \quad
B(p) := q_{21} - (q_{12} + q_{21}) p, \nonumber\\
& C(p) := \delta_2 + (\delta_1 - \delta_2) p, \quad \quad \ \
D(p) := \mu_2 + (\mu_1 - \mu_2) p.
\end{align}
The coefficient functions $A(p), C(p), D(p)$ are positive on $(0,1)$, belong to $\mathbb{C}^2((0,1))$ with bounded first and second derivatives, and $A(p)$ is locally uniformly elliptic, i.e. for every $\varepsilon\in(0,\tfrac12)$, there exists a constant $A_{\varepsilon}>0$ such that $\inf_{p\in[\varepsilon,1-\varepsilon]} A(p) \ge A_{\varepsilon}>0$.
Since $v_x(x,p) \to 0$ as $x\to\infty$, it is natural to approximate $f$ by its first-order Taylor expansion around $0$, namely, $f_{\lambda}(y) = f_{\lambda}(0) + f_{\lambda}'(0) y + o(y)$, which yields
\begin{align*}
&f_{\lambda}(v_x) = f_{\lambda}\big(\kappa_1 g_1(p) e^{-\kappa_1 x}+\kappa_2 g_2(p) e^{-\kappa_2 x}\big) \approx f_{\lambda}(0) + f_{\lambda}'(0)\big(\kappa_1 g_1(p) e^{-\kappa_1 x} + \kappa_2 g_2(p) e^{-\kappa_2 x}\big).
\end{align*}
Collecting the constant terms (independent of $x$) and the coefficients of $e^{-\kappa_i x}$ separately, we arrive at two coupled pairs of ODEs, for $i=1, 2$,
\begin{align}
&\tfrac12 A(p) g_i''(p) + B(p) g_i'(p) - C(p)g_i(p) + \hat{\varpi}_i(p) f_{\lambda}(0) = 0,\label{eq:V_const_eq00}\\
&\tfrac12 A(p) g_i''(p) + \big(B(p) -(\mu_1 - \mu_2)p(1-p)\kappa_i \big) g_i'(p)\nonumber\\
 & \ \ \ \ \ \ \  - \big(-\tfrac{\sigma^2}{2} \kappa_i^2 + f_{\lambda}'(0) \kappa_i + C(p) + D(p) \kappa_i \big) g_i(p)
= 0,\label{eq:V_exp_eq}
\end{align}
where $\hat{\varpi}_1(p):=(1-p) \varpi_1^{(0)} + p \varpi_1^{(1)} \in(0, 1)$ and $\hat{\varpi}_2(p) := 1- \hat{\varpi}_1(p)$.

A closed-form solution to the elliptic ODE \eqref{eq:V_const_eq00} is not available, as the coefficient of $g_i''(p)$ is a quartic polynomial in $p$. Proposition \ref{theosoluodse} establishes the well-posedness of \eqref{eq:V_const_eq00} and its key qualitative properties, with the proof given in Appendix \ref{theosoluodse00}. For further details on \eqref{eq:V_const_eq00} see Section 6 of \cite{gilbarg1977}.
Notably, once the model parameters are specified, it can be solved numerically, e.g., by finite difference, see Section \ref{subsec610}. However, in our RL algorithm, it is not necessary to solve \eqref{eq:V_const_eq00} to approximate $g_i$.

\begin{proposition}\label{theosoluodse}
ODE \eqref{eq:V_const_eq00} satisfies the following properties.
\begin{enumerate}
\item[(i)]
It admits a solution $g_i\in\mathbb{C}^2([\varepsilon, 1-\varepsilon])$ for any $\varepsilon\in(0,\tfrac12)$, which is positive whenever $f_{\lambda}(0)\geq0$.
\item[(ii)] With boundary conditions \eqref{boucong1}-\eqref{addicondi}, it admits a unique solution $g_i\in\mathbb{C}^2((0, 1))\cap \mathbb{C}([0, 1])$, which is bounded and satisfies $0\leq g_i(p)\leq \tfrac{f_{\lambda}(0)}{\delta_1\wedge\delta_2}$ if $f_{\lambda}(0)\geq 0$, and $\tfrac{f_{\lambda}(0)}{\delta_1\wedge\delta_2}\leq g_i(p)\leq 0$ if $f_{\lambda}(0)< 0$.
\item[(iii)] Both $g_i'(p)$ and $g_i''(p)$ are bounded on $(0,1)$.
\end{enumerate}
\end{proposition}

Combining \eqref{eq:V_const_eq00} and \eqref{eq:V_exp_eq} yields a quadratic equation in $\kappa_i$ $(i=1, 2)$ of the form
\begin{align}\label{eqkap00}
& \tfrac{\sigma^2}{2} g_i(p) \kappa_i^2
- \Big((\mu_1 - \mu_2)p(1-p)g_i'(p) + \big(f_{\lambda}'(0) + D(p)\big)g_i(p) \Big) \kappa_i
- \hat{\varpi}_i(p) f_{\lambda}(0) = 0.
\end{align}
Since \eqref{eqkap00} cannot hold pointwise in $p$ unless $\kappa_i$ is allowed to depend on $p$, which would violate the separable ansatz \eqref{eq:single_exp_ansatz}, we analyze \eqref{eqkap00} by a Fredholm-Galerkin projection approach in the following proposition and establish sufficient conditions for $\kappa_i>0$. The proof is omitted, as it follows standard arguments.

\begin{proposition}\label{propkappa}
To eliminate the $p$-dependence in \eqref{eqkap00}, we project the identity against a suitable weight function $w$ (e.g., an invariant density for the adjoint operator of \eqref{eq:V_const_eq00}, see Appendix \ref{prop:adjoint-weight}). An integration-by-parts argument then reduces \eqref{eqkap00} to the scalar quadratic equation
\footnote{\tiny{For numerical robustness in the computations that follow, we solve \eqref{kappasolu} using Kahan's scaled quadratic formula.} }
\begin{align}\label{kappasolu}
&F_{2,i} \kappa_i^2 + F_{1,i} \kappa_i + F_{0, i} = 0,
\end{align}
where $\Phi(p)=(p(1-p)w(p))'$, and
\begin{align*}
&F_{2, i} = \tfrac{\sigma^2}{2}\int_0^1 g_i(p)w(p) \mathrm{d}p, \qquad
F_{0,i}=- f_{\lambda}(0) \int_0^1 \hat{\varpi}_i(p) w(p)\mathrm{d}p,\\
&F_{1,i} = (\mu_1-\mu_2)\int_0^1 \Phi(p) g_i(p) \mathrm{d}p
- \int_0^1 \big(f_{\lambda}'(0)+D(p)\big)g_i(p)w(p) \mathrm{d}p.
\end{align*}
Moreover, \eqref{kappasolu} admits a positive real root $\kappa_i>0$ under either of the following sufficient conditions:
\begin{align}\label{kapgeqsuff}
\text{(i)} \ F_{2,i} F_{0,i} <0, \quad
\text{or} \quad \text{(ii)} \  F_{2,i} F_{0,i} \geq 0, \ F_{1,i}^2 \geq 4 F_{2,i} F_{0,i} \ \text{and} \ F_{1,i}<0,
\end{align}
with the larger of the two positive roots selected in case (ii).
Formally, $\kappa_i = \tfrac{-F_{1,i}+\mathrm{sgn}(F_{2,i})\sqrt{F_{1,i}^2-4 F_{2,i} F_{0,i}}}{2 F_{2,i}}>0$.
\end{proposition}

The following theorem verifies that the value function for problem \eqref{opdipro} is of the form \eqref{eq:single_exp_ansatz}, the proof is deferred to Appendix \ref{verthem2200}.
\begin{theorem}\label{verthem22}
The value function of the regime-switching exploratory optimal dividend problem \eqref{opdipro} admits the form \eqref{eq:single_exp_ansatz}, where $g_i(p)$ $(i=1,2)$ solve ODE \eqref{eq:V_const_eq00} subject to boundary conditions \eqref{boucong1}-\eqref{addicondi}, and $\kappa_i>0$ are determined by \eqref{kappasolu} under the sufficient condition \eqref{kapgeqsuff}, with optimal distributional control given by \eqref{densityop}.
\end{theorem}

\section{RL algorithm design}\label{sec5}
This section develops a continuous-time actor-critic RL algorithm to solve the problem \eqref{opdipro}, based on the semi-analytical representations \eqref{densityop} and \eqref{eq:single_exp_ansatz} whose optimality is established in Theorem \ref{verthem22}.
Most RL algorithms iterate between PE and PI (see \cite{sutton2018}): PE estimates the current policy by estimating its objective value, while PI updates the policy to improve this value. Under suitable conditions, this alternation produces an objective-value sequence that exhibits an overall improving trend (e.g., tending to increase in maximization problems) and converges to the optimal value.
Beginning with the PI theorem and its convergence implication, we subsequently design an actor-critic RL algorithm, termed partially observable regime-switching exploratory optimal dividend (PO-RSEOD) algorithm, to learn the value function and the optimal distributional control for problem \eqref{opdipro}. Throughout, $\partial_x f(x, p)$ denotes the partial derivative of $f$ w.r.t. $x$, and $G(\cdot, \cdot)$ is the Gibbs function defined by \eqref{gu1y}.

\subsection{Policy improvement and policy convergence}\label{subsec 5.1}
Following Section 4 of \cite{bai2023}, we introduce the notion of strong admissibility and state the following lemma, with the proof given in Appendix \ref{lemdencovproof}.
\begin{definition}\label{def510}
A policy is said to be strongly admissible if its density $u\mapsto \pi(u| x, p)$ satisfies:
\begin{enumerate}
\item[(i)] uniform bounds: there exist constants $0<l_- \leq l_+ <\infty$ such that $l_- \leq \pi(u|x,p)\leq l_+$ for all $(x, p)\in \mathbb{R}_+ \times (0,1)$ and $u\in[0, \mathbf{a}]$.
\item[(ii)] local Lipschitz continuity in $(x, p)$, uniformly in $u$: for every compact set $K\Subset\mathbb{R}_+\times (0, 1)$, there exists $L_K>0$ such that, for $(x_i, p_i)\in K$, $i=1,2$ and $u\in[0, \mathbf a]$, \\
    $|\pi(u| x_1, p_1) - \pi(u| x_2, p_2)| \leq L_K(|x_1 -x_2| + |p_1 -p_2|)$.
\end{enumerate}
\end{definition}
If $\pi$ is strongly admissible, then $b^{\pi}(x, p)$ and $\widehat{\mathcal{H}}_{\lambda}^{\pi}(x, p)$, defined in \eqref{defbpixp} and \eqref{enth1}, are uniformly bounded on $\mathbb{R}_+\times (0,1)$ and locally Lipschitz in $(x, p)$.
More precisely, $\tfrac{\mathbf{a}^2}{2} l_- \le b^{\pi} \le \min\{\mathbf a, \tfrac{\mathbf{a}^2}{2}l_+\}<\infty$, and $-\lambda \ln l_+ \leq \widehat{\mathcal{H}}_{\lambda}^{\pi} \le \min\{\mathbf{a} + \lambda \ln\mathbf{a}, \, \mathbf{a} -\lambda \ln l_-\}<\infty$, with bounds that are uniform in $(x, p)$ and uniformly over all strongly admissible $\pi$.

For any strongly admissible policy $\pi$, by Feynman-Kac formula, $J(x, p; \pi)$ satisfies
\begin{align}\label{defjxpde}
&\big( \mathcal{G}^{\pi} - C(p) \big) J(x,p; \pi) + \widehat{\mathcal{H}}_{\lambda}^{\pi}(x, p)\nonumber\\
&:=\tfrac{\sigma^2}{2} J_{xx}(x,p; \pi)
+ \big(D(p) - b^{\pi}(x, p)\big) J_x(x,p; \pi)
+ \tfrac{1}{2} A(p) J_{pp}(x,p; \pi)
+ B(p) J_p(x,p; \pi) \nonumber \\
& \ \ + (\mu_1 - \mu_2) p (1-p) J_{xp}(x,p; \pi) -  C(p) J(x,p; \pi) + \widehat{\mathcal{H}}_{\lambda}^{\pi}(x, p) = 0,
\end{align}
with $J(0, \cdot; \pi) = 0$, and $(A, B, C, D)$, $b^{\pi}$, $\widehat{\mathcal{H}}_{\lambda}^{\pi}$ and $\mathcal{G}^{\pi}$ as defined in \eqref{apbpcpdp}, \eqref{defbpixp}, \eqref{enth1} and \eqref{def:mathcalgpi}, respectively.

\begin{lemma}\label{lemdencov}
Suppose that $\pi(u| x, p) = G(u, c_x(x,p))$ for some function $c\in \mathbb{C}^2((0,\infty) \times (0,1))$ with $|c|$, $|c_x|$ and $|c_p|$ bounded,
then $\pi$ is strongly admissible.
Moreover, for any strongly admissible policy $\pi$, $J(x, p; \pi) \in \mathbb{C}^2((0,\infty) \times (0,1))$ and $|J|$, $|J_x|$ and $|J_p|$ are bounded, with bounds depending only on those of $b^{\pi}$ and $\widehat{\mathcal{H}}_{\lambda}^{\pi}$.
\end{lemma}

We now establish the policy improvement scheme, see Appendix \ref{policyimpro00} for the proof.
\begin{theorem}[Policy improvement theorem]\label{policyimpro}
Suppose that $\pi$ is strongly admissible. If the updated policy $\tilde{\pi}$ is defined by $\tilde{\pi}(u| x, p) := G\big(u, 1-J_x(x, p; \pi)\big)$ for $(x, p)\in \mathbb{R}_+\times (0, 1)$,
then $\tilde{\pi}$ is strongly admissible, and satisfies $J(x, p; \tilde{\pi}) - J(x, p; \pi)  \geq 0$.
\end{theorem}

By Theorem \ref{policyimpro}, any strongly admissible policy admits an improving policy in the Gibbs class (truncated exponential family). We therefore restrict our learning to Gibbs policies.
The updated policy $\tilde{\pi}$ depends only on $(J, \pi)$ and not on the model parameters $(\mu_i, \sigma_i)_{i\in E}$, in this sense, the improvement step is parameter-free. However, with unknown parameters, $J(x, p; \pi)$ cannot be evaluated in closed form from the PDE \eqref{defjxpde}, so $\tilde{\pi}$ is not directly implementable.
We define the following `learning sequence' $((\pi_n, v_n))_{n\geq 0}$: let $c_{0}(x,p)\in\mathbb{C}^2((0,\infty)\times(0,1))$ with $|c_{0}|$, $|\partial_x c_0|$ and $|\partial_p c_0|$ uniformly bounded on $(0,\infty) \times (0,1)$, and set $\pi_{0}(u| x, p):=G(u, 1-\partial_x c_{0}(x,p))$, $v_{0}(x,p) := J(x,p; \pi_{0})$. Recursively, for $n\geq 1$,
\begin{align}\label{iteravalue}
& \pi_{n}(u| x,p):=G(u,,1-\partial_x v_{n-1}(x,p))=G(u, 1-J_x(x,p; \pi_{n-1})), \nonumber \\
& v_{n}(x, p):=J(x, p; \pi_{n}).
\end{align}
By Lemma \ref{lemdencov}, for every $n\geq 0$, $\pi_{n}$ is strongly admissible, $v_{n}\in\mathbb{C}^2((0,\infty)\times(0,1))$, and $|v_n|$, $|\partial_x v_n|$ and $|\partial_p v_n|$ bounded on $(0,\infty)\times(0,1)$, with bounds depending only on those of $b^{\pi_n}$ and $\widehat{\mathcal{H}}_{\lambda}^{\pi_n}$, and hence uniform in $n$ provided that $b^{\pi_n}$ and $\widehat{\mathcal{H}}_{\lambda}^{\pi_n}$ are uniformly bounded in $n$.
Consequently, starting from a strongly admissible seed $\pi_0$, the update \eqref{iteravalue} preserves strong admissibility.

We then analyze the convergence of the value and policy iterates as $n\to \infty$.
The proof follows from an argument analogous to Theorem 5 of \cite{bai2023}, based on the Arzel\`{a}-Ascoli theorem, and is therefore omitted.

\begin{theorem}[Policy convergence theorem]
The value iterates $(v_n)_{n\ge 0}$ in \eqref{iteravalue} form a maximizing sequence and converge to the value function $V$ for problem \eqref{opdipro}, and the policy iterates $(\pi_n)_{n\geq 0}$ converge to the optimal distributional control $\pi^*$ in \eqref{densityop}.
\end{theorem}

\subsection{PO-RSEOD algorithm}
We now develop an actor-critic algorithm, termed PO-RSEOD, for the problem \eqref{opdipro}.

\subsubsection{Parametrization of value function and optimal policy}
Let $v^{(\gamma, \varphi)}$ and $\pi^{(\gamma, \varphi)}$ denote the parameterized value function and policy, respectively, where $(\gamma, \varphi)$ are parameter vectors to be learned via an $m$-th order polynomial approximation. Specifically, Proposition \ref{theosoluodse} (ii) motivates approximating $g_i(p)$ by
\begin{align}\label{paragi0}
&g_i^{\varphi}(p) := \tfrac{f_{\lambda}(0)}{\delta_i} \Big(\sum_{0\leq j,k\leq m} p^j (1-p)^k e^{\varphi_{j,k}^{(i)}} \Big), \ \ \ \ i=1,2,
\end{align}
with $\varphi^{(i)} := (\varphi_{0, 0}^{(i)}, \ldots, \varphi_{m, m}^{(i)})\in\mathbb{R}^{(m+1)\times (m+1)}$, where the exponential parametrization is used to avoid negative values during training.
The weight function $w(p)$ in \eqref{eq:w-closed-2} is parameterized as
\begin{align}\label{parawgamp0}
&w^{\gamma}(p) := \tfrac{\hat{K}^{\gamma}}{\gamma_{\beta_0}}
p^{\gamma_{\beta_1} -2} (1-p)^{-\gamma_{\beta_1} -2} \exp\Big\{-\tfrac{2}{\gamma_{\beta_0}}
\big(\tfrac{e^{\gamma_3}}{p}+\tfrac{e^{\gamma_4}}{1-p}\big)\Big\},
\end{align}
with $\gamma:=(\gamma_0, \ldots, \gamma_4)$, where $\hat{K}^{\gamma}$ denote the $\gamma$-dependent normalizing constant, $\gamma_{\beta_0} := \tfrac{(e^{\gamma_1} - e^{\gamma_2})^2}{e^{\gamma_0}}$,
$\gamma_{\beta_1} := \tfrac{2}{\gamma_{\beta_0}}(e^{\gamma_3} - e^{\gamma_4})$.
By \eqref{kappasolu}, $\kappa_i^{(\gamma, \varphi)}$ $(i=1,2)$ is determined as the positive root of $F_{2,i}^{(\gamma, \varphi)} \kappa_i^2 + F_{1,i}^{(\gamma, \varphi)} \kappa_i + F_{0, i}^{(\gamma, \varphi)} = 0$, where $\Phi^{\gamma}(p):=(p(1-p)w^{\gamma}(p))'$, and
\begin{align*}
&F_{2, i}^{(\gamma, \varphi)} := \tfrac{e^{\gamma_0}}{2}\int_0^1 g_i^{\varphi}(p)w^{\gamma}(p) \mathrm{d}p, \qquad
F_{0,i}^{(\gamma, \varphi)} := -\int_0^1 C(p) g_i^{\varphi}(p)  w^{\gamma}(p) \mathrm{d}p, \\
&F_{1,i}^{(\gamma, \varphi)} := (e^{\gamma_1}- e^{\gamma_2})\int_0^1 \Phi^{\gamma}(p) g_i^{\varphi}(p) \mathrm{d}p
- \int_0^1 \big(f_{\lambda}'(0)+ e^{\gamma_2} +(e^{\gamma_1} - e^{\gamma_2}) p \big)g_i^{\varphi}(p)w^{\gamma}(p) \mathrm{d}p.
\end{align*}
That is, $\kappa_i^{(\gamma, \varphi)} := \tfrac{-F_{1,i}^{(\gamma, \varphi)}+\mathrm{sgn}(F_{2,i}^{(\gamma, \varphi)}) \sqrt{\textbf{F}_i^{(\gamma,\varphi)}}}{2 F_{2,i}^{(\gamma, \varphi)}} $ with $\textbf{F}_i^{(\gamma,\varphi)}
:=(F_{1,i}^{(\gamma,\varphi)})^{2}
-4\,F_{2,i}^{(\gamma,\varphi)}\,F_{0,i}^{(\gamma,\varphi)}>0$.
Let $\Theta\subset\mathbb{R}^{d_\gamma}\times\mathbb{R}^{d_\varphi}$ be compact and denote $\boldsymbol{\theta}=(\gamma,\varphi)\in\Theta$ and $\theta\in\{\gamma,\varphi\}$.
Under standard regularity (dominated-convergence) conditions,
$F_{k,i}^{(\gamma,\varphi)}$ and $\partial_{\theta} F_{k,i}^{(\gamma,\varphi)}$ $(k\in\{0,1,2\})$ are well defined and continuous on $\Theta$. We state the following proposition.

\begin{proposition}\label{prop51para}
The implicit function theorem yields that
$\kappa_i^{(\gamma,\varphi)}$ is differentiable w.r.t. $\theta\in\{\gamma,\varphi\}$ since $2\,F_{2,i}^{(\gamma,\varphi)}\,\kappa_i^{(\gamma,\varphi)}
+F_{1,i}^{(\gamma,\varphi)}
=\sqrt{\textbf{F}_i^{(\gamma,\varphi)}}\neq 0$, and its derivative is
\begin{align*}
&\partial_{\theta} \kappa_i^{(\gamma,\varphi)}
=-
\tfrac{\big(\kappa_i^{(\gamma,\varphi)}\big)^2\,\partial_{\theta} F_{2,i}^{(\gamma,\varphi)}
+\kappa_i^{(\gamma,\varphi)}\,\partial_{\theta}  F_{1,i}^{(\gamma,\varphi)}
+\partial_{\theta} F_{0,i}^{(\gamma,\varphi)}
}{\sqrt{\textbf{F}_i^{(\gamma,\varphi)}}},
\end{align*}
where for $0\leq j, k\leq m$, the partial derivatives of $F_{k, i}^{(\gamma,\varphi)}$ w.r.t. $\varphi_{j,k}^{(i)}$ are denoted by
\begin{align*}
\partial_{\varphi_{j, k}^{(i)}} F_{2,i}^{(\gamma,\varphi)}
&= \tfrac{e^{\gamma_0}}{2}\int_{0}^{1} \big(\partial_{\varphi_{j, k}^{(i)}} g_i^{\varphi}(p)\big)\,w^{\gamma}(p)\,\mathrm{d} p, \\
\partial_{\varphi_{j, k}^{(i)}} F_{0,i}^{(\gamma,\varphi)}
&= -\int_{0}^{1} C(p)\,\big(\partial_{\varphi_{j, k}^{(i)}} g_i^{\varphi}(p)\big)\,w^{\gamma}(p)\,\mathrm{d}p, \\
\partial_{\varphi_{j, k}^{(i)}} F_{1,i}^{(\gamma,\varphi)}
&= \big(e^{\gamma_1}-e^{\gamma_2}\big)\int_{0}^{1} \Phi^{\gamma}(p)\,\big(\partial_{\varphi_{j, k}^{(i)}} g_i^{\varphi}(p)\big)\,\mathrm{d}p
\end{align*}
\begin{align*}
& \ \ \ - \int_{0}^{1}\Big(f'_{\lambda}(0)+e^{\gamma_2}+\big(e^{\gamma_1}-e^{\gamma_2}\big)p\Big)
\,\big(\partial_{\varphi_{j, k}^{(i)}} g_i^{\varphi}(p)\big)\,w^{\gamma}(p)\,\mathrm{d}p,
\end{align*}
with $\partial_{\varphi_{j, k}^{(i)}} g_i^{\varphi}(p) = \tfrac{f_{\lambda}(0)}{\delta_i} p^j (1-p)^k e^{\varphi_{j,k}^{(i)}}$,
and the derivatives of $F_{k, i}^{(\gamma,\varphi)}$ w.r.t. $\gamma_j$ are
\begin{align*}
\partial_{\gamma_j} F_{2,i}^{(\gamma,\varphi)}
&= \mathbf 1_{\{j=0\}}\,F_{2,i}^{(\gamma,\varphi)}
   + \tfrac{e^{\gamma_0}}{2}\int_{0}^{1} g_i^{\varphi}(p)\,\partial_{\gamma_j}w^{\gamma}(p)\,\mathrm{d}p, \\
\partial_{\gamma_j} F_{0,i}^{(\gamma,\varphi)}
&= -\int_{0}^{1} C(p)\,g_i^{\varphi}(p)\,\partial_{\gamma_j}w^{\gamma}(p)\,\mathrm{d}p,\nonumber\\
\partial_{\gamma_j} F_{1,i}^{(\gamma,\varphi)}
&= \big(\mathbf 1_{\{j=1\}}e^{\gamma_1}-\mathbf 1_{\{j=2\}}e^{\gamma_2}\big)
      \int_{0}^{1}\Phi^{\gamma}(p)\,g_i^{\varphi}(p)\,\mathrm{d}p
+ \big(e^{\gamma_1}-e^{\gamma_2}\big)
      \int_{0}^{1}\partial_{\gamma_j}\Phi^{\gamma}(p)\,g_i^{\varphi}(p)\,\mathrm{d}p \\
&\ \ \ - \int_{0}^{1}\Big(\mathbf 1_{\{j=2\}}e^{\gamma_2}
   + \big(\mathbf 1_{\{j=1\}}e^{\gamma_1}-\mathbf 1_{\{j=2\}}e^{\gamma_2}\big)p\Big)
   g_i^{\varphi}(p)\,w^{\gamma}(p)\,\mathrm{d}p \\
&\ \ \
- \int_{0}^{1}\Big(f'_{\lambda}(0)+e^{\gamma_2}
   +\big(e^{\gamma_1}-e^{\gamma_2}\big)p\Big)\,
   g_i^{\varphi}(p)\,\partial_{\gamma_j}w^{\gamma}(p)\,\mathrm{d}p,
\end{align*}
with
$\partial_{\gamma_j} w^{\gamma}(p)= w^{\gamma}(p)\,\Xi_j(p)$ and $\partial_{\gamma_j}\Phi^{\gamma}(p) = \big(p(1-p)\,w^{\gamma}(p)\,\Xi_j(p)\big)'$,
where $\Xi_j$ denotes the log-derivative,
\begin{align*}
\Xi_j(p)
:= \partial_{\gamma_j}\big(\ln w^{\gamma}(p)\big) \notag
&= \partial_{\gamma_j}\ln \widehat K^{\gamma}
- \partial_{\gamma_j}\ln \gamma_{\beta 0}
+ \big(\partial_{\gamma_j}\gamma_{\beta 1}\big)\big(\ln p-\ln(1-p)\big)\notag\\
&\quad
+ \tfrac{2\,\partial_{\gamma_j}\gamma_{\beta 0}}{\gamma_{\beta 0}^{2}}
    \left(\tfrac{e^{\gamma_3}}{p}+\tfrac{e^{\gamma_4}}{1-p}\right) - \tfrac{2}{\gamma_{\beta 0}}
  \left(\tfrac{\mathbf 1_{\{j=3\}}\,e^{\gamma_3}}{p}
       +\tfrac{\mathbf 1_{\{j=4\}}\,e^{\gamma_4}}{1-p}\right),
\end{align*}
and $\partial_{\gamma_0}\gamma_{\beta 0}=-\gamma_{\beta 0}$,
$\partial_{\gamma_1}\gamma_{\beta 0}=\frac{2e^{\gamma_1}}{e^{\gamma_1}-e^{\gamma_2}}\gamma_{\beta 0}$,
$\partial_{\gamma_2}\gamma_{\beta 0}=-\frac{2e^{\gamma_2}}{e^{\gamma_1}-e^{\gamma_2}}\gamma_{\beta 0}$,
$\partial_{\gamma_i}\gamma_{\beta 1}= -\frac{\gamma_{\beta 1}}{\gamma_{\beta 0}}\partial_{\gamma_i}\gamma_{\beta 0}$ for $i=0,1,2$, and
$\partial_{\gamma_3}\gamma_{\beta 1}=\frac{2e^{\gamma_3}}{\gamma_{\beta 0}}$,
$\partial_{\gamma_4}\gamma_{\beta 1}=-\frac{2e^{\gamma_4}}{\gamma_{\beta 0}}$.
\end{proposition}
Consequently, the value function and the policy are parameterized by $(\gamma,\varphi)$ as
\begin{align}
&v^{\boldsymbol{\theta}}(x, p)=v^{(\gamma, \varphi)}(x, p) := g_1^{\varphi}(p) \big(1 - e^{-\kappa_1^{(\gamma, \varphi)}  x}\big) + g_2^{\varphi}(p) \big(1 - e^{-\kappa_2^{(\gamma, \varphi)} x}\big), \label{eq:single_exp_ansatz00}\\
&\pi^{\boldsymbol{\theta}}(u| x, p) = \pi^{(\gamma, \varphi)}(u| x, p) := G\big(u, 1-v_x^{(\gamma, \varphi)}(x, p)\big).
\label{eq:single_exp_ansatz11}
\end{align}
The state dynamics \eqref{sde5} under $\pi^{\boldsymbol{\theta}}=\pi^{(\gamma, \varphi)}$ take the form
\begin{align}\label{sde50012}
	\mathrm{d} \widehat{X}_t^{\boldsymbol{\theta}}
:=& \Big(\widehat{\mu}_{t}
- \int_0^{\mathbf{a}} u \pi^{\boldsymbol{\theta}}(u| \widehat{X}_t^{\boldsymbol{\theta}}, p_t)\mathrm{d}u\Big)\mathrm{d} t + \sigma \mathrm{d} \widehat{W}_t,
\end{align}
with $\widehat{X}_0^{\boldsymbol{\theta}}= x$ and $\widehat{X}_{\widehat{T}_{x, p}^{\boldsymbol{\theta}}}^{\boldsymbol{\theta}}= 0$, where $\widehat{T}_{x, p}^{\boldsymbol{\theta}} :=\inf\{t\geq 0: \widehat{X}_t^{\boldsymbol{\theta}}\leq0\}$ is the ruin time.
The function $\widehat{\mathcal{H}}_{\lambda}^{\pi}(x, p)$ in \eqref{h526}, evaluated at $\pi=\pi^*$, can be parameterized as
\begin{align}\label{h526pare}
\widehat{\mathcal{H}}_{\lambda}^{\boldsymbol{\theta}}(x, p)
& = \Big[
v_x^{\boldsymbol{\theta}} \Big(\mathbf{a} - \tfrac{\lambda}{1- v_x^{\boldsymbol{\theta}}}
+ \tfrac{\mathbf{a}}{e^{\frac{\mathbf{a}}{\lambda}(1-v_x^{\boldsymbol{\theta}})} -1}\Big)
- \lambda \ln \Big(\tfrac{1- v_x^{\boldsymbol{\theta}}}{\lambda(e^{\frac{\mathbf{a}}{\lambda}(1-v_x^{\boldsymbol{\theta}})} -1)}\Big) \Big]\mathbf{1}_{\{v_x^{\boldsymbol{\theta}}\neq 1\}} \nonumber\\
& \ \ \ + \big[\tfrac{\mathbf{a}}{2} + \lambda \ln \mathbf{a}\big] \mathbf{1}_{\{v_x^{\boldsymbol{\theta}}= 1 \}},
\end{align}
with $v_x^{\boldsymbol{\theta}} = v_x^{\boldsymbol{\theta}}(x, p)$.
Express $\widehat{\mathcal H}_\lambda^{\boldsymbol{\theta}}$ in \eqref{h526pare} as a scalar function of
$v_x^{\boldsymbol{\theta}}$:
$\widehat{\mathcal H}_\lambda^{\boldsymbol{\theta}}(x,p)=\mathscr H_\lambda(v_x^{\boldsymbol{\theta}})$. Its sensitivity is $\mathscr{H}_\lambda'(v_x^{\boldsymbol{\theta}}):=\Upsilon_\lambda(v_x^{\boldsymbol{\theta}})$ with $\Upsilon_\lambda(1):=\lim_{v_x^{\boldsymbol{\theta}}\to 1}\Upsilon_\lambda(v_x^{\boldsymbol{\theta}})
= -\tfrac{\mathbf{a}^{2}}{12\,\lambda}$, for $v_x^{\boldsymbol{\theta}}\neq 1$,
\begin{align}\label{scalarsenh}
\Upsilon_\lambda(v_x^{\boldsymbol{\theta}})
:= \mathbf{a}+\tfrac{\mathbf{a}}{e^{\frac{\mathbf{a}}{\lambda}(1-v_x^{\boldsymbol{\theta}})}-1}
  + v_x^{\boldsymbol{\theta}} \Big(-\tfrac{\lambda}{(1-v_x^{\boldsymbol{\theta}})^2}
             +\tfrac{\mathbf{a}^2}{\lambda}\,
               \tfrac{e^{\frac{\mathbf{a}}{\lambda}(1-v_x^{\boldsymbol{\theta}})}}
               {\big(e^{\frac{\mathbf{a}}{\lambda}(1-v_x^{\boldsymbol{\theta}})}-1\big)^2}\Big)
  - \tfrac{\mathbf{a}\,e^{\frac{\mathbf{a}}{\lambda}(1-v_x^{\boldsymbol{\theta}})}}
  {e^{\frac{\mathbf{a}}{\lambda}(1-v_x^{\boldsymbol{\theta}})}-1}.
\end{align}

Alternative parametrizations are possible, notably neural-network approximation as in model-free RL algorithms for continuous state-action problems (e.g., DDPG, SAC).
Such approaches typically entail high-dimensional parametrizations, which can slow learning and aggravate optimization instability and out-of-sample degradation, see \cite{wang20202} and \cite{wang2023}. In contrast, our semi-analytical, structure-exploiting parametrization captures both the value function and the optimal policy with only a few parameters, thereby enhancing learning stability and statistical efficiency.

\subsubsection{Interaction and filtering}\label{sec522}
RL updates the policy and its value estimate via iterative, discrete-time interactions with the environment, even when the control problem is continuous-time. We first introduce a finite-horizon, discrete-time interaction scheme.
\begin{definition}[Discrete-time interaction scheme]\label{discrete-time-scheme}
Fix a step size $\Delta t > 0$ and set $t_k = k \Delta t$ $(k=0, \ldots, K)$ over the finite time horizon $[0, T\wedge \widehat{T}_{x, p}^{\pi}]$, where $T>0$ is a fixed sufficiently large truncation horizon and $K:= K_1 \wedge K_2$ with $K_1:= \lfloor T / \Delta t \rfloor$ and $K_2 := \lfloor \widehat{T}_{x, p}^{\pi} / \Delta t \rfloor$.
The set $\{t_k\}_{k=0}^K$ is called the discrete-time interaction grid.
\end{definition}
When all states are observable, interaction-based learning obviates the need to estimate environmental parameters. Under partial observability, however, policy learning must be coupled with filtering of the unobserved state, which in turn requires parameter estimates. Prior to learning, we estimate the parameters $\mu_i$, $\sigma$ and $q_{ij}$ $(i,j=1,2)$ from historical data, which can be done in various ways (e.g., \cite{elliott1995}, \cite{krishnamurthy2016}), and empirically, filtering accuracy is relatively insensitive to moderate estimation error, see \cite{wu2024}.

In our problem \eqref{opdipro}, at each grid time $t_k$, the surplus $\widehat{X}_{t_k}^{\pi}$ is observed, whereas the regime $I_{t_k}$ is unobserved and must be inferred via filtering, yielding the belief $p_{t_k}$. To this end, let $\widehat{\mu}_{t_{k-1}} = (\mu_1 - \mu_2)p_{t_{k-1}} + \mu_2$, we discretize the innovation process \eqref{hatwt} with increment
\begin{align}\label{hatwt221}
& \Delta \widehat{W}_{t_k} = \widehat{W}_{t_k} - \widehat{W}_{t_{k-1}}
= \tfrac{X_{t_k} - X_{t_{k-1}} - \widehat{\mu}_{t_{k-1}} \Delta t}{\sigma}.
\end{align}
The belief $p_{t_k}$ is updated via the Wonham filter \eqref{wondpt}.
Notably, direct Euler discretization of \eqref{wondpt} may produce iterates outside $[0,1]$, invalidating probabilistic interpretation. To enforce the unit-interval constraint, we use the transformations $\tilde{p}_t^{(1)}:=\ln p_t$ and $\tilde{p}_t^{(2)}:=\ln (1-p_t)$, apply It\^{o}'s formula to obtain SDEs for $\tilde{p}_t^{(1)}$ and $\tilde{p}_t^{(2)}$, and discretize them via the Euler-Maruyama scheme to yield the discrete-time recursion
\begin{align*}
\tilde{p}_{t_k}^{(1)} &= \tilde{p}_{t_{k-1}}^{(1)}
+\big( q_{21} e^{-\tilde{p}_{t_{k-1}}^{(1)}}
- (q_{12}+ q_{21})
-\tfrac{1}{2} \big(\tfrac{\mu_1 - \mu_2}{\sigma}\big)^2
(1- e^{\tilde{p}_{t_{k-1}}^{(1)}})^2 \big) \Delta t \\
& \ \ \ +\tfrac{\mu_1 - \mu_2}{\sigma} (1- e^{\tilde{p}_{t_{k-1}}^{(1)}}) \Delta \widehat{W}_{t_k}, \\
 \tilde{p}_{t_k}^{(2)} &= \tilde{p}_{t_{k-1}}^{(2)}
+\big( q_{12} e^{-\tilde{p}_{t_{k-1}}^{(2)}}
- (q_{12}+ q_{21})
-\tfrac{1}{2} \big(\tfrac{\mu_1 - \mu_2}{\sigma}\big)^2
(1- e^{\tilde{p}_{t_{k-1}}^{(2)}})^2 \big) \Delta t \\
& \ \ \
- \tfrac{\mu_1 - \mu_2}{\sigma} (1- e^{\tilde{p}_{t_{k-1}}^{(2)}}) \Delta \widehat{W}_{t_k},
\end{align*}
then $p_{t_k}$ is recovered as $p_{t_k} = e^{\tilde{p}_{t_k}^{(1)}}/(e^{\tilde{p}_{t_k}^{(1)}} +e^{\tilde{p}_{t_k}^{(2)}})\in(0, 1)$.
Excellent performance of this scheme was shown in \cite{yin2004}.

An episode $\mathcal{D}$ is generated by interacting with the environment under a given policy $\pi$, as follows. Given $(\widehat{X}_{t_k}^{\pi}, p_{t_k})$, we sample action $u_{t_k}\in[0, \mathbf{a}]$ from the stochastic policy
$\pi(u | \widehat{X}_{t_k}^{\pi}, p_{t_k}) = G(u, 1-J_x(\widehat{X}_{t_k}^{\pi}, p_{t_k}; \pi))$, and keep it constant on $[t_k, t_{k+1})$. The surplus state is then advanced by the Euler-Maruyama update,
\begin{align*}
&\widehat{X}_{t_{k+1}}^{\pi}
= \widehat{X}_{t_k}^{\pi} + (\widehat{\mu}_{t_k} - u_{t_k}) \Delta t
+ \sigma \Delta \widehat{W}_{t_{k+1}},
\end{align*}
where $\Delta\widehat W_{t_{k+1}}$ denotes the innovation increment in \eqref{hatwt221}.
Iterating until $T\wedge \widehat{T}_{x,p}^{\pi}$ yields the state-action trajectory (episode) $\mathcal{D}=\{(t_k, \widehat{X}_{t_k}^{\pi}, p_{t_k}, u_{t_k})\}_{k=0}^{K-1}$.
As noted in Section \ref{subsec32}, $u_{t_k}$ is sampled by the inverse of cumulative distribution function $(F^{\pi})^{-1}(Z_{t_k})$ with $Z_{t_k}\sim \mathrm{Unif}[0,1]$, and by \eqref{densityop}, specifies the sampling rule
\begin{align*}
&u_{t_k} := (F^{\pi})^{-1}(Z_{t_k})
=\tfrac{\lambda \ln\big(Z_{t_k} (e^{\frac{\mathbf{a} (1-v_x^{\pi})}{\lambda}} -1) +1\big)}{1-v_x^{\pi}} \mathbf{1}_{\{v_x^{\pi}\neq 1\}}
+\mathbf{a}\, Z_{t_k} \mathbf{1}_{\{v_x^{\pi}= 1\}} \in [0, \mathbf{a}].
\end{align*}

\subsubsection{Policy evaluation (PE)}
Following \cite{jia20222}, we employ two martingale-based methods for PE. One method is to minimize the martingale loss (ML) function, which yields the optimal approximation to the true value function in the mean-square sense and leads to an offline algorithm. The other method is to solve a system of moment equations arising from martingale orthogonality conditions w.r.t. (finitely many) test functions, which leads to online continuous-time temporal-difference (CTD) algorithms, with the specific variant determined by the chosen solution scheme.

Fix a strongly admissible policy $\pi$ as defined in Definition \ref{def510}, and for notational convenience write $J(\widehat{X}_t^{\pi}, p_t; \pi) = v^{\pi}(\widehat{X}_t^{\pi}, p_t)$. By Lemma \ref{lemdencov}, $v^{\pi} \in \mathbb{C}^2((0,\infty) \times (0,1))$ and $|v^{\pi}|$, $|v_x^{\pi}|$ and $|v_p^{\pi}|$ are bounded on $(0,\infty) \times (0,1)$. Since $v^{\pi}$ solves the PDE \eqref{defjxpde}, and furthermore, It\^{o}'s formula yields \eqref{elamvxpitoexp}, we have, for $t\geq 0$,
\begin{align}\label{mtpi12}
&M_t^{\pi} := e^{-\Lambda_t} v^{\pi}(\widehat{X}_t^{\pi}, p_t) + \int_0^t e^{-\Lambda_s} \widehat{\mathcal{H}}_{\lambda}^{\pi}(\widehat{X}_s^{\pi}, p_s) \mathrm{d} s
\end{align}
is a martingale w.r.t. $\mathcal{F}_t$. For $t\in[0, T]$, by optional sampling theorem, one has $M_{T\wedge \widehat{T}_{x, p}^{\pi}\wedge t} = M_{\widehat{T}_{x, p}^{\pi}\wedge t}$, and it is a martingale w.r.t. $\mathcal{F}_{T\wedge \widehat{T}_{x, p}^{\pi}\wedge t}$.
We now state the following standard result, which implies that the martingale condition \eqref{mtpi12} uniquely identifies $v^{\pi}$ for fixed $\pi$, with the proof deferred to Appendix \ref{prop52a00}.

\begin{proposition}\label{prop52a}
Suppose that $\pi$ is strongly admissible, $\tilde{v}^{\pi}\in \mathbb{C}^2((0,\infty) \times (0,1))$ is bounded with $\tilde{v}^{\pi}(0, \cdot) = 0$, and that the process $\tilde{M}^{\pi} := (\tilde{M}_t^{\pi})_{t\geq 0}$ with $\tilde{M}_t^{\pi} := e^{-\Lambda_t} \tilde{v}^{\pi}(\widehat{X}_t^{\pi}, p_t) + \int_0^t e^{-\Lambda_s} \widehat{\mathcal{H}}_{\lambda}^{\pi}(\widehat{X}_s^{\pi}, p_s) \mathrm{d} s$ is an $\mathbb{F}=\{\mathcal{F}_t\}_{t\geq 0}$-martingale. Then $v^{\pi}\equiv \tilde{v}^{\pi}$.
\end{proposition}
To learn the parameter vector $\boldsymbol{\theta}=(\gamma,\varphi)\in\Theta$, given the parametric value function ${v}^{\boldsymbol{\theta}}$ in \eqref{eq:single_exp_ansatz00} associated with the policy $\pi^{\boldsymbol{\theta}}$ in \eqref{eq:single_exp_ansatz11}, in view of \eqref{sde50012} and \eqref{h526pare}, we parameterize $M_t^{\pi}$ in \eqref{mtpi12} as
\begin{align}\label{martingavar}
&M_t^{\boldsymbol{\theta}} := e^{-\Lambda_t} v^{\boldsymbol{\theta}}(\widehat{X}_t^{\boldsymbol{\theta}}, p_t)
+ \int_0^t e^{-\Lambda_s} \widehat{\mathcal{H}}_{\lambda}^{\boldsymbol{\theta}}(\widehat{X}_s^{\boldsymbol{\theta}}, p_s) \mathrm{d} s.
\end{align}
\textbf{Mode 1: ML Minimization.}
The first method minimizes the ML function $\mathrm{ML}(\boldsymbol{\theta})$,
\begin{align}\label{defmlfunc}
\mathrm{ML}(\boldsymbol{\theta})
&:=\frac{1}{2} \mathbb{E}
\big[\int_0^{T \wedge \widehat{T}_{x, p}^{\boldsymbol{\theta}}} \big| M_{T \wedge\widehat{T}_{x, p}^{\boldsymbol{\theta}}}^{\boldsymbol{\theta}} - M_t^{\boldsymbol{\theta}} \big| ^2 \mathrm{d} t \big]
=\frac{1}{2} \mathbb{E}
\big[\int_0^{T \wedge \widehat{T}_{x, p}^{\boldsymbol{\theta}}} \big| m^{\boldsymbol{\theta}}(\widehat{X}_t^{\boldsymbol{\theta}}, p_t) \big| ^2 \mathrm{d} t \big],
\end{align}
where
$m^{\boldsymbol{\theta}}(\widehat{X}_t^{\boldsymbol{\theta}}, p_t) := e^{-\Lambda_t} v^{\boldsymbol{\theta}}(\widehat{X}_t^{\boldsymbol{\theta}}, p_t) - \int_t^{T \wedge \widehat{T}_{x, p}^{\boldsymbol{\theta}}} e^{-\Lambda_s} \widehat{\mathcal{H}}_{\lambda}^{\boldsymbol{\theta}}(\widehat{X}_s^{\boldsymbol{\theta}}, p_s) \mathrm{d} s$.

By Theorem 6 of \cite{bai2023}, the minimizers of $\mathrm{ML}(\boldsymbol{\theta})$ coincide with those of the mean-squared value error (MSVE), i.e., $\arg\min_{\boldsymbol{\theta}\in\Theta}\mathrm{ML}(\boldsymbol{\theta})
=\arg\min_{\boldsymbol{\theta}\in\Theta}\mathrm{MSVE}(\boldsymbol{\theta})$, where $\mathrm{MSVE}(\boldsymbol{\theta})
:=\mathbb{E}[\int_{0}^{T\wedge \widehat{T}_{x,p}^{\boldsymbol{\theta}}}
|v^{\boldsymbol{\theta}}(\widehat X_{t}^{\boldsymbol{\theta}},p_t)
- v^{\pi^*}(\widehat X_{t}^{\pi^*},p_t)|^{2} \mathrm{d}t]$.
This implies that, for any $\boldsymbol{\theta}^* \in \arg\min_{\boldsymbol{\theta}\in\Theta} \mathrm{ML}(\boldsymbol{\theta})$, the parametric value function $v^{\boldsymbol{\theta}^*}$ is the $\mathbb{L}^2$-best (MSVE-optimal) approximation to the true value function $v^{\pi^*}$ within the class $(v^{\boldsymbol{\theta}})_{\boldsymbol{\theta}\in \Theta}$.

In the discrete-time interaction scheme (Definition \ref{discrete-time-scheme}), let $e^{-\Lambda_{t_i}}:=e^{-\sum_{j=0}^{i} \widehat{\delta}_{t_j} \Delta t}$ with $\widehat{\delta}_{t_j} = (\delta_1 - \delta_2)p_{t_j} + \delta_2$, and $\widehat{\mathcal{H}}_{\lambda}^{\boldsymbol{\theta}}$ as in \eqref{h526pare}, $\mathrm{ML}(\boldsymbol{\theta})$ in \eqref{defmlfunc} takes the discrete form
\begin{align}\label{defmlfunc01}
\mathrm{ML}_{\Delta t}(\boldsymbol{\theta})
& = \frac{1}{2} \mathbb{E}
\big[\sum_{k=0}^{K-1}
\big(m^{\boldsymbol{\theta}}(\widehat{X}_{t_k}^{\boldsymbol{\theta}}, p_{t_k})\big)^2 \, \Delta t \big],
\end{align}
with $m^{\boldsymbol{\theta}}(\widehat{X}_{t_k}^{\boldsymbol{\theta}}, p_{t_k})
:=e^{-\Lambda_{t_k}} \,
v^{\boldsymbol{\theta}}(\widehat{X}_{t_k}^{\boldsymbol{\theta}}, p_{t_k})
- \sum_{j=k}^{K - 1} e^{-\Lambda_{t_j}} \widehat{\mathcal{H}}_{\lambda}^{\boldsymbol{\theta}}(\widehat{X}_{t_j}^{\boldsymbol{\theta}}, p_{t_j}) \, \Delta t$ the (HJB) martingale residual.
By Lemma 1.1 of \cite{jia2022} and Theorem 7 of \cite{bai2023}, the following proposition holds, which provides a principled basis for approximating the value function $v^{\pi^*}$.
\begin{proposition}
Let $\Theta\subset\mathbb{R}^{d_\gamma}\times\mathbb{R}^{d_\varphi}$ be a nonempty compact set, and suppose that
$\mathrm{ML}_{\Delta t}(\boldsymbol{\theta})$ and $\mathrm{ML}(\boldsymbol{\theta})$ are continuous on $\boldsymbol{\theta}=(\gamma, \varphi)\in \Theta$. Then,
\begin{enumerate}
\item[(i)] $\mathrm{ML}_{\Delta t}(\boldsymbol{\theta}) \to \mathrm{ML}(\boldsymbol{\theta})$ uniformly in $\boldsymbol{\theta} \in\Theta$ as $\Delta t\to0$;
\item[(ii)] there exists
$\boldsymbol{\theta}_{\Delta t}^{*}\in\arg\min_{\boldsymbol{\theta}\in\Theta}\mathrm{ML}_{\Delta t}(\boldsymbol{\theta})$ for any $\Delta t>0$ by Weierstrass' theorem;
\item[(iii)] $\lim_{\Delta t^{(n)}\downarrow 0} \boldsymbol{\theta}_{\Delta t^{(n)}}^*
= \boldsymbol{\theta}^* \in \arg\min_{\boldsymbol{\theta}\in\Theta} \mathrm{ML}(\boldsymbol{\theta})$ for a sequence $(\Delta t^{(n)})_{n\geq 0}\downarrow 0$, and $v^{\boldsymbol{\theta}^*}$ is the MSVE-optimal approximation to the true value function $v^{\pi^*}$.
\end{enumerate}
\end{proposition}

To stabilize the joint estimation of the policy and environment parameters and to prevent unrealistic excursions of $\gamma$ and $(g(0), g(1))$, we augment the discretized ML function \eqref{defmlfunc01} with two quadratic regularizers:
an $\mathbb{L}^2$-penalty on the environment parameters $e^{\gamma} = (\mathrm{e}^{\gamma_0},\dots,\mathrm{e}^{\gamma_4})$, and a boundary-matching penalty on the functions $g_1, g_2$.
Specifically, for given weights $w^{\mathrm{env}}=(w_0^{\mathrm{env}},\dots,w_4^{\mathrm{env}})\ge 0$ and $w^{\mathrm{bc}}=(w_0^{\mathrm{bc}},w_1^{\mathrm{bc}})\ge 0$, and boundary targets $g_{\mathrm{bc}}:=(g(0),g(1))$, the regularized ML function is defined by
\begin{align}\label{defmlfunc0100}
\mathrm{ML}^{\mathrm{reg}}_{\Delta t}(\boldsymbol{\theta})
& := \frac{1}{2} \mathbb{E}
\big[\sum_{k=0}^{K-1}
\big(m^{\boldsymbol{\theta}}(\widehat{X}_{t_k}^{\boldsymbol{\theta}}, p_{t_k})\big)^2 \, \Delta t \big] +\frac{1}{2} \sum_{j=0}^4 w_j^{\mathrm{env}}(e^{\gamma_j} - \mathrm{env}_j)^2 \\
& \ \ \ + \frac{1}{2}\big[
w_0^{\mathrm{bc}}\big(g_1^{\varphi}(0)+g_2^{\varphi}(0)-g(0)\big)^2
+ w_1^{\mathrm{bc}}\big(g_1^{\varphi}(1)+g_2^{\varphi}(1)-g(1)\big)^2\big], \nonumber
\end{align}
where $\mathrm{env}:=(\sigma^2,\mu_1,\mu_2,q_{21},q_{12})$ is the reference (true/estimated) environment vector.
After sampling (see Section \ref{sec522}), we approximate \eqref{defmlfunc0100} by data in an episode $\mathcal{D}=\{(t_k, \widehat{X}_{t_k}^{\boldsymbol{\theta}}, p_{t_k}, u_{t_k})\}_{k=0}^{K-1}$,
and utilize \emph{Stochastic Gradient Descent} (SGD) algorithm to minimize $\mathrm{ML}_{\Delta t}^{\text{reg}}(\boldsymbol{\theta}; \mathcal{D})$ and design the parameter updating rule for $\boldsymbol{\theta}$. That is, for $\boldsymbol{\theta}=(\gamma, \varphi)$, $\nabla_{\boldsymbol{\theta}}=(\partial_{\gamma}, \partial_{\varphi})$ and $(t_k, \widehat{X}_{t_k}^{\boldsymbol{\theta}}, p_{t_k}, u_{t_k})\in \mathcal{D}$, the gradient of $\mathrm{ML}_{\Delta t}^{\text{reg}}$ is given by
\begin{align}\label{theta01}
& \nabla_{\boldsymbol{\theta}} \mathrm{ML}_{\Delta t}^{\text{reg}}(\boldsymbol{\theta}; \mathcal{D})
=\sum_{k=0}^{K-1}
m^{\boldsymbol{\theta}}(\widehat{X}_{t_k}^{\boldsymbol{\theta}}, p_{t_k}) \, \Delta t \cdot \nabla_{\boldsymbol{\theta}} m^{\boldsymbol{\theta}}(\widehat{X}_{t_k}^{\boldsymbol{\theta}}, p_{t_k})
\\
& + \Big(0,\big(w_j^{\mathrm{env}}(\mathrm{e}^{\gamma_j}-\mathrm{env}_j) \mathrm{e}^{\gamma_j}\big)_{j=0,\dots,4}\Big)
+ \Big(
\big(w_0^{\mathrm{bc}} e_0 \nabla_{\varphi_i}g_i^{\varphi}(0)
+ w_1^{\mathrm{bc}} e_1 \nabla_{\varphi_i}g_i^{\varphi}(1) \big)_{i=1,2},\;
0 \Big),\nonumber
\end{align}
where $e_i := g_1^{\varphi}(i)+g_2^{\varphi}(i)-g(i)$ for $i=0,1$, and
\begin{align*}
\nabla_{\boldsymbol{\theta}} m^{\boldsymbol{\theta}}(\widehat{X}_{t_k}^{\boldsymbol{\theta}}, p_{t_k})
= e^{-\Lambda_{t_k}}\,\nabla_{\boldsymbol{\theta}} v^{\boldsymbol{\theta}}(\widehat X_{t_k}^{\boldsymbol{\theta}},p_{t_k})
-\sum_{j=k}^{K-1}e^{-\Lambda_{t_j}}\,
\nabla_{\boldsymbol{\theta}} \widehat{\mathcal H}_\lambda^{\boldsymbol{\theta}}(\widehat X_{t_j}^{\boldsymbol{\theta}}, p_{t_j}) \, \Delta t,
\end{align*}
with $\nabla_{\boldsymbol{\theta}} \widehat{\mathcal H}_\lambda^{\boldsymbol{\theta}}(x,p)
= \Upsilon_\lambda\big(v_x^{\boldsymbol{\theta}}(x,p)\big)\;
  \nabla_{\boldsymbol{\theta}} v_x^{\boldsymbol{\theta}}(x,p)$ and $\Upsilon_\lambda$ defined in \eqref{scalarsenh}, and by \eqref{eq:single_exp_ansatz00},
\begin{align}
\nabla_{\boldsymbol{\theta}} v^{\boldsymbol{\theta}}(x,p)
&=\sum_{i=1}^{2}\big[
  \nabla_{\boldsymbol{\theta}} g_i^{\varphi}(p) \,(1-e^{-\kappa_i^{\boldsymbol{\theta}} x})
  + g_i^{\varphi}(p)\,\nabla_{\boldsymbol{\theta}} \kappa_i^{\boldsymbol{\theta}}\,x\,e^{-\kappa_i^{\boldsymbol{\theta}} x}
  \big], \label{eq:d-theta-v}\\
\nabla_{\boldsymbol{\theta}} v_x^{\boldsymbol{\theta}}(x,p)
&=\sum_{i=1}^{2}\big[
  \nabla_{\boldsymbol{\theta}} g_i^{\varphi}(p)\,\kappa_i^{\boldsymbol{\theta}}\,e^{-\kappa_i^{\boldsymbol{\theta}} x}
  + g_i^{\varphi}(p)\,\nabla_{\boldsymbol{\theta}} \kappa_i^{\boldsymbol{\theta}}\,e^{-\kappa_i^{\boldsymbol{\theta}} x}\,(1-\kappa_i^{\boldsymbol{\theta}} x)\big], \nonumber
\end{align}
with $\nabla_{\boldsymbol{\theta}} g_i^{\varphi}$ and $\nabla_{\boldsymbol{\theta}} \kappa_i^{\boldsymbol{\theta}}$ $(i=1,2)$ given in Proposition \ref{prop51para}.
Note that $\partial_{\gamma} g_i^{\varphi}(p) = 0$.
Then, $\boldsymbol{\theta}$ is updated by
$\boldsymbol{\theta}\leftarrow \boldsymbol{\theta} \, - \, \eta \, \odot \nabla_{\boldsymbol{\theta}}\mathrm{ML}_{\Delta t}^{\mathrm{reg}}(\boldsymbol{\theta}; \mathcal{D})$,
where $\eta:=\big(\eta_1, \eta_2, (\eta_{\gamma,0},\dots,\eta_{\gamma,4})\big)>0$ denotes the learning rate, and $\odot$ is the componentwise multiplication of learning rates and the gradients.
To save computations, we update each parameter only once for a given policy, rather than performing multiple updates to approach the true policy. This practice is common in actor-critic algorithms, see \cite{sutton2018}.

\begin{remark}
Estimating \eqref{defmlfunc0100} from a single on-policy trajectory may make the training process noisy. We therefore use mini-batches, averaging the estimate over multiple independent episodes. Each update is on-policy: after updating $\pi$, a new batch of finite-horizon trajectories $\mathcal{D}$ is sampled and used to further update the parameters.
\end{remark}
\textbf{Mode 2: Episodic Online CTD($\rho$).}
The second method imposes the martingale property through orthogonality: $M^{\boldsymbol{\theta}}$ is an $\mathbb{F}$-martingale if and only if, for any square-integrable $\mathbb{F}$-predictable test function $\xi$, the martingale orthogonality condition holds:
\begin{align}\label{orthos}
&\mathbb{E}\big[\int_0^t \xi_s \mathrm{d} M_s^{\boldsymbol{\theta}}\big] = 0.
\end{align}
Thus, approximating the optimizer $\boldsymbol{\theta}^*$ amounts to solving \eqref{orthos}.
Since \eqref{orthos} entails infinitely many conditions, we enforce it numerically only against a finite set of test functions, leading to a tractable finite system of moment conditions.
Although \eqref{opdipro} is intrinsically an infinite-horizon problem, truncation to $[0,T]$ yields a finite-horizon formulation in which \eqref{orthos} can be solved via the offline CTD methods of \cite{jia2022}, i.e., continuous analogues of the well-known discrete-time TD family, including CTD($\rho$) with $\rho\in[0, 1]$ and its least-squares variant CLSTD (particularly convenient for linear parametrizations).
To better match the infinite-horizon objective, we instead adopt an online CTD scheme that updates parameters along the realized trajectory.

By \eqref{martingavar}, the martingale increment associated with $\boldsymbol{\theta}$ is given by
\begin{align*}
\mathbb{E}[\mathrm{d}M_t^{\boldsymbol{\theta}}]
:= \mathrm{e}^{-\Lambda_t}\,
\mathbb{E}\big[
\mathrm{d}v^{\boldsymbol{\theta}}(\widehat X_t^{\boldsymbol{\theta}},p_t)
- \delta_t\,v^{\boldsymbol{\theta}}(\widehat X_t^{\boldsymbol{\theta}},p_t)\,\mathrm{d}t
+ \widehat{\mathcal H}_\lambda^{\boldsymbol{\theta}}(\widehat X_t^{\boldsymbol{\theta}},p_t)\,\mathrm{d}t
\big].
\end{align*}
In the discrete-time interaction scheme (Definition~\ref{discrete-time-scheme}),
\cite{jia2022} suggests $M_{t_{k+1}}^{\boldsymbol{\theta}} - M_{t_k}^{\boldsymbol{\theta}}
\approx e^{-\Lambda_{t_k}}\, \Delta_k^{\boldsymbol{\theta}}$, where $\Delta_k^{\boldsymbol{\theta}}$ is the discrete TD/martingale residual at time $t_k$, defined by
\begin{align}\label{disappmt}
&\Delta_k^{\boldsymbol{\theta}}:=\Delta v^{\boldsymbol{\theta}}(\widehat{X}_{t_k}^{\boldsymbol{\theta}}, p_{t_k}) - \delta_{t_k} v^{\boldsymbol{\theta}}(\widehat{X}_{t_k}^{\boldsymbol{\theta}}, p_{t_k}) \Delta t + \widehat{\mathcal{H}}_{\lambda}^{\boldsymbol{\theta}}(\widehat{X}_{t_k}^{\boldsymbol{\theta}}, p_{t_k}) \Delta t,
\end{align}
with $\Delta v^{\boldsymbol{\theta}}(\widehat{X}_{t_k}^{\boldsymbol{\theta}}, p_{t_k})
= v^{\boldsymbol{\theta}}(\widehat{X}_{t_{k+1}}^{\boldsymbol{\theta}}, p_{t_{k+1}})
-v^{\boldsymbol{\theta}}(\widehat{X}_{t_k}^{\boldsymbol{\theta}}, p_{t_k})$.
For episodic online CTD($\rho$), fix $\rho\in[0,1]$ and choose the test process (eligibility trace)
$\xi_t^{\boldsymbol{\theta}} := \int_0^t \rho^{t-s} \nabla_{\boldsymbol{\theta}}
v^{\boldsymbol{\theta}}(\widehat X_{s}^{\boldsymbol{\theta}}, p_{s})\mathrm{d}s$,
or in discrete time
\begin{align}\label{deltakindex}
&\xi_{t_k}^{\boldsymbol{\theta}}:=\sum_{i=0}^k \rho^{(k-i)\Delta t}\,\nabla_{\boldsymbol{\theta}}v^{\boldsymbol{\theta}}(\widehat X_{t_i}^{\boldsymbol{\theta}}, p_{t_i})\,\Delta t,
\end{align}
with the convention $0^z:=\mathbf{1}_{\{z=0\}}$ (so that $\xi_{t_k}^{\boldsymbol{\theta}}=\nabla_{\boldsymbol{\theta}} v^{\boldsymbol{\theta}}(\widehat X_{t_k}^{\boldsymbol{\theta}},p_{t_k}) \Delta t$ when $\rho=0$), where a continuous decay rate $\varrho>0$ corresponds to $\rho^{\Delta t}=e^{-\varrho\Delta t}$. Given an episode $\mathcal{D}=\{(t_k, \widehat{X}_{t_k}^{\boldsymbol{\theta}}, p_{t_k}, u_{t_k})\}_{k=0}^{K-1}$,
the CTD search direction takes the Robbins-Monro increment
\begin{align}\label{eq:G-TD}
G_{\mathrm{TD}}(\boldsymbol{\theta};\mathcal{D})
:= \sum_{k=0}^{K-1}
\mathrm{e}^{-\Lambda_{t_k}}\,
\xi_{t_k}^{\boldsymbol{\theta}}\,
\Delta_k^{\boldsymbol{\theta}}.
\end{align}
As in the regularized ML objective \eqref{defmlfunc0100}, we regularize the CTD direction \eqref{eq:G-TD} as
\begin{align}\label{eq:G-reg}
G_{\mathrm{TD}}^{\mathrm{reg}}(\boldsymbol{\theta};\mathcal{D})
&:= G_{\mathrm{TD}}(\boldsymbol{\theta};\mathcal{D})
- \Big(0,\big(w_j^{\mathrm{env}}(\mathrm{e}^{\gamma_j}-\mathrm{env}_j) \mathrm{e}^{\gamma_j}\big)_{j=0,\dots,4}\Big)\\
&\quad
- \Big(
\big(w_0^{\mathrm{bc}} e_0 \nabla_{\varphi_i}g_i^{\varphi}(0)
+ w_1^{\mathrm{bc}} e_1 \nabla_{\varphi_i}g_i^{\varphi}(1) \big)_{i=1,2},\;
0 \Big).\nonumber
\end{align}
With a learning-rate vector $\eta$, the episodic CTD($\rho$) update is $\boldsymbol{\theta}
\leftarrow \boldsymbol{\theta}
+ \eta \odot G_{TD}^{\mathrm{reg}}(\boldsymbol{\theta};\mathcal{D})$.
The choices $\rho=0$ and $\rho\approx1$ recover CTD(0) and near Monte-Carlo updates, respectively. In our implementation we specialize to CTD(0), so \eqref{eq:G-TD} reduces to a sum of one-step contributions along the episode.

Lemma 3 of \cite{jia2022} implies the following convergence result. Let $\Theta\subset\mathbb{R}^{d_\gamma}\times\mathbb{R}^{d_\varphi}$ be nonempty compact. Let $\mathcal{X}_\varepsilon(\Theta)$ be a class of $\mathbb{F}$-predictable processes $\xi=(\xi_t)_{t\in[0,T]}$ such that $\mathbb{E}[\int_0^{T} \xi_t^2 \mathrm{d} \langle M^{\boldsymbol{\theta}}\rangle_t]<\infty$ and $\mathbb{E}|\xi_t-\xi_s|^2 \le r(\boldsymbol{\theta}) |t-s|^{\varepsilon}$ for $\varepsilon\in(0,1)$ and $r\in\mathbb{C}(\Theta)$.
\begin{proposition}
Suppose that there exists $\boldsymbol{\theta}^*\in \Theta$ such that $\mathbb{E}[\int_{0}^{T}\xi_t\,\mathrm{d} M_t^{\boldsymbol{\theta}^*}]=0$ for all $\xi\in\mathcal{X}_\varepsilon(\Theta)$, and let $(\Delta t^{(n)})_{n\ge 0}\downarrow 0$, for each $n$ there exists $\boldsymbol{\theta}_{\Delta t^{(n)}}^*\in\Theta$ satisfying $\mathbb{E}[\sum_{k=0}^{K^{(n)}-1}\xi_{t_k^{(n)}}
(M_{t_{k+1}^{(n)}}^{\boldsymbol{\theta}_{\Delta t^{(n)}}^*}
-M_{t_k^{(n)}}^{\boldsymbol{\theta}_{\Delta t^{(n)}}^*})]=0$, and $(\boldsymbol{\theta}_{\Delta t^{(n)}}^*)_{n\ge0}$ converges in $\Theta$ as $n\to\infty$.
Then $\boldsymbol{\theta}_{\Delta t^{(n)}}^*\to \boldsymbol{\theta}^*$ as $n\to\infty$. Moreover, $|\mathbb{E}[\int_0^{T} \xi_t \mathrm{d} M_t^{\boldsymbol{\theta}_{\Delta t}^*}]|\le \hat{r} (\Delta t)^{\varepsilon/2}$ for some $\hat r>0$.
\end{proposition}

{
\SetAlFnt{\scriptsize}          
\SetAlCapFnt{\scriptsize}       
\SetAlCapNameFnt{\scriptsize}   
\SetAlgoNlRelativeSize{-3} 
\setlength{\algomargin}{0.3em}  
\SetInd{0.5em}{0.5em}           
\begin{algorithm*}[H]\label{Algorithm1}
\caption{PO-RSEOD Algorithm}
\KwIn{Market, initial surplus $X_0$, initial belief state $p_0$, initial market state $I_0$, fixed time horizon $T$, time step $\Delta t$, exploration weight $\lambda$, learning rates $\eta$ (ML/CTD), trace parameter $\rho\in[0,1]$, tolerance $\varepsilon$, number of iterations $N$
}
\textbf{Estimation of environment parameters:} {Obtain historical data to estimate $\mu_i$, $\sigma$ and $q_{ij}$ for $i,j=1,2$}\\
\KwOut{Estimated optimal dividend policy $\pi^{\boldsymbol{\theta}^*}$ and value function $v^{\boldsymbol{\theta}^*}$}
\textbf{Actor-critic Learning Procedure:} Initialization $\boldsymbol{\theta}^{(0)} = (\varphi^{(0)}, \gamma^{(0)})\in \Theta$
\BlankLine
\For{$n=1$ \KwTo $N$}{
  \uIf{\textbf{\textup{Mode 1:}} \textnormal{ML Minimization}}{
    \For{$k=0$ \KwTo $K^{(n)}-1$}
             {Sample episode
$\mathcal{D}^{(n)}=\{(t_k,\widehat X_{t_k}^{\boldsymbol{\theta}^{(n)}},p_{t_k},u_{t_k})\}_{k=0}^{K^{(n)}-1}$ from market under $\pi^{\boldsymbol{\theta}^{(n)}}$, where each $p_{t_k}$ is calculated by the
discretized Wonham filter using the estimated environment parameters, the stopping index
$K^{(n)}:=K_1 \wedge \min\{ k\ge0: \widehat X_{t_k}^{\boldsymbol{\theta}^{(n)}}<\varepsilon\}$ with
$K_1=\lfloor T/\Delta t\rfloor$ and $\varepsilon>0$ taken sufficiently small
}
\textbf{PE Update}\\
    Update $\boldsymbol{\theta}^{(n+1)}\leftarrow \boldsymbol{\theta}^{(n)} - \eta \odot \nabla_{\boldsymbol{\theta}^{(n)}}\mathrm{ML}_{\Delta t}^{\text{reg}}(\boldsymbol{\theta}^{(n)}; \mathcal{D}^{(n)})$ using \eqref{theta01}, the classical Robbins-Monro conditions that ensure convergence of the stochastic approximation under standard assumptions
  }
  \Else{
  \textbf{Mode 2: \textnormal{Episodic Online CTD(0)}}\\
  \For{$k=0$ \KwTo $K^{(n)}-1$}{
    Sample action $u_{t_k} \sim \pi(\cdot|\widehat X_{t_k}^{\boldsymbol{\theta}^{(n)}}, p_{t_k}) = G(u, 1-v_x^{\boldsymbol{\theta}^{(n)}}(\widehat{X}_{t_k}^{\boldsymbol{\theta}^{(n)}}, p_{t_k}))$,
    step environment to $(\widehat X_{t_{k+1}}^{\boldsymbol{\theta}^{(n)}}, p_{t_{k+1}})$, and continue until the stopping index $K^{(n)}$ ($K^{(n)}$ as noted in Mode 1) \\
    Compute gradient $\nabla_{\boldsymbol{\theta}^{(n)}}v^{\boldsymbol{\theta}^{(n)}}(\widehat X_{t_k}^{\boldsymbol{\theta}^{(n)}}, p_{t_k})$ in \eqref{eq:d-theta-v}, residual $\Delta_k^{\boldsymbol{\theta}^{(n)}}$ in \eqref{disappmt}, and eligibility trace $\xi_{t_k}^{\boldsymbol{\theta}^{(n)}}$ in \eqref{deltakindex} \\
  }
  \textbf{PE Update}\\
    Update $\boldsymbol{\theta}^{(n+1)} \leftarrow \boldsymbol{\theta}^{(n)}
+ \eta\odot G_{TD}^{\mathrm{reg}}(\boldsymbol{\theta}^{(n)};\mathcal{D}^{(n)})$ by \eqref{eq:G-reg}
  }
}
\end{algorithm*}
}

\section{Numerical analysis}\label{sec6}
In this section, we implement PO-RSEOD (Algorithm \ref{Algorithm1}) to estimate the unknown environment parameters and learn the optimal policy \eqref{densityop} along with the value function \eqref{eq:single_exp_ansatz}. As a convergence benchmark, we also compute these via an FD scheme, assuming known environment parameters. Some parameters are specified as follows: initial states $X_0=1$, $p_0=0.5$, fixed time horizon $T=10$ (years), time step $\Delta t=1/252$, tolerance $\varepsilon=10^{-8}$, discount rates $\delta_1= 0.1$, $\delta_2= 0.3$, and exploration weight $\lambda=1$.
For benchmarking we fix the environment at the ground-truth values $\mu_1= 1.2$, $\mu_2= 0.5$ (annual net premium inflow rates), and $(q_{12},q_{21})=(0.36,2.89)$ per year (from \cite{wu2024}).

\subsection{FD and polynomial approximation}\label{subsec610}

The endpoint splits $\varpi^{(j)}=(\varpi^{(j)}_1,\varpi^{(j)}_2)\in[0, 1]^2$, $j \in\{0,1\}$ in \eqref{addicondi}, can be calibrated by matching the endpoint-derivative targets: let $\nu(\cdot,i)$ be the completely observed value function (see Section \ref{sec23}), we enforce
\begin{equation*}
v_x(0,0)=g(0)(\varpi^{(0)}_1\kappa_1+\varpi^{(0)}_2\kappa_2)\approx \nu_x(0,2), \
v_x(0,1)=g(1)(\varpi^{(1)}_1\kappa_1+\varpi^{(1)}_2\kappa_2)\approx \nu_x(0,1).
\end{equation*}
Define the normalized targets $s^{(0)}:=\tfrac{\nu_x(0,2)}{g(0)}$, $s^{(1)}:=\tfrac{\nu_x(0,1)}{g(1)}$, and the residual
$r^{(j)}(\kappa,\varpi^{(j)}):=s^{(j)}-(\varpi^{(j)}_1\kappa_1+\varpi^{(j)}_2\kappa_2)$.  Then, $\varpi^{(j)}$ are determined by the following outer fixed-point iteration.
\begin{enumerate}
\item[(i)] \emph{Initialization.} Choose admissible initial splits $\varpi^{(j)}\in[0,1]^2$ with $\varpi^{(j)}_1+\varpi^{(j)}_2=1$.
\item[(ii)] \emph{State update.} With the current splits $\varpi^{(j)}$, solve the boundary-value problems
\eqref{eq:V_const_eq00}-\eqref{eq:V_exp_eq} using the FD scheme in Appendix \ref{finitediff00} to obtain $g_i$ and update $\kappa_i$.
\item[(iii)] \emph{Split update.}
For fixed $\kappa=(\kappa_1,\kappa_2)$, updating $\varpi^{(j)}$ amounts to the least-squares problem $\min_{\varpi^{(j)}}\ \big(r^{(j)}(\kappa,\varpi^{(j)})\big)^2$. Its solution is $\varpi^{(j),\star}_1
=\operatorname{proj}_{[0,1]}\!\big(\tfrac{s^{(j)}-\kappa_2}{\kappa_1-\kappa_2}\big)$ and $\varpi^{(j),\star}_2=1-\varpi^{(j),\star}_1$, where $\operatorname{proj}_{[0,1]}$ denotes projection onto $[0,1]$. In the nearly degenerate case $|\kappa_1-\kappa_2|\ll 1$, set $\varpi^{(j),\star}_1 = 1/2$ for stability.
Then, we perform a relaxed update with $\rho_{re}\in(0,1]$: $\varpi^{(j)}\leftarrow (1-\rho_{re}) \varpi^{(j)}+\rho_{re} \varpi^{(j),\star}$.
\item[(iv)]
\emph{Stopping.} Repeat (ii)-(iii) until
$|r^{(j)}(\kappa,\varpi^{(j)})|\le \varepsilon_{\mathrm{ep}}$ and
$\|\varpi^{(j)}_{\text{new}}-\varpi^{(j)}\|_\infty\le \varepsilon_{\varpi}$.
\end{enumerate}

\begin{figure}[htbp]
	\centering
	\includegraphics[width=3in]{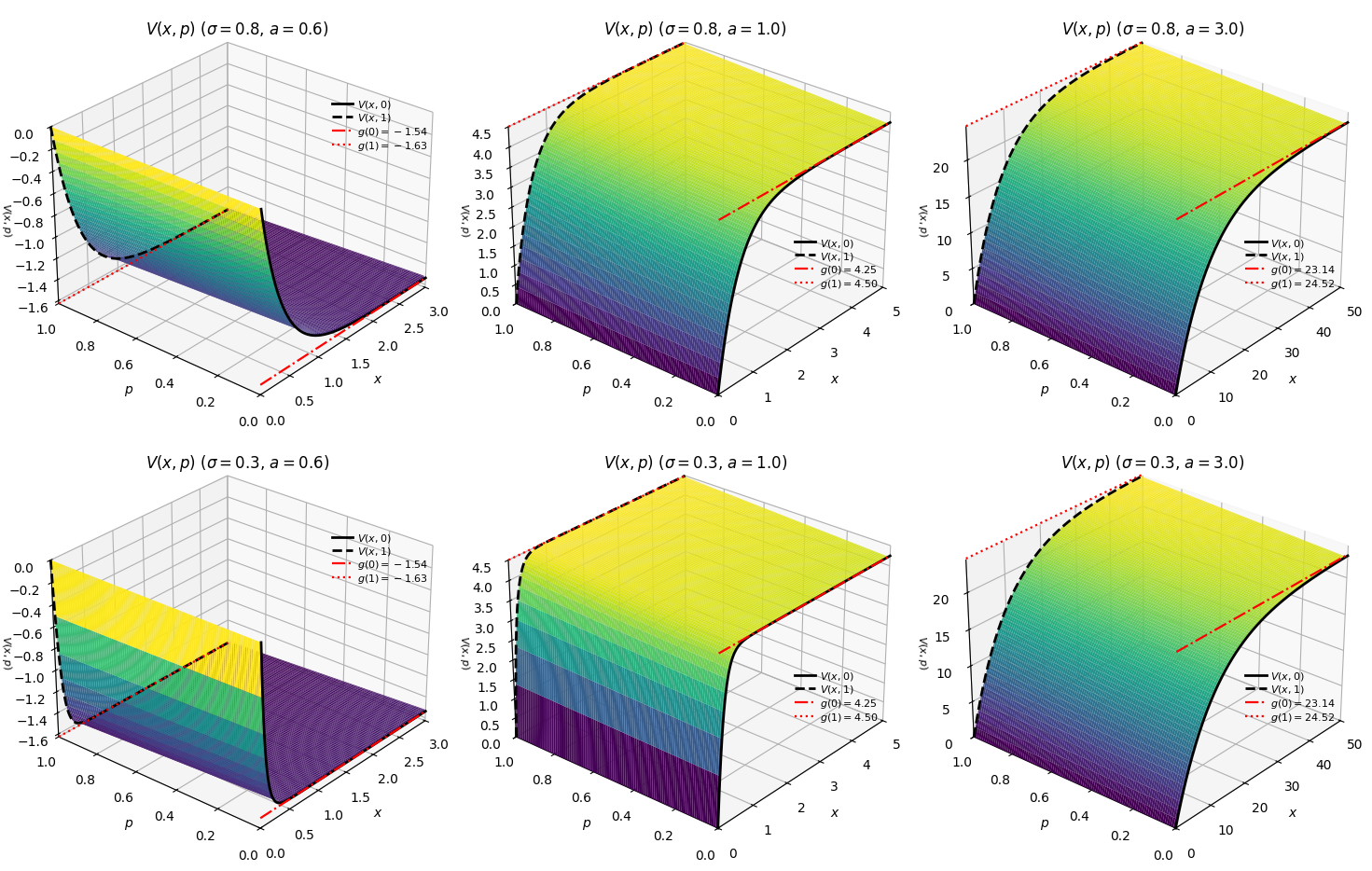}
	\caption{Surfaces of $V(x, p)$ from FD for dividend caps $\mathbf{a}\in\{0.6, 1, 3\}$ and volatilities $\sigma\in\{0.3, 0.8\}$. }
	\label{fig:1}
\end{figure}

Fix the true environment parameters and set $\rho_{re}=0.01$, $\varepsilon_{\mathrm{ep}} = \varepsilon_{\varpi} = 10^{-12}$ and $(\nu_x(0,1), \nu_x(0, 2))=(1.2, 1.6)$. We first calibrate the endpoint splits $\varpi^{(j)}$ via steps (i)-(iv). We then solve the boundary-value problems \eqref{eq:V_const_eq00}-\eqref{eq:V_exp_eq} by the FD scheme to obtain FD approximations of $g_i(p)$ and $\kappa_i$ $(i=1,2)$. Substituting these into \eqref{densityop} and \eqref{eq:single_exp_ansatz} yields $\pi^*(u|x, p)$ and $V(x,p)$ for problem \eqref{opdipro}. Fig.\ref{fig:1} plots the resulting FD approximation of $V(x, p)$ for dividend caps $\mathbf{a}\in\{0.6, 1, 3\}$ and volatilities $\sigma\in\{0.3, 0.8\}$.

From Fig.\ref{fig:1}, we first observe that for fixed $p$, $V(x,p)$ is monotone in $x$: it decreases when $\mathbf{a}=0.6$ and increases when $\mathbf{a}\in\{1,3\}$. This accords with our analytic characterization: for $\lambda=1$, the critical level $\mathbf{a}_c:=\lambda\ln(1+1/\lambda)=\ln 2\approx 0.6931$ determines the sign of $f_{\lambda}(0)$ (Prop.~\ref{propfuncf00}(iv)), if $\mathbf{a}>\mathbf{a}_c$, then $f_{\lambda}(0)>0$ and $g_i(p)\ge 0$ (Prop.~\ref{theosoluodse}(ii)), implying $\partial_x V(x,p)=\sum_{i=1}^2 \kappa_i g_i(p) e^{-\kappa_i x}\ge 0$ and hence $V(x,p)$ is nondecreasing in $x$; conversely, if $0<\mathbf{a}<\mathbf{a}_c$, then $V(x,p)$ is nonincreasing in $x$.
Second, for fixed $x$, $V(x,p)$ varies only marginally with $p$.
Third, $V(x, p)$ flattens as $x\to \infty$, with $p\in\{0,1\}$-edge slices converging to $g(0)$ and $g(1)$, respectively; these endpoint limits increase with $\mathbf{a}$. Finally, larger $\sigma$ slows the approach of $V$ to its asymptote $g(p)$, whereas smaller $\sigma$ leads to earlier saturation. Economically, higher volatility amplifies downside risk and weakens the marginal value of additional surplus, thereby delaying convergence in $x$.

\begin{figure}[htbp]
	\centering
	\includegraphics[width=3in]{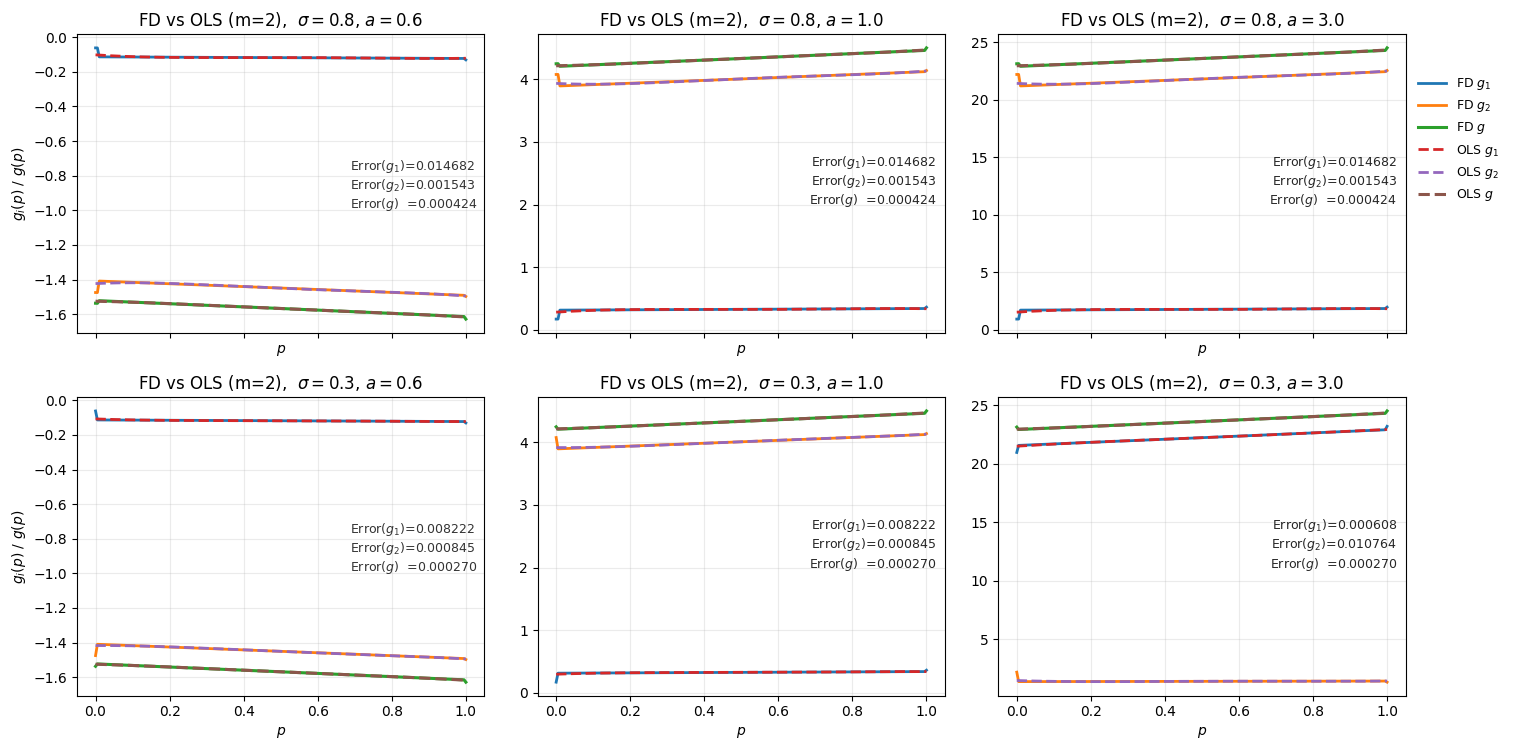}
	\caption{$g_i(p)$ $(i=1,2)$ and $g(p)$ from FD versus their OLS polynomial fits with order $m=2$.}
	\label{fig:2}
\end{figure}

In our RL algorithm, we parameterize $g_i(p)$ as in \eqref{paragi0} and choose the polynomial degree $m$ by matching FD benchmarks. Concretely, we set $m=2$. Fig.\ref{fig:2} contrasts the FD benchmarks for $g_i(p)$ $(i=1,2)$ and $g(p)$ with their quadratic approximation fitted by ordinary least squares (OLS). The (annotated) OLS fit errors, measured by the mean absolute error normalized by the mean absolute value of the benchmark curve, are all below 1.5\%, and the OLS and FD curves are virtually indistinguishable on the interior of $[0,1]$, up to minor endpoint effects. We therefore adopt $m=2$ throughout.

\subsection{Environment estimation using two methods} \label{sec62}
Our RL algorithm requires estimating the environment before interacting with it for filtering purposes.
Following the heuristic method of \cite{dai2010}, we classify regimes using an additive window-threshold rule with lookback $M=252$: for each $t_k$, define the window change $\Delta^{(M)}X_{t_k}:=X_{t_k}-X_{t_{b(k)}}$ with $b(k):=\max\{0,k-M\}$. Given thresholds $U_{+},U_{-}>0$, we label $I_{t_k} = 1$ (bull) if $\Delta^{(M)}X_{t_k}\ge U_{+}$, and $I_{t_k} = 2$ (bear) if $\Delta^{(M)}X_{t_k}\le - U_{-}$. Instead of the fixed percentage thresholds in \cite{dai2010} ($U_+=24\%$, $U_-=19\%$), we adopt the adaptive choice $U_{+}=U_{-}=\eta_u \hat{\sigma}_{0}\sqrt{M \Delta t}$, where $\hat{\sigma}_{0}$ is a robust pilot scale (MAD-based) estimator $\mathrm{MAD}(\Delta X_{t_k})/\sqrt{\Delta t}$ and $\eta_u =0.15$ empirically. Conditional on the inferred regimes, we estimate the environment parameters from the increments $\Delta X_{t_k}:=X_{t_{k+1}}-X_{t_k}$, the regime-dependent drift is estimated by $\mu_i =\frac{\sum_{k=0}^{K-1}\mathbf{1}\{I_{t_k}=i\}\,\Delta X_{t_k}}{\Delta t\,\sum_{k=0}^{K-1}\mathbf{1}\{I_{t_k}=i\}}$,
the diffusion coefficient by the annualized residual variance $\sigma^2=(K\Delta t)^{-1}\sum_{k=0}^{K-1}(\Delta X_{t_k}-\mu_{I_{t_k}}\Delta t)^2$,
and the transition rates $q_{ij}$ $(i\neq j$) by the reciprocal of the average (yearly) holding time in regime $i$, with $q_{ii}=-q_{ij}$.
For comparison, parameters are also estimated via a standard EM algorithm (see e.g., \cite{elliott1995}).

We simulate 100 (independent) twenty-year sample paths under the true parameters. On each path, we apply both methods to obtain parameter estimates. Table \ref{Table 1} reports the averaged estimates across these 100 paths.
The results indicate that the heuristic method produces reasonably accurate estimates of $\mu_i$ $(i=1,2)$ and $\sigma$. In comparison, the EM algorithm tends to overestimate both drifts and to distort the transition rates $q_{12}$ and $q_{21}$, leading to spuriously frequent regime switches. The mechanism is clear: with daily increments dominated by variance (low signal-to-noise), the EM algorithm may identify regime switches mainly from noises, resulting in biased estimates. Similar failures have been documented in \cite{sass2004} and \cite{wu2024}.
Overall, the heuristic method delivers more reliable transition-rate estimates, despite a modest upward bias in $q_{12}$, and we therefore adopt it for estimating the environment parameters.

\begin{table}[!ht]
\caption{Comparison between the heuristic estimation method and the EM algorithm.}
\label{Table 1}
\centering
\tabstyle
\resizebox{\linewidth}{!}{%
\begin{tabular}{l| l l l l l}
\toprule
$\text{Parameters}$ & $\mu_1$ & $\mu_2$ & $\sigma$ & $q_{12}$ & $q_{21}$ \\
\cline{1-6}
$\text{True}$ & $1.200$ & $0.500$ & $0.300$ & $0.360$ & $2.890$\\
$\text{Heuristic}$ & $1.130 \ (\text{sd:} 0.073)$ & $0.190\ (\text{sd:} 0.419)$ & $0.300\ (\text{sd:} 0.003)$ & $1.002\ (\text{sd:} 0.005)$ & $3.070\ (\text{sd:} 6.782)$\\
$\text{EM}$ & $1.328\ (\text{sd:}0.198) $ & $0.734\ (\text{sd:} 0.401)$ & $0.300\ (\text{sd:} 0.003)$ & $1.756\ (\text{sd:} 1.621)$ & $3.592\ (\text{sd:} 3.211)$ \\
\bottomrule
\multicolumn{6}{l}{\tiny Note: `sd' denotes the sample standard deviation of the estimates.}
\end{tabular}%
}
\end{table}

\begin{figure}[htbp]
	\centering
	\includegraphics[width=3in]{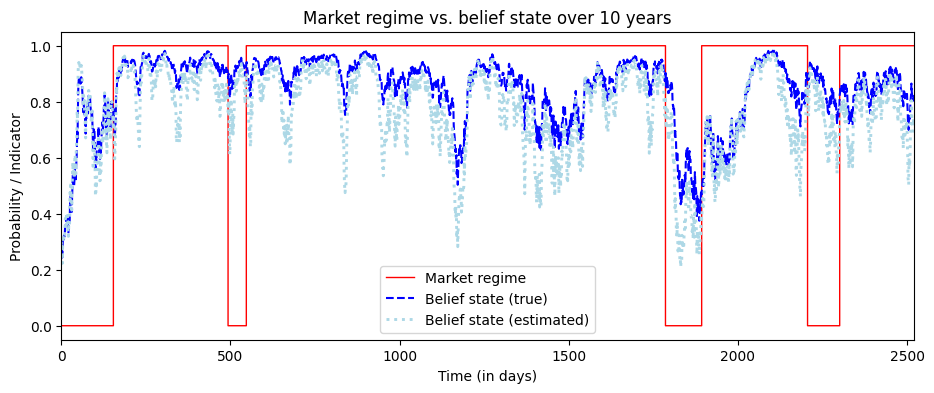}
	\caption{Market regime vs. belief states constructed with true and estimated parameters over 10 years.}
	\label{figsurpath}
\end{figure}

We now consider two filtering settings. In the true setting, the discretized Wonham filter is run under the true environment parameters to produce the belief process. In the estimated setting, the filter uses the heuristic estimates.
Fig.\ref{figsurpath} compares these belief trajectories.
The estimated belief closely tracks that under the true parameters and shows that parameter estimation error has only a limited effect on filtering accuracy.

\subsection{Convergence of PO-RSEOD}\label{sec63}
We analyze the convergence of PO-RSEOD, with particular attention to $\gamma_i$ for $i=0, \ldots, 4$. Comparing \eqref{parawgamp0} and \eqref{eq:w-closed-2}, we see that if
\begin{align*}
&e^{\gamma_0} \to e^{\gamma_0^*}:= \sigma^2, \quad
e^{\gamma_1} \to e^{\gamma_1^*}:= \mu_1, \quad
e^{\gamma_2} \to e^{\gamma_2^*}:= \mu_2,\\
&e^{\gamma_3} \to e^{\gamma_3^*}:= q_{21}, \quad
e^{\gamma_4} \to e^{\gamma_4^*}:= q_{12},
\end{align*}
then PO-RSEOD learns the optimal policy $\pi^*(u|x, p)$ for problem \eqref{opdipro}.
Initialize $\varphi = (\varphi_1,\varphi_2)$ with $\varphi_1,\varphi_2\in\mathbb{R}^9$ set to the constant vector $-3\cdot \mathbf{1}_9$, and initialize $\gamma=(\ln 0.07,0,0,0,0)$.
The base learning rates are $\eta= (3\times 10^{-4}, 3\times 10^{-4}, \eta_{\gamma})$ with $\eta_{\gamma} = (3\times 10^{-2}, 5\times 10^{-3}, 5\times 10^{-3}, 5\times 10^{-3}, 5\times 10^{-3})$, with a power-law decay of exponent $0.1$ applied during training (see \cite{bottou2018}), where $\varphi$ is updated with a smaller step size due to its higher sensitivity to data variability.
We empirically set the regularization weights to $w^{\mathrm{env}}=(7.0,0.5,0.5,0.2,0.2)$ and $w^{\mathrm{bc}}=(60,60)$, and fix $\sigma =0.3$ (true) and $\mathbf{a}=1$.

\begin{figure}[htbp]
	\centering
	\includegraphics[width=2.4in]{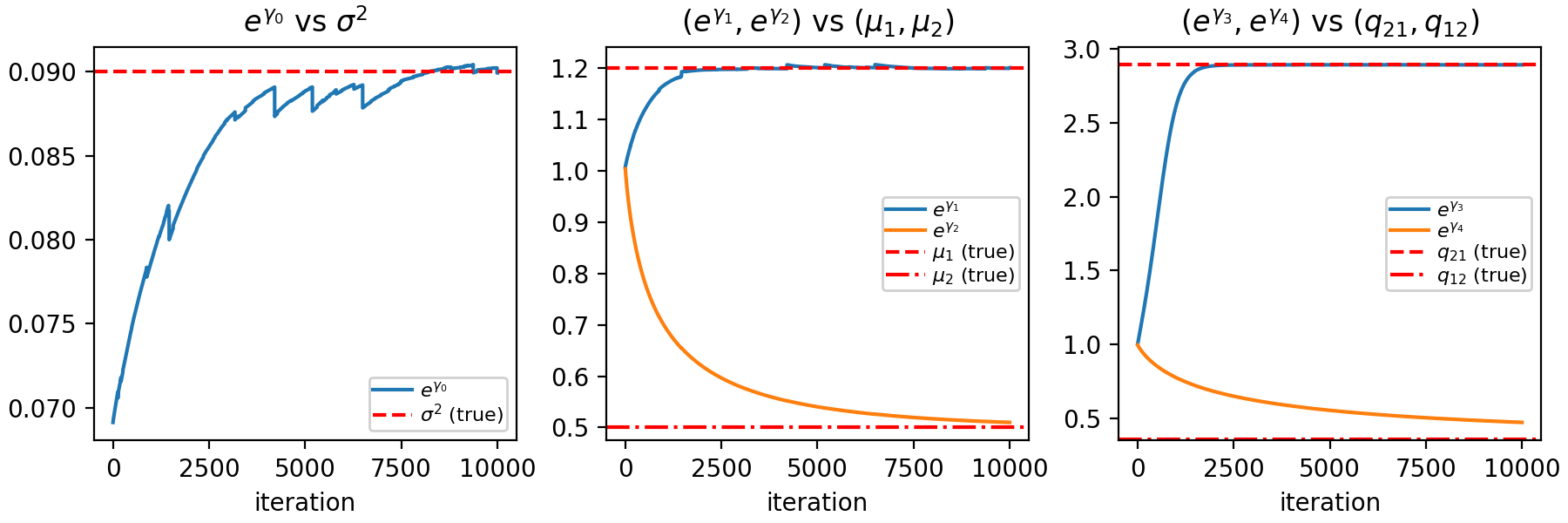} \ \ \ \ \ \ \ \
\includegraphics[width=1.2in]{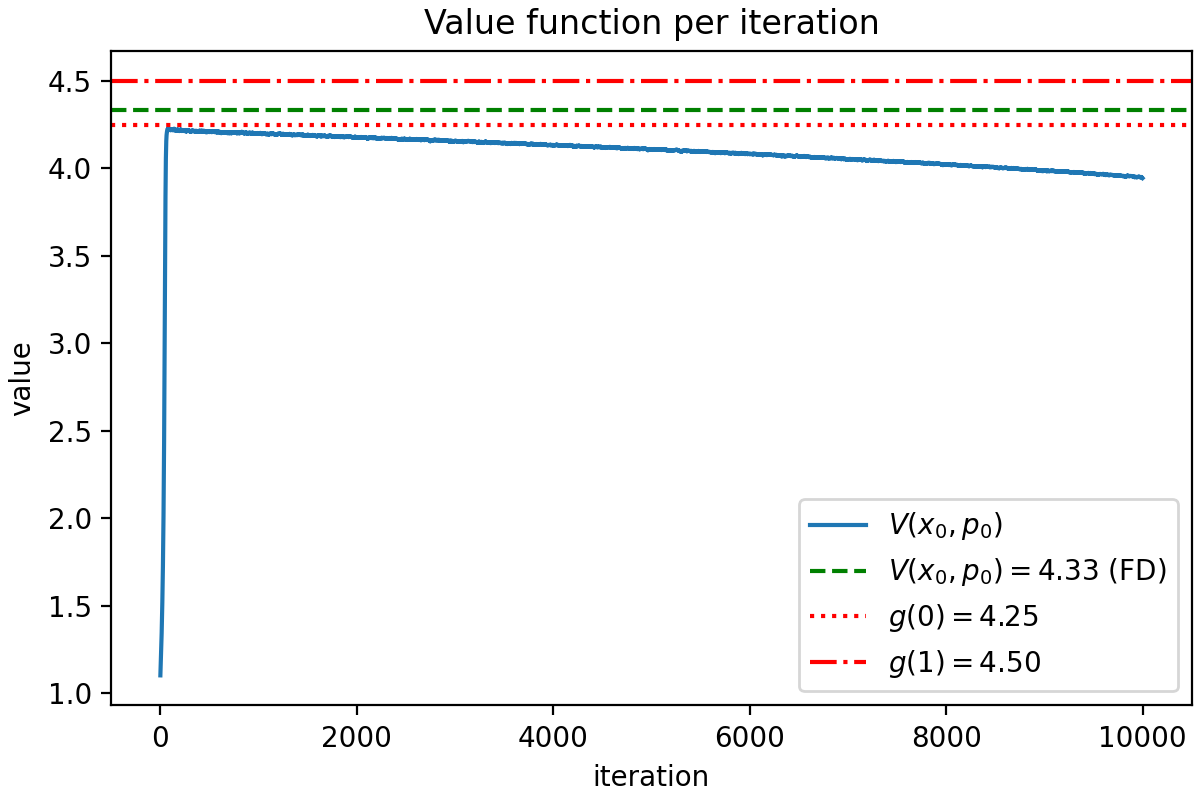}
\includegraphics[width=1.2in]{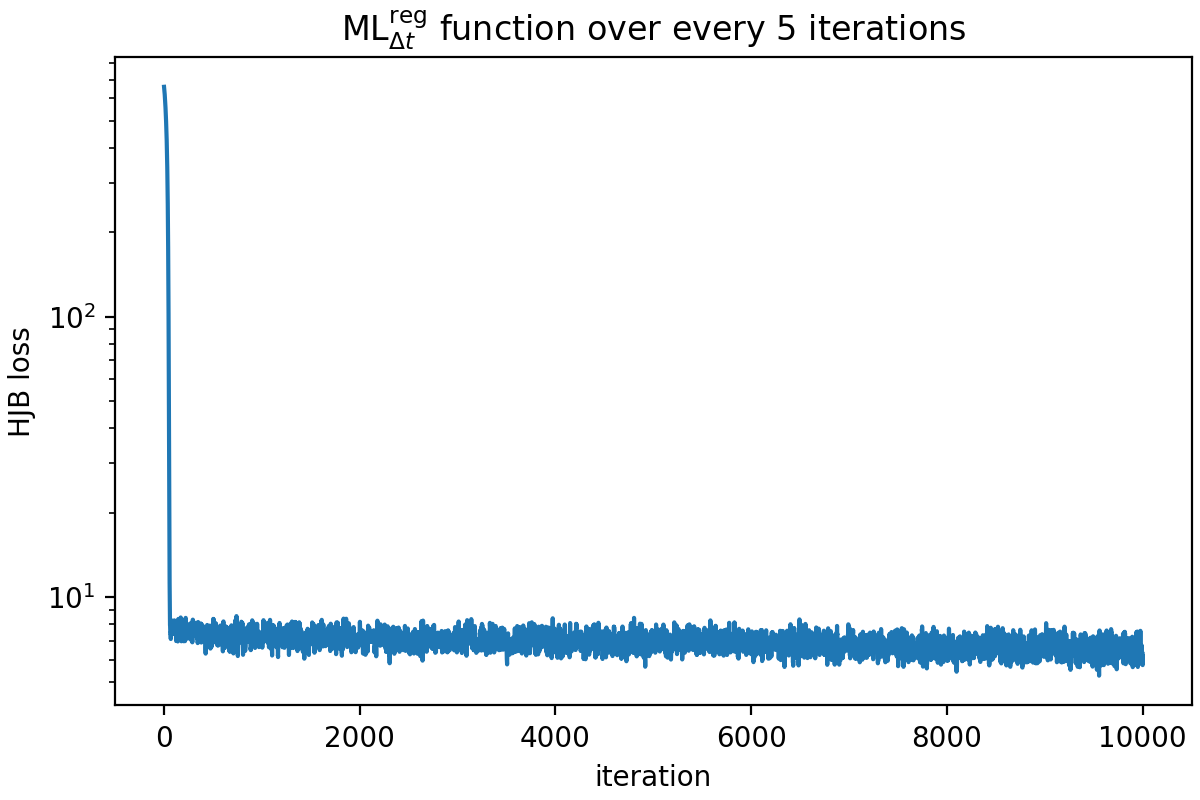} \\
\subfloat{\scriptsize{(1) ML Minimization (est)}}\\
\includegraphics[width=2.4in]{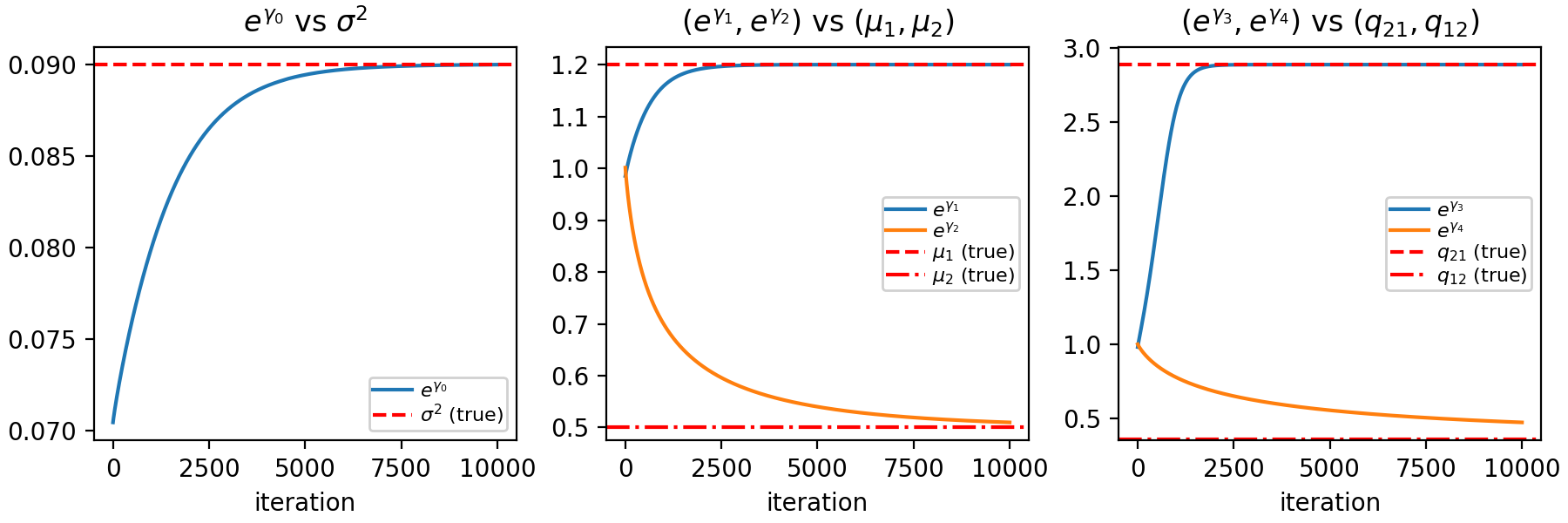} \ \ \ \ \ \ \ \
\includegraphics[width=1.2in]{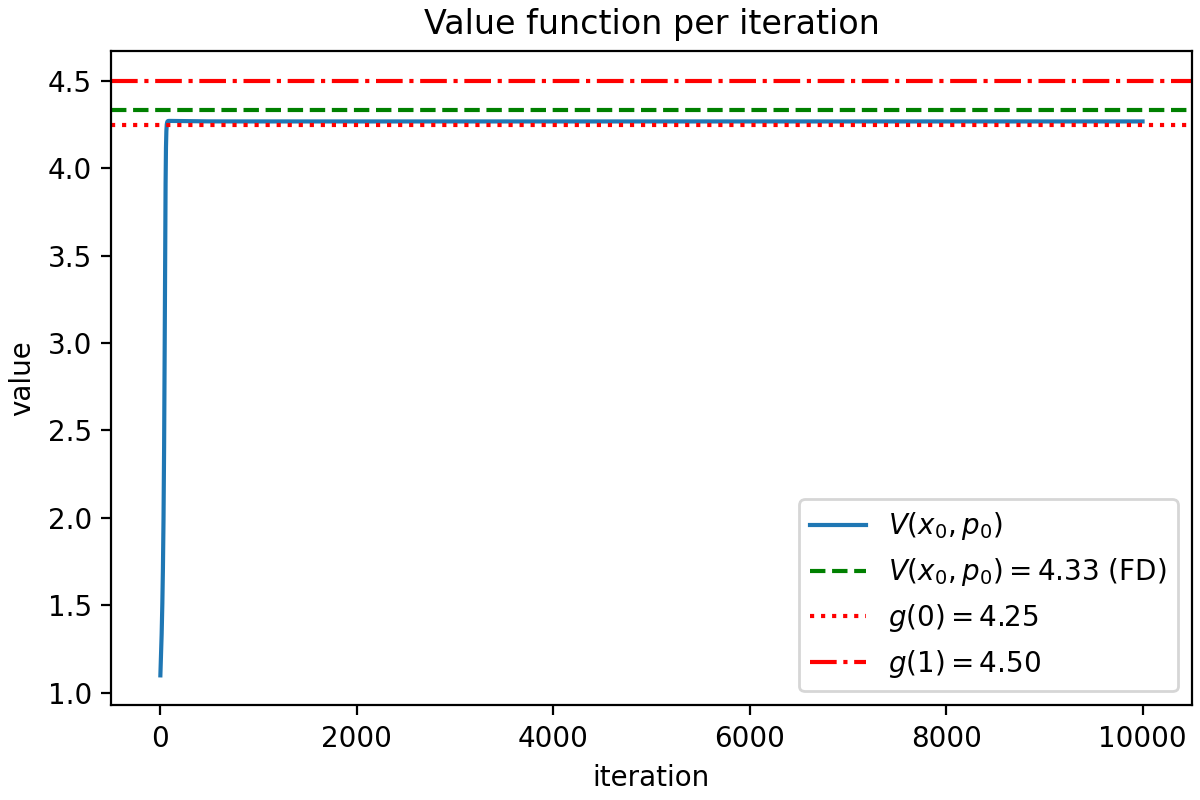}
\includegraphics[width=1.2in]{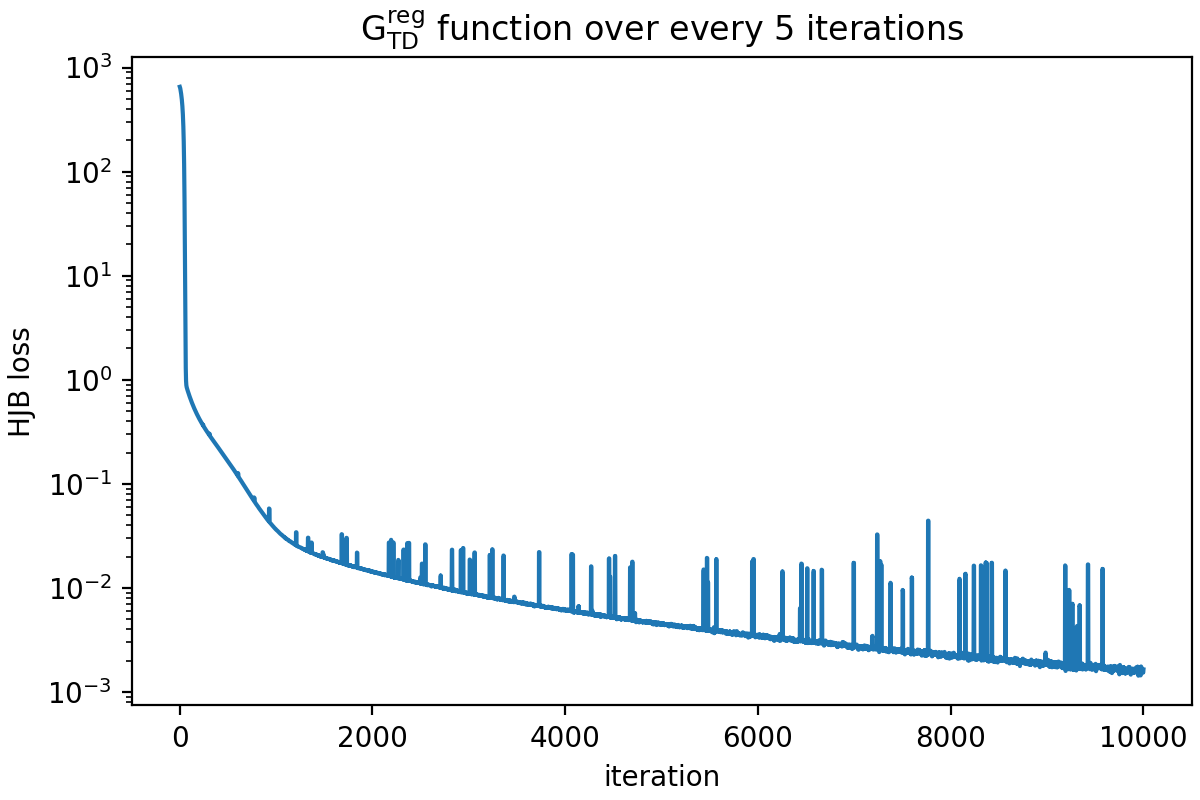}\\
\subfloat{\scriptsize{(2) Episodic Online CTD(0) (est)}}\\
\includegraphics[width=2.4in]{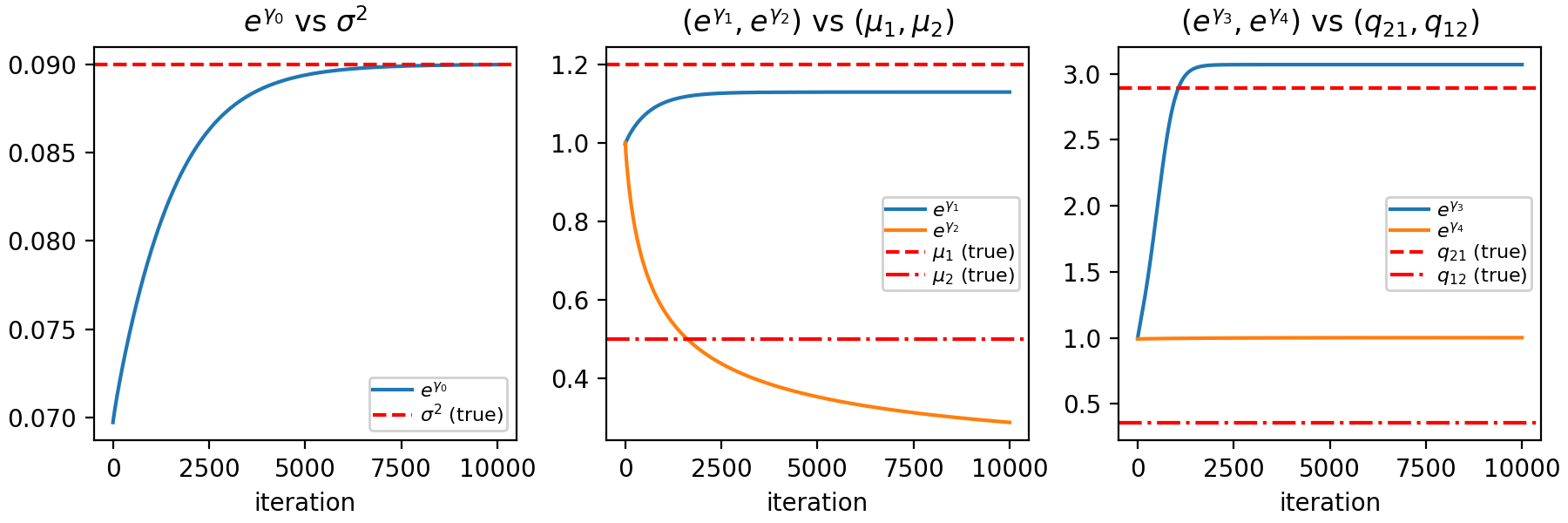} \ \ \ \ \ \ \ \
\includegraphics[width=1.2in]{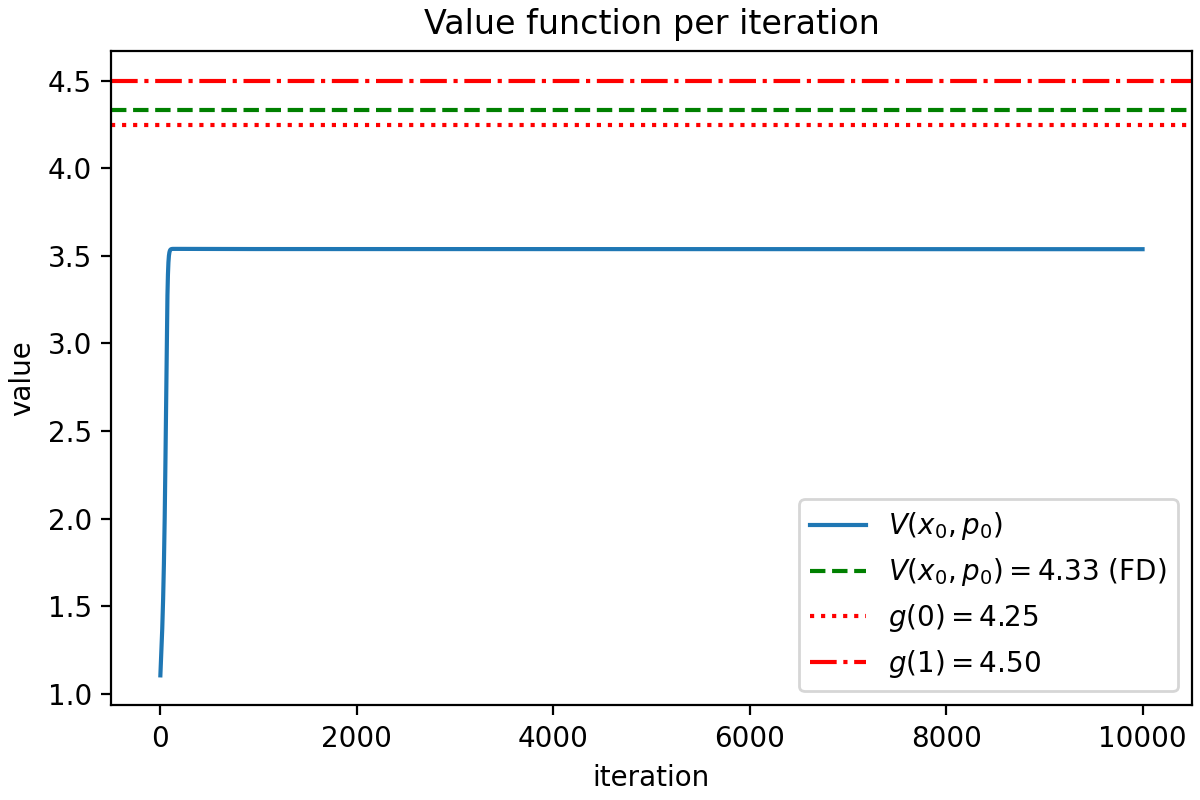}
\includegraphics[width=1.2in]{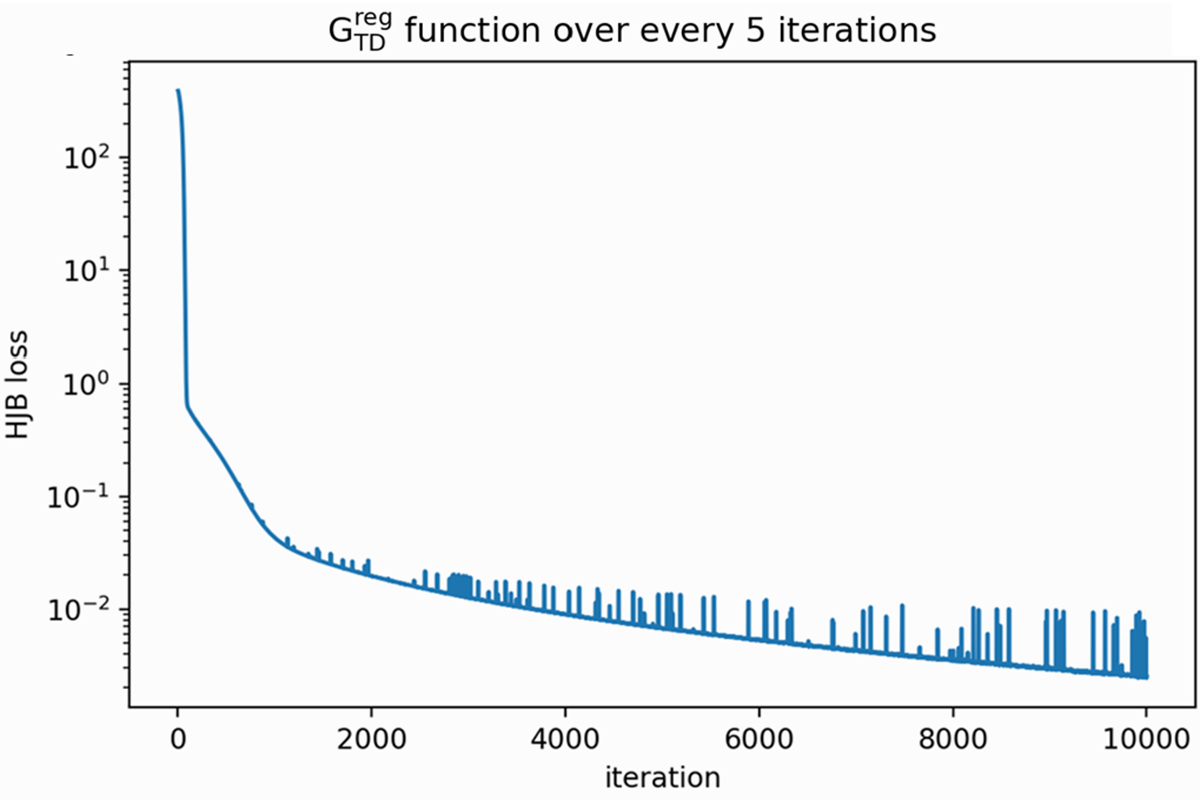}\\
\subfloat{\scriptsize{(3) Episodic Online $\mathrm{CTD}^*$(0) (est)}}\\
	\caption{\textbf{Left:} training convergence of $e^{\gamma_i}$, $i=0,\ldots, 4$ for ML Minimization/Episodic Online CTD(0)/ $\mathrm{CTD}^*(0)$; \textbf{Right:} training paths of value function and regularized HJB loss (averaged over every 5 iterations).}
	\label{figgamma}
\end{figure}

For training, in each iteration, we simulate a 30-year surplus path of \eqref{sde1} under the true environment parameters in Table \ref{Table 1}, where the first 20 years of data are used for estimating the environment parameters, and the last 10 years of data are used to generate a single on-policy episode via filtering under the resulting estimates. We train PO-RSEOD with a total of $N=10,000$ iterations, which corresponds to 10,000 episodes (batche=1), and compare three training schemes: (1) ML Minimization with filtering under the estimated parameters and regularization toward the true parameters; (2) Episodic Online CTD(0) under the same filtering-regularization configuration; (3) Episodic Online CTD(0) with filtering under the estimated parameters and regularization toward the heuristic estimates in Table \ref{Table 1}, denoted $\mathrm{CTD}^*$(0). In what follows, the tags (est) and (true) denote filtering under the estimated and true parameters, respectively. The learning processes of $e^{\gamma_i}$, $i=0, 1, \ldots, 4$ under the three training modes are shown in Fig.\ref{figgamma}(left). The corresponding estimates of the value function $V(X_0, p_0)$ and the moving-average HJB loss are displayed in Fig.\ref{figgamma}(right).

As shown in Fig.\ref{figgamma}(left), when the regularization targets the true parameters, our actor-critic learning converges close to $\gamma_i^*$ $(i=0,\ldots,4)$ under both ML and CTD(0). This behavior is essentially the same under (true) and (est) filtering, so we omit the (true) results for brevity. In contrast, when regularization is imposed toward the heuristic estimates ($\mathrm{CTD}^*(0)$), $\gamma_0$ converges accurately to $\gamma_0^*$, whereas $\gamma_i$ for $i=1, \ldots, 4$ converge to noticeably biased limits relative to $\gamma_i^*$. This indicates that parameter estimation error, propagating through the regularization target, can affect policy learning, even though Fig.\ref{figsurpath} suggests that its impact on filtering accuracy is limited.
Fig.\ref{figgamma}(right) highlights that CTD(0) effectively bring the value estimate $V$ to approach the FD benchmark within finitely many iterations, consistent with the policy convergence result (Section \ref{subsec 5.1}). By comparison, ML attains a good early fit but then drifts downward, while $\mathrm{CTD}^*(0)$ converges to a value with a significant gap relative to the FD benchmark, reflecting heightened sensitivity to the regularization target. Across all three modes the HJB residual decreases, but the steady-state magnitudes differ substantially: CTD(0) and $\mathrm{CTD}^*(0)$ stabilize around $10^{-2}$, whereas ML remains around $10^1$, pointing to much tighter HJB consistency under CTD-style updates.

Table~\ref{Table 2} reports the learned values of $e^{\gamma_i}$ $(i=0,\ldots,4)$ after 10,000 iterations, along with the heuristic environment estimates averaged over 10,000 simulated paths. Row 1 gives the benchmark under the true parameters, rows 2-3 correspond to ML minimization with (true) and (est) filtering, respectively, while rows 4-6 summarize episodic online CTD(0) and $\mathrm{CTD}^*(0)$.
It shows that regularization toward the true parameters (ML and CTD(0)) drives the learned $e^{\gamma_i}$ close to their targets, with only a slight upward bias in $e^{\gamma_4}$. By contrast, regularization toward the heuristic estimates ($\mathrm{CTD}^*(0)$) propagates estimation error into the learned limits: $e^{\gamma_0}$ remains accurate, $e^{\gamma_1}$-$e^{\gamma_2}$ are biased downward, and $e^{\gamma_3}$-$e^{\gamma_4}$ are biased upward. The parameter estimates reported in the right panel exhibit the same bias pattern as in Table \ref{Table 1}.

\begin{table}[!ht]
\caption{Convergence results for parameter $e^{\gamma}$ and the estimated environment parameters during training.}
\label{Table 2}
\centering
\tabstyle
\scriptsize
\resizebox{\linewidth}{!}{%
\begin{tabular}{l| c c c c c| c c c c c }
\toprule
$\ $ & $e^{\gamma_0}$ & $e^{\gamma_1}$ & $e^{\gamma_2}$ & $e^{\gamma_3}$ & $e^{\gamma_4}$ & ${\sigma^2}_{\text{est}}$ & ${\mu_{1}}_{\text{est}}$ & ${\mu_{2}}_{\text{est}}$ & ${q_{21}}_{\text{est}}$ & ${q_{12}}_{\text{est}}$ \\
\cline{1-11}
$\text{Optimal}$ & $0.0900$ & $1.2000$ & $0.5000$ & $2.8900$ & $0.3600$
& $0.0900$ & $1.2000$ & $0.5000$ & $2.8900$ & $0.3600$\\
\midrule
$\text{ML (true)}$ & $0.0901$ & $1.1996$ & $0.5098$ & $2.9988$ & $0.4755$
& $0.0900$ & $1.2000$ & $0.5000$ & $2.8900$ & $0.3600$\\
$\text{ML (est)}$ & $0.0899$ & $1.2006$ & $0.5100$ & $2.8902$ & $0.4750$
& $0.0901$ & $1.1234$ & $0.5398$ & $4.8114$ & $1.0026$\\
\midrule
$\text{CTD(0) (true)}$ & $0.0900$ & $1.2000$ & $0.5100$ & $2.8900$ & $0.4749$
& $0.0900$ & $1.2000$ & $0.5000$ & $2.8900$ & $0.3600$ \\
$\text{CTD(0) (est)}$ & $0.0900$ & $1.2000$ & $0.5099$ & $2.8900$ & $0.4753$
& $0.0901$ & $1.1233$ & $0.5265$ & $4.6918$ & $1.0027$\\
$\mathrm{CTD}^*(0) \text{(est)}$ & $0.0900$ & $1.1300$ & $0.2876$ & $3.0700$ & $1.0019$
& $0.0901$ & $1.1219$ & $0.5507$ & $4.8916$ & $1.0026$ \\
\bottomrule
\end{tabular}
}
\end{table}

We further assess the out-of-sample performance of PO-RSEOD under both the true and the estimated environments, using the parameter sets in Table \ref{Table 2}. For each setting, we generate 100,000 10-year testing paths from the model \eqref{sde1} that are independent of those used for training, and implement the learned policies on them.
As a benchmark, we also implement the optimal policy obtained by the FD scheme under the true parameters.
Performance is summarized by five criteria: the mean and variance of the estimated value $V(X_0, p_0)$, the signal-to-noise ratio $\mathrm{SNR}:=\mathbb{E}[V]/\sqrt{\text{Var}(V)}$, and two annualized Sharpe-type ratios computed from per-step series: Sharpe-SR uses surplus returns $r_k = (\widehat{X}_{t_{k+1}} - \widehat{X}_{t_k})/ \widehat{X}_{t_k}$, while Sharpe-RI uses discounted reward increments $\Delta V_{t_k} = e^{-\Lambda_{t_k}} \widehat{\mathcal{H}}_{\lambda}(\widehat{X}_{t_k}, p_{t_k}) \Delta t$. Each ratio is the Monte Carlo mean of the corresponding series divided by its Monte Carlo standard deviation, annualized by $\sqrt{252}$.
The resulting estimates based on 100,000 simulated paths are displayed in Table \ref{Table 3}.

\begin{table}[!ht]
\caption{Out-of-sample performance.}
\label{Table 3}
\centering
\footnotesize
\begin{tabular}{l| c c c c c }
\toprule
$\ $ & $\text{Mean}$ & $\text{Variance}$ & $\text{SNR}$ & $\text{Sharpe-SR}$ & $\text{Sharpe-RI}$\\
\cline{1-6}
$\text{Optimal}$ & $4.356172$ & $0.018554$ & $31.980767$ & $1.462071$ & $45.254949$\\
\midrule
$\text{ML (true)}$ & $4.305606$ & $0.019616$ & $30.741601$ & $1.407813$ & $45.256492$\\
$\text{ML (est)}$ & $4.083415$ & $0.007700$ & $46.536182$ & $1.367718$ & $42.129792$ \\
\midrule
$\text{CTD(0) (true)}$ & $4.354974$ & $0.019281$ & $31.363528$ & $1.468726$ & $45.253587$ \\
$\text{CTD(0) (est)}$ & $4.128010$ & $0.007929$ & $46.358537$ & $1.471748$ & $41.981312$ \\
$\mathrm{CTD}^*(0)\text{(est)}$ & $3.943364$ & $0.006076$ & $50.587958$ & $1.484987$ & $42.233154$ \\
\bottomrule
\end{tabular}
\end{table}

\normalsize

Table \ref{Table 3} shows that under (true) filtering, PO-RSEOD attains out-of-sample moments of $V$ that are very close to the optimal benchmark. The (est) setting also performs well and, consistent with the stronger effective discounting implied by (est) filtering,
\footnote{{\tiny{Estimated-environment filtering produces a downward-shifted and more concentrated belief process $p_t$ in \eqref{wondpt} than (true) filtering (e.g., TD(est) yields $\mathbb{E}[p_t]=0.8363$ and $\mathrm{Var}(p_t)=0.006665$, versus (true) $\mathbb{E}[p_t]=0.8800$ and $\mathrm{Var}(p_t)=0.011605$). The effective discount rate $\widehat\delta_t$ in \eqref{hatmu11} increases as $p_t$ decreases, this shift moves $\widehat\delta_t$ closer to $\delta_2=0.3$, strengthening discounting and thereby lowering both the level and (especially) the variability of the discounted dividend stream and hence of $V$.}}}
it yields a slightly smaller mean but a markedly lower variance.
Moreover, in the (true) setting, the SNR and both annualized Sharpe-type ratios closely match the optimal benchmark, suggesting that any residual policy-approximation error has only a minor impact on risk-adjusted performance. In the (est) setting, the higher SNR reflects a more concentrated cross-path distribution of $V$; Sharpe-SR increases slightly under CTD but decreases marginally under ML, indicating modest method dependence of risk-adjusted surplus dynamics; whereas Sharpe-RI is essentially invariant across methods, pointing to stable discounted reward under the stronger effective discounting.
In addition, although $\mathrm{CTD}^*(0)$ exhibits noticeable deviations in the learned policy and value function in Fig.\ref{figgamma}, its out-of-sample performance remains competitive, with a still-high mean and low variance of $V$. This highlights the PO-RSEOD algorithm's effectiveness in capturing the optimal policy with reasonable accuracy, positioning it as a viable alternative to mathematically optimized solutions.

\appendix
\section{Proofs and more}\label{Section 5}
In what follows, we provide proofs of the main results for completeness.

\subsection{Proof of Proposition \ref{prop21}}
\begin{proof}\label{prop21pf}
(i) For $x> y\geq 0$ and $\pi_{y,p}\in \Pi_{y,p}[0,\mathbf{a}]$, to show that $V$ is strictly increasing in $x$, we consider an admissible policy $\pi_{x,p}\in\Pi_{x,p}[0,\mathbf{a}]$ such that
$\pi_{x,p}(u| \widehat{X}_t^{\pi}, p_t):
= \pi_{y,p}(u| \widehat{X}_t^{\pi}, p_t) \mathbf{1}_{\{t<\widehat{T}_{y,p}^{\pi}\}}
+\tfrac{e^{\frac{u}{\lambda}}}{\lambda (e^{\frac{\mathbf{a}}{\lambda}}-1)} \mathbf{1}_{\{t\geq \widehat{T}_{y,p}^{\pi}\}}$ with $(\widehat{X}_t^{\pi}, p_t)\in \mathbb{R}_+ \times (0,1)$ and $t\in[0, \widehat{T}_{x,p}^{\pi}]$. Since $e^x-1>x$, one can verify that $J(x, p; \pi_{x,p})> J(y, p; \pi_{y,p})$ for $\mathbf{a}>1$, and then by the arbitrariness of $\pi_{y,p}$ we get $V(x,p)\geq J(x, p; \pi_{x,p})> V(y,p)$.\\
(ii) The admissible policy $\hat{\pi}\equiv \tfrac{1}{\mathbf{a}}$ satisfies $J(x, p; {\hat{\pi}}) \geq (\tfrac{\mathbf{a}}{2}+\lambda \ln{\mathbf{a}}) \tfrac{1- e^{-(\delta_1 \vee \delta_2) \widehat{T}_{x,p}^{\pi}}}{\delta_1 \vee \delta_2}  \geq 0$ for $\mathbf{a}>1$, then we have $V(x,p)\geq 0$. From $\int_0^{\mathbf{a}} u \pi(u)\mathrm{d} u \leq \mathbf{a}$ and $-\int_0^{\mathbf{a}} \pi(u) \ln \pi(u) \mathrm{d} u \leq \ln \mathbf{a}$ (by the well-known Kullback-Leibler divergence property or Jensen's inequality), one has $\widehat{\mathcal{H}}_{\lambda}^{\pi}(\widehat{X}_t^{\pi}, p_t)\leq \mathbf{a} +\lambda\ln\mathbf{a}$, and thus $V(x, p)\leq \tfrac{1}{\delta_1 \wedge \delta_2} (\mathbf{a}+\lambda \ln \mathbf{a})$.
\end{proof}
\subsection{Proof of Proposition \ref{propfuncf00}}
\begin{proof}\label{propfuncf0011}
We first show that $f_{\lambda}(y)\in\mathbb{C}^{\infty}(\mathbb{R})$. For $y\neq 1$, \eqref{funcfyp} implies $f_{\lambda}(y)\in\mathbb{C}^{\infty}(\mathbb{R}\setminus \{1\})$, as its logarithmic argument is an entire function. By L'H\^{o}pital's rule, one has $\lim_{y\to 1} f_{\lambda}(y) = \lambda\ln\mathbf a$, and since $f_{\lambda}(1) = \lambda\ln\mathbf a$, we obtain $f_{\lambda}$ is continuous at $y=1$ (hence on $\mathbb{R}$).
Then, there exists $r>0$ such that $f_{\lambda}(y)\in(\lambda \ln\tfrac{\mathbf{a}}{2}, \lambda \ln\tfrac{3\mathbf{a}}{2})$ for all $|y-1|<r$, and $f_{\lambda}(y)\in\mathbb{C}^{\infty}\big((1-r, 1+r)\big)$. Combining the $\mathbb{C}^{\infty}$ regularity of $f_{\lambda}$ on $\mathbb{R}\setminus\{1\}$ with its smoothness in a neighborhood of $y=1$ yields $f_{\lambda}(y) \in \mathbb{C}^{\infty}(\mathbb{R})$.

We now prove (i). Let $z=\tfrac{1-y}{\lambda}$ with $y\in\mathbb{R}$, and rewrite $\ln \tfrac{e^{z \mathbf{a}} -1}{z} = \ln (\int_0^{\mathbf{a}} e^{z r}\mathrm{d} r)$. For any $z_1, z_2\in\mathbb{R}$ and $\zeta\in(0,1)$, since $\tfrac{1}{\zeta}, \tfrac{1}{1-\zeta}>1$, the H\"{o}lder inequality states that
\begin{align*}
&(\int_0^{\mathbf{a}} |e^{\zeta z_1 r}|^{\tfrac{1}{\zeta}} \mathrm{d} r)^{\zeta}
 \cdot (\int_0^{\mathbf{a}} |e^{(1-\zeta) z_2 r}|^{\tfrac{1}{1-\zeta}}\mathrm{d} r)^{(1- \zeta)}
\leq \int_0^{\mathbf{a}} |e^{\zeta z_1 r+(1-\zeta)z_2 r}| \mathrm{d} r,
\end{align*}
and it is strict whenever $z_1 \neq z_2$.
Then, for any $y_1, y_2\in \mathbb{R}$ with $y_1 \neq y_2$, we have
\begin{align*}
&\zeta f_{\lambda}(y_1) + (1-\zeta) f_{\lambda}(y_2)
= \zeta \lambda \ln (\int_0^{\mathbf{a}} e^{\frac{1-y_1}{\lambda} r}\mathrm{d} r)
+ (1- \zeta)\lambda \ln (\int_0^{\mathbf{a}} e^{\frac{1-y_2}{\lambda} r}\mathrm{d} r)\\
&=\lambda \ln \big((\int_0^{\mathbf{a}} e^{\frac{1-y_1}{\lambda} r}\mathrm{d} r)^{\zeta}
\cdot (\int_0^{\mathbf{a}} e^{\frac{1-y_2}{\lambda} r}\mathrm{d} r)^{(1- \zeta)}\big)
> \lambda \ln (\int_0^{\mathbf{a}} e^{\zeta \tfrac{1-y_1}{\lambda} r+(1-\zeta) \tfrac{1-y_2}{\lambda} r} \mathrm{d} r)\\
& =f_{\lambda}\big(\zeta y_1 + (1-\zeta) y_2\big).
\end{align*}
This shows that $f_{\lambda}$ is strictly convex, equivalently, $f_{\lambda}''(y)> 0$, and $f_{\lambda}'(y)$ is strictly increasing for all $y\in\mathbb{R}$.
We then examine the monotonicity of $f_{\lambda}$.
By L'H\^{o}pital's rule, we have $\lim_{y\rightarrow -\infty} f_{\lambda}'(y) = \lim_{z\rightarrow \infty} - (\tfrac{\mathbf{a} z -1}{z} +\tfrac{\mathbf{a}}{e^{z \mathbf{a}} -1}) = -\mathbf{a}$, and $\lim_{y\rightarrow \infty} f_{\lambda}'(y) = \lim_{z\rightarrow -\infty} - (\tfrac{\mathbf{a} z -1}{z} +\tfrac{\mathbf{a}}{e^{z \mathbf{a}} -1})=0$. Then, since $f_{\lambda}'$ is strictly increasing, we have $f_{\lambda}'$ is bounded with $f_{\lambda}'\in(-\mathbf{a}, 0)$. Consequently, $f_{\lambda}$ is strictly decreasing on $\mathbb{R}$. Moreover, using Taylor expansion of $e^{\mathbf{a} z}$ about $z$, \eqref{firdefy} is equivalent to
\begin{align*}
f_{\lambda}'(y)&= -\big(\tfrac{\mathbf{a} z - 1}{z} + \tfrac{\mathbf{a}}{e^{\mathbf{a} z}-1}\big)
=  -\big(\mathbf{a} - \tfrac{1}{z}\big)
   - \tfrac{1}{z}\big(1 - \tfrac{1}{2} \mathbf{a} z + \tfrac{1}{12} \mathbf{a}^2 z^2 + O(z^3)\big)\\
   &
   = -\big(\mathbf{a} - \tfrac{1}{z}\big)
    -\big(\tfrac{1}{z} - \tfrac{\mathbf{a}}{2} + \tfrac{\mathbf{a}^2}{12} z + O(z^2)\big)
   = -\tfrac{\mathbf{a}}{2}-\tfrac{\mathbf{a}^2}{12} z + O(z^2).
\end{align*}
Thus, $\lim_{y\to 1{\pm}} f_{\lambda}'(y)= \lim_{z\to 0{\pm}} \big(-\tfrac{\mathbf{a}}{2}-\tfrac{\mathbf{a}^2}{12} z + O(z^2)\big)= -\tfrac{\mathbf{a}}{2}$.

We proceed to prove (ii). Since (i) implies $|f_{\lambda}'(y)| \leq \mathbf{a}$ for all $y\in\mathbb{R}$, the mean value theorem yields, for any $y,z\in\mathbb{R}$, $|f_\lambda(y)-f_\lambda(z)|\le \mathbf a |y-z|$, so $f_\lambda$ is uniformly Lipschitz continuity on $\mathbb{R}$ with Lipschitz constant at most $\mathbf{a}$. Moreover, since $f_\lambda(1)=\lambda\ln \mathbf{a}$ is finite (for fixed $\lambda$) and $|f_\lambda(y)|\le |f_\lambda(1)| + \mathbf{a} |y-1|$, we have $f_\lambda$ is locally bounded.

The conclusion in (iii) follows directly from \eqref{funcfyp}.

As for (iv), since $f_{\lambda}$ is continuous and strictly decreasing on $\mathbb{R}$ with $\lim_{y\to\infty} f_{\lambda}(y)<0$ and $\lim_{y\to -\infty} f_{\lambda}(y)>0$, there exists a unique $y_0\in\mathbb{R}$ such that $f_{\lambda}(y_0)=0$.
For $\mathbf{a}>1$, since $f_{\lambda}(1)= \lambda \ln \mathbf{a}>0$ and $f_{\lambda}(1+\lambda)= \lambda \ln (1-e^{-\mathbf{a}})<0$, we get $y_0\in(1, 1+\lambda)$.

With regard to (v), since $f_{\lambda}'(y)> -\mathbf{a}$, we derive, for $y<y_0$ and for $y\geq y_0$, the following estimates
\begin{align*}
&|f_{\lambda}(y)| = f_{\lambda}(y) = f_{\lambda}(y_0) -\int_y^{y_0} f_{\lambda}'(z) \mathrm{d} z = -\int_y^{y_0} f_{\lambda}'(z) \mathrm{d} z <\int_y^{y_0} \mathbf{a} \mathrm{d} z = \mathbf{a} (y_0 -y), \\
&f_{\lambda}(y)| = -f_{\lambda}(y) = -f_{\lambda}(y_0) -\int_{y_0}^y f_{\lambda}'(z) \mathrm{d} z = -\int_{y_0}^y f_{\lambda}'(z) \mathrm{d} z <\int_{y_0}^y \mathbf{a} \mathrm{d} z = \mathbf{a} (y-y_0).
\end{align*}
Then, we have $|f_{\lambda}(y)|<\mathbf{a}|y-y_0|\leq \mathbf{a}|y| + \mathbf{a}|y_0|$ for any $y\in\mathbb{R}$.

Finally, we prove (vi). Let $s=1-y$. If $y<1$, then $s>0$, and as $\lambda\downarrow 0$,
\begin{align*}
&f_\lambda(y) = \lambda\big(\ln\lambda-\ln s+\ln(e^{\frac{\mathbf{a}s}{\lambda}}-1)\big)
= \mathbf{a}s+\lambda\big(\ln\lambda-\ln s+\ln(1-\mathrm{e}^{-\frac{\mathbf{a}s}{\lambda}})\big)\to \mathbf{a}s.
\end{align*}
If $y>1$, then $f_\lambda(y) = \lambda\big(\ln\lambda-\ln (-s)
+ \ln(1-\mathrm{e}^{\frac{\mathbf{a}s}{\lambda}})\big)\to 0$ as $\lambda\downarrow 0$. Further, $f_\lambda(1)=\lambda\ln \mathbf{a}\to 0$ as $\lambda\downarrow 0$. Thus $f_\lambda(y)\to f_0(y)=\mathbf{a}(1-y)_+$ pointwise on $\mathbb{R}$ as $\lambda\downarrow 0$.
\end{proof}

\subsection{Proof of Theorem \ref{verthem}}
\begin{proof}\label{verpf}
We first show that the distributional control $\pi^*=\big(\pi^*(\cdot|\widehat{X}_t^{\pi^*}, p_t)\big)_{t\geq 0}$ defined by \eqref{densityop} is admissible by verifying conditions (i)-(iii) in Definition \ref{def1}. Condition (i) holds because the control actions are $\mathbb{F}^{X}$-adapted and $\pi^*$ is normalized by construction.
For condition (ii), under $\pi^*$ the state dynamics \eqref{sde5} take the form \eqref{sde56}, with bounded measurable drift and constant (nondegenerate) volatility, \eqref{sde56} admits a unique strong solution.
Condition (iii) follows directly from the uniform bound $\widehat{\mathcal{H}}_{\lambda}^{\pi^*}(\widehat{X}_t^{\pi^*}, p_t) \leq \mathbf{a} +\lambda\ln\mathbf{a}$, proved in Proposition \ref{prop21}(ii). Thus, $\pi^*$ is admissible.

Consider an arbitrary admissible distributional control $\pi$. Define the infinitesimal generator $\mathcal{G}^{\pi}$ of the process $(\widehat{X}_t^{\pi}, p_t)$, acting on $v(x,p) \in \mathbb{C}^2((0, \infty)\times(0,1))$, as
\begin{align}\label{def:mathcalgpi}
\mathcal{G}^{\pi} v(x,p)
&:= \int_0^{\mathbf{a}} -u v_x(x,p) \pi(u|x, p) \mathrm{d} u + \tfrac{\sigma^2}{2} v_{xx}(x,p)
+\big(\mu_2+ (\mu_1-\mu_2)p\big) v_x(x,p)\nonumber
\\
& + \tfrac{1}{2} \tfrac{(\mu_1 -\mu_2)^2}{\sigma^2} p^2 (1-p)^2 v_{pp}(x,p)
+ \big(q_{21} -(q_{12}+q_{21}) p\big) v_p(x,p) \nonumber \\
& + (\mu_1 - \mu_2) p (1-p) v_{xp}(x,p).
\end{align}
Define stopping time $\tau_n$ as $\tau_n:=\inf\{t\geq 0: \widehat{X}_t^{\pi} \notin [0, n] \ \mathrm{or} \ p_t\in(0, 1/n]\cup [1-1/n, 1) \}$.
Using It\^{o}'s formula to $e^{- \Lambda_{t\wedge \tau_n}} v(\widehat{X}_{t\wedge \tau_n}^{\pi}, p_{t\wedge \tau_n})$ with $\Lambda_{t\wedge \tau_n} :=\int_0^{t\wedge\tau_n} \widehat{\delta}_s \mathrm{d} s$, we obtain
\begin{align}\label{elamvxpitoexp}
e^{-\Lambda_{t\wedge \tau_n} } v(\widehat{X}_{t\wedge \tau_n}^{\pi}, p_{t\wedge \tau_n})
&=v(x, p) + \int_0^{t\wedge \tau_n} e^{-\Lambda_s} (\mathcal{G}^{\pi} - \widehat{\delta}_s) v(\widehat{X}_s^{\pi}, p_s) \mathrm{d} s\\
&\ \
+\int_0^{t\wedge \tau_n} e^{-\Lambda_s} v_x \sigma \mathrm{d} \widehat{W}_s
+\int_0^{t\wedge \tau_n} e^{-\Lambda_s} v_p \tfrac{\mu_1 - \mu_2}{\sigma} p_s (1-p_s) \mathrm{d} \widehat{W}_s, \nonumber
\end{align}
with $v_x = v_x(\widehat{X}_s^{\pi}, p_s)$ and $v_p = v_p(\widehat{X}_{s}^{\pi}, p_s)$.
The two stochastic integrals driven by Brownian motion are martingales because $|v_x|$, $|v_p|$ are bounded and the integrand processes are square integrable on every finite horizon. Thus, after taking expectations on both sides, we obtain
\begin{align*}
&\mathbb{E}[e^{-\Lambda_{t\wedge \tau_n} } v(\widehat{X}_{t\wedge \tau_n}^{\pi}, p_{t\wedge \tau_n})]
= v(x, p) + \mathbb{E}_{x,p}\big[\int_0^{t\wedge \tau_n} e^{-\Lambda_s}
(\mathcal{G}^{\pi} -\widehat{\delta}_s) v(\widehat{X}_s^{\pi}, p_s) \mathrm{d} s\big].
\end{align*}
Using the exploratory HJB equation \eqref{hjbeq}, we further obtain the following inequality
\begin{align*}
v(x, p) & \geq  \mathbb{E}[e^{-\Lambda_{t\wedge \tau_n}} v(\widehat{X}_{t\wedge \tau_n}^{\pi}, p_{t\wedge \tau_n})] \\
& \ \ \
+ \mathbb{E}_{x,p}\big[\int_0^{t\wedge \tau_n} e^{-\Lambda_s} \big(\int_0^{\mathbf{a}} \big(u-\lambda \ln\pi(u|\widehat{X}_{s}^{\pi}, p_{s})\big) \pi(u|\widehat{X}_{s}^{\pi}, p_{s}) \mathrm{d} u \big) \mathrm{d} s\big].
\end{align*}
Note that $\tau_n\rightarrow \widehat{T}^{\pi}_{x,p}$ for $\mathbb{P}$-a.s. as $n$ goes to infinity. Letting $t, n\uparrow \infty$, since $|v|$ is bounded and $v(0, \cdot)=0$, by the bounded convergence theorem, we have
\begin{align*}
&\lim_{t, n\uparrow \infty} \mathbb{E}[e^{-\Lambda_{t\wedge \tau_n}} v(\widehat{X}_{t\wedge \tau_n}^{\pi}, p_{t\wedge \tau_n})]
= \mathbb{E}[e^{-\Lambda_{\widehat{T}^{\pi}_{x,p}}} v(\widehat{X}_{\widehat{T}^{\pi}_{x,p}}^{\pi}, p_{\widehat{T}^{\pi}_{x,p}})]=0.
\end{align*}
Furthermore, since $\pi$ is admissible and satisfies condition (iii) of Definition \ref{def1}, we get
\begin{align*}
&\mathbb{E}_{x,p}\big[\int_0^{t\wedge \tau_n} e^{-\Lambda_s} \big(\int_0^{\mathbf{a}} \big(u-\lambda \ln\pi(u|\widehat{X}_{s}^{\pi}, p_{s})\big) \pi(u|\widehat{X}_{s}^{\pi}, p_{s}) \mathrm{d} u \big) \mathrm{d} s\big] \\
&
= \mathbb{E}_{x,p}\big[\int_0^{t\wedge \tau_n} e^{-\Lambda_s} \int_0^{\mathbf{a}} u \pi(u|\widehat{X}_{s}^{\pi}, p_{s}) \mathrm{d} u \mathrm{d} s\big] \\
&
-\lambda \mathbb{E}_{x,p}\big[\int_0^{t\wedge \tau_n} e^{-\Lambda_s} \int_0^{\mathbf{a}} \big(\ln\pi(u|\widehat{X}_{s}^{\pi}, p_{s}) \pi(u|\widehat{X}_{s}^{\pi}, p_{s}) +1\big) \mathrm{d} u \mathrm{d} s\big]
+ \lambda \mathbb{E}_{x,p}\big[\int_0^{t\wedge \tau_n} e^{-\Lambda_s} \mathbf{a} \mathrm{d} s \big].
\end{align*}
Together with the non-negativity of $\int_0^{\mathbf{a}} u \pi(u|\widehat{X}_{s}^{\pi}, p_{s}) \mathrm{d} u$ and $y\ln y +1\geq y >0$ for any $y\in(0,\infty)$, allows us to use the monotone convergence theorem to obtain
\begin{align}\label{vvalueineq}
&v(x, p) \geq \mathbb{E}_{x,p}\big[\int_0^{\widehat{T}^{\pi}_{x,p}} e^{-\Lambda_s} \big(\int_0^{\mathbf{a}} \big(u-\lambda \ln\pi(u|\widehat{X}_{s}^{\pi}, p_{s})\big) \pi(u|\widehat{X}_{s}^{\pi}, p_{s}) \mathrm{d} u \big) \mathrm{d} s\big]=J(x, p; \pi).
\end{align}
As $\pi$ is chosen arbitrarily, we obtain $v(x, p)\geq V(x, p)$ for all $(x,p)\in\mathbb{R}_+ \times (0,1)$.
In addition, \eqref{vvalueineq} becomes an equality if the supremum in \eqref{hjbeq} is achieved, i.e., $\pi=\pi^*$. Thus, $v(x,p)$ is the value function and $\pi^*$ is the optimal distributional control.
\end{proof}

\subsection{Proof of Proposition \ref{vbound2}}
\begin{proof}\label{vbound22}
Substituting $\pi^*$ given by \eqref{densityop} into \eqref{enth1}, we obtain
\begin{align}\label{h526}
&\widehat{\mathcal{H}}_{\lambda}^{\pi^*}(\widehat{X}_t^{\pi^*}, p_t)
= \int_0^{\mathbf{a}}\big(u - \lambda \ln\pi^*(u| \widehat{X}_t^{\pi^*},p_t)\big) \pi^*(u| \widehat{X}_t^{\pi^*}, p_t)\mathrm{d}u \\
&= \big[
v_x \big(\mathbf{a} - \tfrac{\lambda}{1- v_x}
+ \tfrac{\mathbf{a}}{e^{\frac{\mathbf{a}}{\lambda}(1-v_x)} -1}\big)
- \lambda \ln \big(\tfrac{1- v_x}{\lambda(e^{\frac{\mathbf{a}}{\lambda}(1-v_x)} -1)}\big) \big]\mathbf{1}_{\{v_x\neq 1\}}
+ \big[\tfrac{\mathbf{a}}{2} + \lambda \ln \mathbf{a}\big] \mathbf{1}_{\{v_x= 1 \}},\nonumber
\end{align}
with $v_x = v_x(\widehat{X}_t^{\pi^*}, p_t)$. By defining
\begin{align*}
\mathcal{K}_{\lambda,\mathbf{a}}(y)
& := \big[\big(\mathbf{a} - \tfrac{1}{y}
+ \tfrac{\mathbf{a}}{e^{\mathbf{a}y} -1}\big)\mathbf{1}_{\{y\neq 0\}}+ \big(\tfrac{\mathbf{a}}{2} \big) \mathbf{1}_{\{y = 0 \}} \big] (1-\lambda y) \\
& \ \ \ +\lambda \ln \big[\big(\tfrac{e^{\mathbf{a}y} -1}{y}\big) \mathbf{1}_{\{y\neq 0\}}
+ \mathbf{a} \mathbf{1}_{\{y = 0 \}}\big],
\end{align*}
the above expression simplifies to $\widehat{\mathcal{H}}_{\lambda}^{\pi^*}(\widehat{X}_t^{\pi^*}, p_t) = \mathcal{K}_{\lambda,\mathbf{a}}\big(\tfrac{1-v_x}{\lambda} \big)$.
By \eqref{opdipro}, we have
\begin{align}\label{vmax45}
&v(x,p)
=\mathbb{E}_{x,p}\big[\int_0^{\widehat{T}_{x,p}^{\pi^*}} e^{-\int_0^{t} \widehat{\delta}_s \mathrm{d}s} \mathcal{K}_{\lambda,\mathbf{a}}\big(\tfrac{1-v_x}{\lambda} \big)\mathrm{d}t\big].
\end{align}
We first determine the maximum of $\mathcal{K}_{\lambda,\mathbf{a}}(y)$. Taking the first-order derivative of $\mathcal{K}_{\lambda,\mathbf{a}}(y)$ w.r.t. $y (\neq 0)$ yields,
\begin{align}\label{mathkpn}
&\mathcal{K}_{\lambda,\mathbf{a}}'(y)
= \tfrac{(1-\lambda y)K_1(y)}{y^2 (e^{\mathbf{a}y}-1)^2},
\end{align}
where $K_1(y) :=e^{2\mathbf{a}y} - (2+\mathbf{a}^2 y^2) e^{\mathbf{a}y} +1$. Differentiating $K_1(y)$ w.r.t. $y$, we obtain $K_1'(y)= \mathbf{a} e^{\mathbf{a}y} K_2(y)$, where $K_2(y):=2 e^{\mathbf{a}y} - \mathbf{a}^2 y^2 - 2\mathbf{a}y -2$. Further differentiation yields $K_2'(y)= 2 \mathbf{a} K_3(y)$, where $K_3(y):=e^{\mathbf{a}y} - \mathbf{a} y - 1$ and $K_3'(y)= \mathbf{a} (e^{\mathbf{a}y} - 1)$.
From $K_3''(y) = \mathbf{a}^2 e^{\mathbf{a}y}>0$ we infer that $K_3'(y)$ is strictly increasing w.r.t. $y$. Together with $K_3'(0)=0$, this yields $K_3'(y)<0$ for $y<0$ and $K_3'(y)>0$ for $y>0$, and further $K_3(y)$ is decreasing for $y<0$ and increasing for $y>0$, with $K_3(y)> K_3(0)=0$. Then, $K_2'(y)>0$, implying that $K_2(y)$ is increasing. Combining this with $K_2(0)=0$, we obtain $K_2(y)<0$ for $y<0$ and $K_2(y)>0$ for $y>0$. As a result, $K_1'(y)<0$ for $y<0$ and $K_1'(y)>0$ for $y>0$, which means $K_1(y)$ is decreasing for $y<0$ and increasing for $y>0$, and $K_1(y)>K_1(0)=0$ for $y\neq 0$.
Accordingly, by \eqref{mathkpn}, the sign of $\mathcal{K}_{\lambda,\mathbf{a}}'(y)$ depends solely on that of $1-\lambda y$. Specifically, $\mathcal{K}_{\lambda,\mathbf{a}}'(y)>0$ when $1-\lambda y>0$ (i.e. $y<1/\lambda$), and $\mathcal{K}_{\lambda,\mathbf{a}}'(y)<0$ when $1-\lambda y<0$ (i.e. $y>1/\lambda$). Therefore, $\mathcal{K}_{\lambda,\mathbf{a}}(y)$ is increasing for $y<1/\lambda$ and decreasing for $y>1/\lambda$, and the maximum of $\mathcal{K}_{\lambda,\mathbf{a}}(y)$ is obtained at $y=1/\lambda$, $\max \mathcal{K}_{\lambda,\mathbf{a}}(y) =  \mathcal{K}_{\lambda,\mathbf{a}}(\tfrac{1}{\lambda})
=\lambda \ln \big(\lambda(e^{\frac{\mathbf{a}}{\lambda}} -1)\big)$.
Note that $\lambda \ln \big(\lambda(e^{\ln(\frac{1}{\lambda} +1)} -1)\big)=0$, and $\mathcal{K}_{\lambda,\mathbf{a}}(\frac{1}{\lambda}) =0$ when $\mathbf{a}=\lambda\ln(\tfrac{1}{\lambda}+1)$. If $\mathbf{a}<\lambda\ln(\tfrac{1}{\lambda}+1)$, then $\mathcal{K}_{\lambda,\mathbf{a}}(\tfrac{1}{\lambda}) <0$, whereas if $\mathbf{a}>\lambda\ln(\tfrac{1}{\lambda}+1)$, then $\mathcal{K}_{\lambda,\mathbf{a}}(\frac{1}{\lambda}) >0$.

We now examine the limiting behavior of $v(x, p)$ as $x\rightarrow \infty$. From
$\mathcal{K}_{\lambda,\mathbf{a}}(y)\leq \mathcal{K}_{\lambda,\mathbf{a}}(1/\lambda)=\lambda \ln \big(\lambda(e^{\mathbf{a}/\lambda} -1)\big)$ and \eqref{vmax45}, it follows that
\begin{align*}
&v(x,p)
\leq\mathbb{E}_{x,p}\big[\int_0^{\widehat{T}_{x,p}^{\pi^*}} e^{-\int_0^{t} \widehat{\delta}_s \mathrm{d}s} \lambda \ln \big(\lambda(e^{\frac{\mathbf{a}}{\lambda}} -1)\big)\mathrm{d}t\big].
\end{align*}
As $x\rightarrow \infty$, the dominated convergence theorem yields
\begin{align}\label{asyv11}
\lim_{x\rightarrow \infty} v(x,p)
&\leq
\lambda \ln \big(\lambda(e^{\frac{\mathbf{a}}{\lambda}} -1)\big) \cdot \mathbb{E}_{x,p}\big[\int_0^{\infty} e^{-\int_0^{t} \widehat{\delta}_s \mathrm{d}s} \mathrm{d}t\big].
\end{align}
Further, under a policy $\hat{\pi}=\{\hat{\pi}(u|\widehat{X}_t^{\hat{\pi}}, p_t)\}_{t\geq 0}$ with $\hat{\pi}(u|\widehat{X}_t^{\hat{\pi}}, p_t)
:= \tfrac{e^{{u}/{\lambda}}}{\lambda (e^{\mathbf{a}/{\lambda}}-1)}$, we get
\begin{align*}
v(x,p)\geq J(x,p; \hat{\pi})
&= \lambda \ln \big(\lambda(e^{\frac{\mathbf{a}}{\lambda}} -1)\big) \cdot \mathbb{E}_{x,p}\big[\int_0^{\widehat{T}_{x,p}^{\hat{\pi}}} e^{-\int_0^{t} \widehat{\delta}_s \mathrm{d}s} \mathrm{d} t\big].
\end{align*}
Letting $x\rightarrow \infty$, the dominated convergence theorem yields
\begin{align}\label{asyv22}
\lim_{x\rightarrow \infty} v(x,p)
&\geq
\lambda \ln \big(\lambda(e^{\frac{\mathbf{a}}{\lambda}} -1)\big) \cdot \mathbb{E}_{x,p}\big[\int_0^{\infty} e^{-\int_0^{t} \widehat{\delta}_s \mathrm{d}s} \mathrm{d}t\big].
\end{align}
Both \eqref{asyv11} and \eqref{asyv22} imply the asymptotic behavior of $v(x,p)$ as $x\to \infty$,
\begin{align*}
\lim_{x\rightarrow \infty} v(x,p)
= f_{\lambda}(0) \cdot
 \mathbb{E}_{x,p}\big[\int_0^{\infty} e^{-\int_0^{t} \widehat{\delta}_s \mathrm{d}s} \mathrm{d}t\big].
\end{align*}
It reduces to \eqref{asyv00} when $\delta=\delta_1=\delta_2$, and Proposition \ref{propfuncf00} (vi) further implies that $\lim_{\substack{x\rightarrow\infty \\ \lambda\downarrow 0}} v(x,p) = \tfrac{\mathbf{a}}{\delta}$.
\end{proof}

\subsection{Proof of Proposition \ref{theosoluodse}}
\begin{proof}\label{theosoluodse00}
Define $\mathcal{L} g_i:= \tfrac12 A(p) g_i'' + B(p) g_i' - C(p)g_i$, $i=1,2$. \\
\text{(i)} Since $A, A', B, C\in \mathbb{C}((0,1))$ with $A(p), C(p)>0$, the existence of a smooth solution to \eqref{eq:V_const_eq00} follows directly from Section 8 of \cite{feller1952}, which also ensures $g_i(p)>0$ whenever $\hat{\varpi}_i(p) f_{\lambda}(0)\geq 0$ (equivalently, $f_{\lambda}(0)\geq 0$ because $\hat{\varpi}_i(p)> 0$ on $(0,1)$). \\
\text{(ii)}
Equation \eqref{eq:V_const_eq00} constitutes the Dirichlet problem: $\mathcal{L}g_i=-\hat{\varpi}_i(p) f_{\lambda}(0)$ on $(0,1)$, with $g_i(0)$, $g_i(1)$ specified by \eqref{addicondi}. Since $A, B, C \in \mathbb{C}^2((0,1))$ and $C(p)>0$, Theorem 6.25 of \cite{gilbarg1977} yields existence and uniqueness of a solution $g_i\in\mathbb{C}^2((0, 1))\cap \mathbb{C}([0, 1])$.

To bound $g_i$, we apply the comparison principle (Theorem 3.3 in \cite{gilbarg1977}): if $g_i, h_i\in \mathbb{C}^2((0, 1))\cap \mathbb{C}([0, 1])$ satisfy $\mathcal{L} g_i\geq \mathcal{L} h_i$ on $(0,1)$ and $g_i(0)\leq h_i(0)$, $g_i(1)\leq h_i(1)$, then $g_i\leq h_i$ on $(0,1)$. Denote by $C_*:= \inf_{p\in(0,1)} C(p)= \delta_1 \wedge \delta_2 > 0$.
Then, for $f_{\lambda}(0)\geq 0$, since $\mathcal{L} (g_i -\tfrac{f_{\lambda}(0)}{C_*})= -\hat{\varpi}_i(p) f_{\lambda}(0) + C(p) \tfrac{f_{\lambda}(0)}{C_*} \geq 0$ on $(0,1)$ and $g_i(0), g_i(1)\leq \tfrac{f_{\lambda}(0)}{C_*}$, by the comparison principle, we obtain $g_i - \tfrac{f_{\lambda}(0)}{C_*}\leq 0$, i.e. $g_i\leq \tfrac{f_{\lambda}(0)}{C_*}$.
Similarly, as $\mathcal{L} g_i= -\hat{\varpi}_i(p)f_{\lambda}(0) \leq 0$ and $g_i(0), g_i(1)\geq 0$, we get $g_i\geq 0$. Hence, $0\leq g_i(p)\leq \tfrac{f_{\lambda}(0)}{C_*}$ for $f_{\lambda}(0)\geq 0$.
For $f_{\lambda}(0)<0$, the same argument gives $\tfrac{f_{\lambda}(0)}{C_*}\leq g_i(p)\leq 0$. \\
\text{(iii)} To bound $g_i'(p)$, we differentiate both sides of \eqref{eq:V_const_eq00} w.r.t. $p$ and obtain
\begin{align*}
\tfrac12 A(p) \tfrac{\mathrm{d}}{\mathrm{d}p^2} \big(g_i'(p)\big) + B_1(p) \tfrac{\mathrm{d}}{\mathrm{d}p}\big(g_i'(p)\big) + C_1(p) g_i'(p) + D_1(p) g_i(p) + E_i = 0,
\end{align*}
with $g_i'(0)=g_i'(1)=0$, where $A$ is as in \eqref{eq:V_const_eq00} and $B_1, C_1, D_1$ are bounded (indeed polynomial) functions on $(0,1)$, and $C_1(p)<0$.
Applying the same comparison principle as for $g_i$ to the equation for $g_i'$ yields bounds for $g_i'(p)$ on $(0,1)$.
An analogous argument applied once more gives bounds for $g_i''(p)$.
\end{proof}

\subsection{Adjoint-weight in Proposition \ref{propkappa}}
\begin{proof}\label{prop:adjoint-weight}
Recall $\mathcal{L} g_i(p) := \tfrac12 A(p) g_i''(p)+B(p) g_i'(p)-C(p) g_i(p)$ for $p\in(0,1)$, with $A, B, C$ as in \eqref{apbpcpdp}. Let $w\in \mathbb{C}^2((0,1))$. The following integration-by-parts identity holds:
\begin{align}\label{eq:ibp-main}
&\textstyle\int_0^1 (\mathcal L g_i)(p) w(p) \mathrm{d} p
=\big[\tfrac12 A(p) w(p)\,g_i'(p)\big]_{0}^{1} \\ &
\ \ \  -\textstyle\int_0^1 \tfrac{1}{2} \big(A(p)w(p) \big)' \, g_i'(p) \mathrm{d} p
+\textstyle\int_0^1 B(p)w(p)\, g_i'(p) \mathrm{d}p
- \textstyle\int_0^1 C(p) w(p) \,g_i(p)\,\mathrm d p \nonumber\\
&=\big[\tfrac12 A(p) w(p)\,g_i'(p)\big]_{0}^{1}
-\big[\big(\tfrac12 \big(A(p) w(p)\big)'-B(p) w(p)\big)\,g_i(p)\big]_{0}^{1} \nonumber\\
&\ \ \ +\textstyle\int_0^1 \big(\tfrac12 \big(A(p) w(p)\big)''-\big(B(p) w(p)\big)'-C(p) w(p)\big)\,g_i(p)\,\mathrm d p. \nonumber
\end{align}
The formal adjoint of $\mathcal L$ w.r.t. the Lebesgue inner product is $\mathcal L^* w =\tfrac12 (A w)''-(B w)'-C w$.
A canonical weight $w$ that eliminates the terms involving $g_i'$ in \eqref{eq:ibp-main} satisfies $(Aw)'=2 B w$, i.e. $\tfrac{w'}{w}=\tfrac{2B-A'}{A}$. Its unique positive solution is, for any fixed $p_0\in(0,1)$,
\begin{equation}\label{eq:w-explicit}
w(p)=\tfrac{K}{A(p)}
\exp\big(\textstyle\int_{p_0}^{p}\tfrac{2B(s)}{A(s)} \mathrm{d}s\big),
\end{equation}
where the constant $K>0$ is determined by the normalization $\int_0^1 w(p)\mathrm{d} p =1$.

Since
$\textstyle\int_{p_0}^{p}\tfrac{2B(s)}{A(s)} \mathrm{d}s
=\frac{2}{\beta_0}\big(
(q_{21}-q_{12})(\ln\tfrac{p}{1-p} - \ln\tfrac{p_0}{1-p_0})
-q_{21}(\tfrac{1}{p}-\tfrac{1}{p_0})
-q_{12} (\tfrac{1}{1-p}-\tfrac{1}{1-p_0}) \big)$ with $\beta_0 := \tfrac{(\mu_1 -\mu_2)^2}{\sigma^2}$,
substituting it into \eqref{eq:w-explicit} gives, for $\beta_1 = \tfrac{2}{\beta_0}(q_{21}-q_{12})$,
\begin{align}\label{eq:w-closed-2}
&w(p)
= \tfrac{\hat{K}}{\beta_0}
p^{\beta_1-2}(1-p)^{-\beta_1 -2}
\exp\big\{-\tfrac{2}{\beta_0}\big(\tfrac{q_{21}}{p}+\tfrac{q_{12}}{1-p}\big)
\big\},
\end{align}
where $\hat{K}:=K \exp\big\{-\tfrac{2}{\beta_0}\big(
(q_{21}-q_{12})\ln\tfrac{p_0}{1-p_0}-\tfrac{q_{21}}{p_0}-\tfrac{q_{12}}{1-p_0}
\big)\big\}$.
For such $w$, one has $\tfrac12 (A w)''-(B w)'=0$ and $\tfrac12 (A w)'-B w=0$, then \eqref{eq:ibp-main} simplifies to
\begin{equation}\label{eq:ibp-final}
\textstyle\int_0^1 (\mathcal L g_i)(p) w(p) \mathrm{d}p
=-\textstyle\int_0^1 C(p) g_i(p) w(p) \mathrm{d}p.
\end{equation}
Note that $w$ depends only on $A$ and $B$ (hence not on $g_i$ and $\hat{\varpi}_i$).
By \eqref{eq:V_const_eq00}, $g_i$ satisfies $\mathcal L g_i(p)+\hat{\varpi}_i(p) f_{\lambda}(0)=0$, integrating against the weight \eqref{eq:w-closed-2} and using \eqref{eq:ibp-final} yield
\begin{align*}
\textstyle\int_0^1 C(p) g_i(p) w(p) \mathrm{d}p
= f_{\lambda}(0)\textstyle\int_0^1 \hat{\varpi}_i(p) w(p)\,\mathrm{d}p,
\end{align*}
which is used repeatedly to eliminate $F_{0,i}$ in the projection formula \eqref{kappasolu}.
\end{proof}

\subsection{Proof of Theorem \ref{verthem22}}
\begin{proof}\label{verthem2200}
Proposition \ref{theosoluodse} guarantees a unique solution $g_i(p)\in\mathbb{C}^2((0, 1))\cap \mathbb{C}([0, 1])$ $(i=1,2)$ to \eqref{eq:V_const_eq00}, subject to the boundary values \eqref{boucong1}-\eqref{addicondi}, and provides uniform bounds on $g_i$ and $g_i'$. Proposition \ref{propkappa} further yields constants $\kappa_i>0$ under the condition \eqref{kapgeqsuff}.
These results specify a candidate value function $v$ of the form \eqref{eq:single_exp_ansatz}. Then $v\in \mathbb{C}^2((0, \infty)\times(0,1))$ solves \eqref{hjb156}, satisfies $v(0,\cdot)=0$, and has uniformly bounded $|v|$, $|v_x|$ and $|v_p|$ on $\mathbb{R}_+ \times (0,1)$. By the verification theorem (Theorem \ref{verthem}), $v$ is the value function of \eqref{opdipro}, and the optimal distributional control is given by \eqref{densityop}.
\end{proof}

\subsection{Proof of Lemma \ref{lemdencov}}
\begin{proof}\label{lemdencovproof}
We first show that $\pi$ is strongly admissible by verifying conditions (i)-(ii) in Definition \ref{def510}. Let $\theta(x, p):=c_x(x, p)$ and set $M:=\|c\|_{\infty}+\|c_x\|_{\infty}+\|c_p\|_{\infty}<\infty$.\\
(i) Since $G(u, \theta)$ is continuous and positive on $[0, \mathbf a] \times [-M, M]$, there exist constants $l_{\pm}$ that depend only on $M$ and $\mathbf{a}$, defined by
$l_-:=\min_{u\in[0, \mathbf{a}], |\theta|\leq M} G(u, \theta)>0$ and $l_+:= \max_{u\in[0, \mathbf{a}], |\theta|\leq M} G(u, \theta)<\infty$,
such that $l_- \leq \pi(u|x,p)\leq l_+$ for all $u\in[0,\mathbf a]$.   \\
(ii) From $c\in\mathbb{C}^2((0, \infty) \times (0,1))$ we have $c_{xx}$ and $c_{xp}$ are continuous, and then bounded on any compact set $K\Subset (0, \infty)\times (0, 1)$: $\|c_{xx}\|_{\mathbb{L}^\infty(K)}<\infty$ and $\|c_{xp}\|_{\mathbb{L}^\infty(K)}<\infty$. Moreover, since $\partial_{\theta} G$ is continuous on $[0, \mathbf{a}]\times [-M, M]$, we have $C_{G}:=\sup_{u\in[0, \mathbf{a}], |\theta|\leq M} |\partial_{\theta} G(u, \theta)|<\infty$. By the chain rule, $\partial_x \pi(u| x, p) = \partial_{\theta} G(u, \theta(x,p))c_{xx}(x, p)$ and $\partial_p \pi(u| x, p) = \partial_{\theta} G(u, \theta(x, p))c_{xp}(x, p)$, and then
\begin{align*}
&\sup_{u\in[0, \mathbf{a}]} \big(
  \|\partial_x \pi(u \mid \cdot,\cdot)\|_{\mathbb{L}^\infty(K)}
  + \|\partial_p \pi(u \mid \cdot,\cdot)\|_{\mathbb{L}^\infty(K)}\big)
\leq C_{G}(\|c_{xx}\|_{\mathbb{L}^\infty(K)} + \|c_{xp}\|_{\mathbb{L}^\infty(K)})
< \infty.
\end{align*}
Thereby, there exists $L_{K}:=C_{G}\big(\|c_{xx}\|_{\mathbb{L}^\infty(K)} + \|c_{xp}\|_{\mathbb{L}^\infty(K)}\big)<\infty$ such that,
\begin{align*}
& |\pi(u| x_1, p_1) - \pi(u| x_2, p_2)| \leq L_K(|x_1 -x_2| + |p_1 -p_2|), \ \ \ \mathrm{for}\ (x_i,p_i)\in K, u\in[0,\mathbf a].
\end{align*}

Fix a strongly admissible policy $\pi$, we now prove that $J(\cdot,\cdot;\pi)$ is the unique bounded viscosity solution to PDE \eqref{defjxpde} with $J(0, \cdot;\pi)=0$. Under such $\pi$, $b^{\pi}$ and $\widehat{\mathcal{H}}_{\lambda}^{\pi}$ are bounded on $\mathbb{R}_+\times (0,1)$ and locally Lipschitz in $(x, p)$, so all coefficients of \eqref{defjxpde} are bounded and locally H\"{o}lder. Since $C(p)\ge \delta_1\wedge\delta_2>0$ and $\widehat{\mathcal H}^{\pi}_{\lambda}$ is bounded, the Feynman-Kac representation yields the uniform bound $\|J\|_{\infty}\leq \tfrac{ \|\widehat{\mathcal{H}}_{\lambda}^{\pi}\|_{\infty}}{\delta_1 \wedge \delta_2} < \infty$.
To establish existence and uniqueness in the viscosity sense, we employ a uniformly elliptic regularization: for $\epsilon>0$, introduce $\mathcal G^{\pi,\epsilon}
:=\mathcal G^{\pi}+\tfrac{\epsilon}{2}(\partial_{xx}+\partial_{pp})$, equivalently, replace the (positive-semidefinite, possibly degenerate) diffusion matrix $\mathbf A(p):=\bigl(\begin{smallmatrix}
\sigma^{2} & (\mu_{1}-\mu_{2})p(1-p)\\
(\mu_{1}-\mu_{2})p(1-p) & A(p)
\end{smallmatrix}\bigr)$ by
$\mathbf A_\epsilon(p):=\mathbf A(p)+\epsilon \mathbf I_2$ with $\mathbf I_2$ the $2\times 2$ identity matrix, so that
$\xi^\top\mathbf A_\epsilon(p)\xi\ge \epsilon|\xi|^2$ for any $\xi\in\mathbb R^2$.
Let $D_{R,\varepsilon}:=(0,R)\times(\varepsilon,1-\varepsilon)$ for $R>1$ and $\varepsilon\in(0,\tfrac12)$ and consider the Dirichlet problem:
\begin{equation}\label{eq:Dirichlet-Re}
\big(\mathcal G^{\pi,\epsilon}-C(p)\big)y +\widehat{\mathcal H}_{\lambda}^{\pi} =0
\ \mathrm{in}\  D_{R,\varepsilon}, \qquad
 y=0 \ \mathrm{on}\ \partial D_{R,\varepsilon}.
\end{equation}
By classical Dirichlet theory, \eqref{eq:Dirichlet-Re} admits a unique solution $y\in \mathbb{C}^{2,\eta}(D_{R,\varepsilon})\cap \mathbb{C}(\overline{D_{R,\varepsilon}})$ for some $\eta\in(0,1)$, where $\mathbb{C}^{2,\eta}$ denotes the standard H\"{o}lder space.
The maximum principle gives the $\mathbb{L}^\infty$ bound $\|y\|_{\mathbb{L}^\infty(D_{R,\varepsilon})}
\le \frac{\|\widehat{\mathcal H}_{\lambda}^{\pi}\|_{\infty}}{\delta_1 \wedge \delta_2}$. For any compact set $K\Subset(0,\infty)\times(0,1)$, choose $R$ large and $\varepsilon$ small such that
$K\Subset D_{R,\varepsilon}$. Interior Schauder estimates (see \cite{gilbarg1977}) then yield
\[
\|y\|_{\mathbb{C}^{2,\eta}(K)}
\le C_{K,\varepsilon}\big(\|y\|_{\mathbb{L}^\infty(D_{R,\varepsilon})}
+\|\widehat{\mathcal H}_{\lambda}^{\pi}\|_{\mathbb{C}^{0,\eta}(D_{R,\varepsilon})}\big),
\]
with $C_{K,\varepsilon}$ independent of $R$. Letting $R\to\infty$ and using a standard diagonal argument, one obtains a bounded function
$J^\epsilon\in \mathbb{C}^{2,\eta}_{\mathrm{loc}}((0,\infty)\times(0,1))$ solving the global regularized equation $(\mathcal G^{\pi,\epsilon}-C(p))\,J^\epsilon+\widehat{\mathcal H}_\lambda^\pi=0$ on $(0,\infty)\times(0,1)$ with $J^\epsilon(0,p)=0$, and $\|J^\epsilon\|_\infty\le \tfrac{\|\widehat{\mathcal H}_\lambda^\pi\|_\infty}{\delta_1 \wedge \delta_2}$.
Then, along a sequence $\epsilon_n\downarrow 0$, the stability of bounded viscosity solutions implies $J^{\epsilon_n}\to J(\cdot,\cdot;\pi)$ locally uniformly, where
$J(\cdot,\cdot;\pi)$ is the unique bounded viscosity solution to \eqref{defjxpde} with $J(0,\cdot;\pi)=0$.

To further establish the regularity of $J$, we invoke hypoelliptic regularity of $\mathcal G^\pi$. The second-order part of \eqref{defjxpde} can be written as
\[
\tfrac{1}{2}\big(\sigma^2\partial_{xx}+2(\mu_1-\mu_2)p(1-p)\partial_{xp}+A(p)\partial_{pp}
\big)
=\tfrac12\,V_1^2
\
\mathrm{with} \
V_1:=\sigma\partial_x+\tfrac{\mu_1-\mu_2}{\sigma}p(1-p)\partial_p,
\]
and let $V_0$ denote the first-order drift field in \eqref{defjxpde}.
Hence, $(\mathcal G^\pi-C(p))y=\frac12 V_1^2 y + V_0 y - C(p)y$.
Moreover, the H\"ormander bracket condition holds on $(0,\infty)\times(0,1)$: the vector fields $V_1$ and $[V_0,V_1]$ span $\mathbb R^2$ at every point, so $\mathcal G^\pi$ is hypoelliptic. By H\"ormander's hypoelliptic regularity theorem, $J(\cdot,\cdot;\pi)$ is smooth in the interior, in particular, $J(\cdot,\cdot;\pi)\in \mathbb{C}^{2}((0,\infty)\times(0,1))$.

Finally, we derive uniform bounds for the first partial derivatives of $J$ via the killed semigroup representation $J=\int_0^\infty P_t\widehat{\mathcal H}^{\pi}_{\lambda}\mathrm{d}t$ with $\|P_t \widehat{\mathcal H}^{\pi}_{\lambda}\|_\infty\le e^{-(\delta_1 \wedge \delta_2) t}\|\widehat{\mathcal H}^{\pi}_{\lambda}\|_\infty$,
where $(P_t)_{t\ge0}$ is the Markov semigroup of $(\widehat X^\pi_t,p_t)$ killed at rate $C(p)$.
For H\"ormander diffusions with bounded Lipschitz coefficients, standard Bismut-Elworthy-Li/Malliavin gradient estimates (see \cite{elworthy1994}) yield
$\|\nabla P_t \widehat{\mathcal H}^{\pi}_{\lambda}\|_\infty \le r\,t^{-1/2}e^{-(\delta_1 \wedge \delta_2) t}\|\widehat{\mathcal H}^{\pi}_{\lambda}\|_\infty$ for bounded $\widehat{\mathcal H}^{\pi}_{\lambda}$,
with $r$ depending on the bounds of $b^\pi$.
Integrating over $t$ gives
\begin{align*}
&\|\nabla J(\cdot,\cdot;\pi)\|_\infty
\le \int_0^\infty \|\nabla P_t \widehat{\mathcal H}_\lambda^\pi\|_\infty\,\mathrm{d}t
\le r\,\|\widehat{\mathcal H}_\lambda^\pi\|_\infty\int_0^\infty t^{-1/2}e^{-(\delta_1 \wedge \delta_2) t}\,\mathrm{d}t
<\infty,
\end{align*}
i.e. $\|J_x\|_\infty+\|J_p\|_\infty <\infty$, which implies $|J_x|$ and $|J_p|$ are bounded on $(0,\infty)\times(0,1)$, with bounds depending on those of $b^{\pi}$ and $\widehat{\mathcal{H}}_{\lambda}^{\pi}$.
\end{proof}

\subsection{Proof of Theorem \ref{policyimpro}}
\begin{proof}\label{policyimpro00}
If $\pi$ is strongly admissible, then by Lemma \ref{lemdencov} the associated objective function $J(x, p; \pi) \in \mathbb{C}^2((0,\infty) \times (0,1))$ with $|J|$, $|J_x|$ and $|J_p|$ bounded on $(0,\infty) \times (0,1)$, further, the updated policy $\tilde{\pi}$ is itself strongly admissible.
Since $J(x, p; \pi)$ is a solution to the PDE \eqref{defjxpde},
in view of \eqref{hjbeq} and \eqref{densityop}, $\tilde{\pi}$ is the maximizer of $\sup_{\bar{\pi}\in{\Pi_{x,p}[0, \mathbf{a}]}} \widehat{\mathcal{H}}_{\lambda}^{\bar{\pi}}(x, p)- b^{\bar{\pi}}(x, p) J_x(x, p; \pi)$. Then, we have
\begin{align}\label{hatpigg}
&\big(\mathcal{G}^{\tilde{\pi}} - \big(\delta_2 +(\delta_1-\delta_2)p\big)\big)\,J(x,p; \pi) + \widehat{\mathcal{H}}_{\lambda}^{\tilde{\pi}}(x, p) \geq 0,
\end{align}
with $\mathcal{G}^{\tilde{\pi}}$ as in \eqref{def:mathcalgpi}. Define stopping time $\tilde{\tau}_n:=\inf\{t\geq 0: \widehat{X}_t^{\tilde{\pi}} \notin [0, n] \ \mathrm{or} \ p_t\in(0, 1/n]\cup [1-1/n, 1) \}$.
Applying It\^{o}'s formula to $e^{- \Lambda_{t\wedge \tilde{\tau}_n}} J(\widehat{X}_{t\wedge \tilde{\tau}_n}^{\tilde{\pi}}, p_{t\wedge \tilde{\tau}_n}; \pi)$ with $\Lambda_{t\wedge \tilde{\tau}_n} := \int_0^{t\wedge \tilde{\tau}_n} \widehat{\delta}_s \mathrm{d} s$, we deduce from \eqref{hatpigg} that
\begin{align*}
& e^{-\Lambda_{t\wedge \tilde{\tau}_n} } J(\widehat{X}_{t\wedge \tilde{\tau}_n}^{\tilde{\pi}}, p_{t\wedge \tilde{\tau}_n}; \pi)
=J(x, p; \pi) + \int_0^{t\wedge \tilde{\tau}_n} e^{-\Lambda_s} (\mathcal{G}^{\tilde{\pi}} - \widehat{\delta}_s) J(\widehat{X}_s^{\tilde{\pi}}, p_s; \pi) \mathrm{d} s\\
&\ \
+\int_0^{t\wedge \tilde{\tau}_n} e^{-\Lambda_s} J_x \sigma \mathrm{d} \widehat{W}_s
+\int_0^{t\wedge \tilde{\tau}_n} e^{-\Lambda_s} J_p \tfrac{\mu_1 - \mu_2}{\sigma} p_s (1-p_s) \mathrm{d} \widehat{W}_s\\
&\geq J(x, p; \pi) - \int_0^{t\wedge \tilde{\tau}_n} e^{-\Lambda_s} \widehat{\mathcal{H}}_{\lambda}^{\tilde{\pi}}(\widehat{X}_s^{\tilde{\pi}}, p_s)\mathrm{d} s
+\int_0^{t\wedge \tilde{\tau}_n} e^{-\Lambda_s} J_x \sigma \mathrm{d} \widehat{W}_s \\
 & \ \
+\int_0^{t\wedge \tilde{\tau}_n} e^{-\Lambda_s} J_p \tfrac{\mu_1 - \mu_2}{\sigma} p_s (1-p_s) \mathrm{d} \widehat{W}_s,
\end{align*}
with $J_x = J_x(\widehat{X}_s^{\tilde{\pi}}, p_s; \pi)$ and $J_p = J_p(\widehat{X}_{s}^{\tilde{\pi}}, p_s; \pi)$. The two stochastic integrals driven by Brownian motion are martingales because $|J_x|$ and $|J_p|$ are bounded and the integrand processes are square integrable on every finite horizon. Note that $\tilde{\tau}_n\rightarrow \widehat{T}^{\tilde{\pi}}_{x,p}$ for $\mathbb{P}$-a.s. as $n$ goes to infinity. Analogously to the derivation from \eqref{elamvxpitoexp} to \eqref{vvalueineq}, taking expectations on both sides and letting $t,n\to\infty$, the bounded/dominated convergence theorem yields $J(x, p; \tilde{\pi}) = \mathbb{E}_{x,p}\big[\int_0^{\widehat{T}^{\tilde{\pi}}_{x,p}} e^{-\Lambda_s} \widehat{\mathcal{H}}_{\lambda}^{\tilde{\pi}}(\widehat{X}_s^{\tilde{\pi}}, p_s)\mathrm{d} s \big] \geq J(x, p; \pi)$.
\end{proof}

\subsection{Proof of Proposition \ref{prop52a}}
\begin{proof}\label{prop52a00}
Recall that $v^{\pi}(0, \cdot) = \tilde{v}^{\pi}(0, \cdot) = 0$ and $\widehat{X}_{\widehat{T}^{\pi}_{x,p}}^{\pi} = 0$. By \eqref{mtpi12} and the definition of $\tilde{M}^{\pi}$ we have, $\tilde{M}_{\widehat{T}^{\pi}_{x,p}}^{\pi} = \int_0^{\widehat{T}^{\pi}_{x,p}} e^{-\Lambda_s} \widehat{\mathcal{H}}_{\lambda}^{\pi}(\widehat{X}_s^{\pi}, p_s) \mathrm{d} s
= M_{\widehat{T}^{\pi}_{x,p}}^{\pi}$.
Then, since $\widehat{\mathcal{H}}_{\lambda}^{\pi}$, $v^{\pi}$ and $\tilde{v}^{\pi}$ are bounded, both $M^{\pi}$ and $\tilde{M}^{\pi}$ are uniformly integrable $\mathbb{F}$-martingales, by optional sampling theorem, we have, for all $(x, p)\in \mathbb{R}_+ \times (0,1)$,
\begin{align*}
&\tilde{v}^{\pi}(x, p) = \tilde{M}_0^{\pi} = \mathbb{E}[ \tilde{M}_{\widehat{T}^{\pi}_{x,p}}^{\pi} | \mathcal{F}_0]
= \mathbb{E}[ M_{\widehat{T}^{\pi}_{x,p}}^{\pi} | \mathcal{F}_0]
= M_0^{\pi}
= v^{\pi}(x, p).
\end{align*}
\end{proof}

\subsection{FD scheme in Section \ref{subsec610}}\label{finitediff00}
A FD scheme (see e.g., \cite{leveque2007}) discretizes the interval $[0, 1]$ into $M+1$ grid points $0=p_0<p_1<\cdots<p_{M}=1$ with uniform step $\Delta p = 1/ M$ and $p_j = j \Delta p$. For interior nodes $j=1, \cdots, M-1$, we set $g_i''(p_j) \approx \tfrac{g_i(p_{j+1}) - 2 g_i(p_j) + g_i(p_{j-1})}{(\Delta p)^2}$.
For $g_i'(p) $, we use the centered difference when convection is mild, otherwise an upwind formula (with direction set by $\operatorname{sign}(B)$) is adopted to preserve monotonicity. Define the local P\'{e}clet number $\mathrm{Pe}_j:=\tfrac{2\,|B(p_j)|\,\Delta p}{A(p_j)}$ and $\mathrm{Pe}_\ast=2$, we set $g_i'(p_j)\approx\tfrac{g_i(p_{j+1}) - g_i(p_{j-1})}{2\Delta p}$ if $\mathrm{Pe}_j \le \mathrm{Pe}_\ast$; $g_i'(p_j)\approx\tfrac{g_{i}(p_j)-g_i(p_{j-1})}{\Delta p}$ if $\mathrm{Pe}_j>\mathrm{Pe}_\ast$ and $B(p_j)\ge 0$; $g_i'(p_j)\approx \tfrac{g_i(p_{j+1})-g_i(p_{j})}{\Delta p}$ if $\mathrm{Pe}_j>\mathrm{Pe}_\ast$ and $B(p_j)< 0$. Substituting into \eqref{eq:V_const_eq00} gives the FD equation at interior nodes,
\begin{align*}
\hat{A}_j g_i(p_{j-1}) + \hat{B}_j g_i(p_j) + \hat{C}_j g_i(p_{j+1})
= -\hat\varpi_i(p_j) f_{\lambda}(0),
\end{align*}
where $\hat{A}_j:= \tfrac{A(p_j)}{2(\Delta p)^2} - \tfrac{B(p_j)}{2 \Delta p} \mathbf{1}_{\{\mathrm{Pe}_j \le \mathrm{Pe}_\ast\}} - \tfrac{B(p_j)}{\Delta p}\mathbf{1}_{\{\mathrm{Pe}_j>\mathrm{Pe}_\ast, B(p_j)\geq 0\}}$, $\hat{B}_j:= - \tfrac{A(p_j)}{(\Delta p)^2} -C(p)+\tfrac{B(p_j)}{\Delta p}\mathbf{1}_{\{\mathrm{Pe}_j>\mathrm{Pe}_\ast\}}(\mathbf{1}_{\{B(p_j)\geq 0\}}- \mathbf{1}_{\{B(p_j)<0\}})$ and $\hat{C}_j:=\tfrac{A(p_j)}{2 (\Delta p)^2} + \tfrac{B(p_j)}{2\Delta p}\mathbf{1}_{\{\mathrm{Pe}_j \le \mathrm{Pe}_\ast\}} + \tfrac{B(p_j)}{\Delta p}\mathbf{1}_{\{\mathrm{Pe}_j>\mathrm{Pe}_\ast, B(p_j)< 0\}}$.
We solve $(g_i(p_0), \ldots, g_i(p_M))$ from the following linear system
\[
{\footnotesize
\begin{pmatrix}
1 & 0 & 0 & 0 & \cdots & 0 & 0\\
\hat{A}_{1} & \hat{B}_{1} & \hat{C}_{1} & 0 & \cdots & 0 & 0\\
0 & \hat{A}_{2} & \hat{B}_{2} & \hat{C}_{2} & \cdots & 0 & 0\\
0 & 0 & \hat{A}_{3} & \hat{B}_{3} & \ddots & \vdots & \vdots\\
\vdots & \vdots & \ddots & \ddots & \ddots & \hat{C}_{M-2} & 0\\
0 & 0 & \cdots & 0 & \hat{A}_{M-1} & \hat{B}_{M-1} & \hat{C}_{M-1}\\
0 & 0 & 0 & 0 & \cdots & 0 & 1
\end{pmatrix}
\begin{pmatrix}
g_{i}(p_0)\\
g_{i}(p_1)\\
g_{i}(p_2)\\
g_{i}(p_3)\\
\vdots\\
g_{i}(p_{M-1})\\
g_{i}(p_M)
\end{pmatrix}
=
\begin{pmatrix}
\varpi_i^{(0)} g(0)\\
-\hat\varpi_i(p_1) f_{\lambda}(0)\\
-\hat\varpi_i(p_2) f_{\lambda}(0)\\
-\hat\varpi_i(p_3) f_{\lambda}(0)\\
\vdots\\
-\hat\varpi_i(p_{M-1}) f_{\lambda}(0)\\
\varpi_i^{(1)} g(1)
\end{pmatrix},
}
\]
with $g(0)$ and $g(1)$ given by \eqref{boucong1}-\eqref{boucong2}.
We first apply the FD scheme to solve $g_i(p)$ $(i=1,2)$, and evaluate $w_j:= w(p_j)$ from \eqref{eq:w-closed-2} on the same grid (normalization is optional, as scaling $w$ does not affect $\kappa_i$).
We then determine $\kappa_i$ by evaluating the coefficients $F_{k, i}$ $(k=0, 1, 2)$ in \eqref{kappasolu} using the composite trapezoidal rule, namely, with samples $y_j:= y(p_j)$, the approximation $\int_0^1 y(p)\mathrm{d}p \approx \Delta t \big(\frac{1}{2} y_0 + y_1 + \cdots + y_{M-1} + \tfrac{1}{2} y_M\big)$. In our implementation, we set $\Delta p = 10^{-4}$ to obtain accurate results.

\section*{Funding}

This research is supported by National Natural Science Foundation of China (Grant No.11671204), and China Scholarship Council (No.202206840089).

\section*{Competing interests}
There were no competing interests to declare which arose during the preparation or publication process of this article.

\renewcommand{\bibsection}{\section*{References}}



\end{document}